\newcount\mgnf\newcount\tipi\newcount\tipoformule\newcount\greco

\tipi=2          
\tipoformule=0   


\global\newcount\numsec
\global\newcount\numfor
\global\newcount\numtheo
\global\advance\numtheo by 1

\def\senondefinito#1{\expandafter\ifx\csname#1\endcsname\relax}

\def\SIA #1,#2,#3 {\senondefinito{#1#2}%
\expandafter\xdef\csname #1#2\endcsname{#3}\else
\write16{???? ma #1,#2 e' gia' stato definito !!!!} \fi}

\def\etichetta(#1){(\veroparagrafo.\veraformula)%
\SIA e,#1,(\veroparagrafo.\veraformula) %
\global\advance\numfor by 1%
\write15{\string\FU (#1){\equ(#1)}}%
\write16{ EQ #1 ==> \equ(#1) }}

\def\letichetta(#1){\veroparagrafo.\verotheo
\SIA e,#1,{\veroparagrafo.\verotheo}
\global\advance\numtheo by 1
\write15{\string\FU (#1){\equ(#1)}}
\write16{ Sta \equ(#1) == #1 }}

\def\tetichetta(#1){\veroparagrafo.\veraformula 
\SIA e,#1,{(\veroparagrafo.\veraformula)}
\global\advance\numfor by 1
\write15{\string\FU (#1){\equ(#1)}}
\write16{ tag #1 ==> \equ(#1)}}

\def\FU(#1)#2{\SIA fu,#1,#2 }\def\etichettaa(#1){(A\veroparagrafo.\veraformula)%
\SIA e,#1,(A\veroparagrafo.\veraformula) %
\global\advance\numfor by 1%
\write15{\string\FU (#1){\equ(#1)}}%
\write16{ EQ #1 ==> \equ(#1) }}

\def\BOZZA{
\def\alato(##1){%
 {\rlap{\kern-\hsize\kern-1.4truecm{$\scriptstyle##1$}}}}%
\def\aolado(##1){%
 {
{
 \rlap{\kern-1.4truecm{$\scriptstyle##1$}}}}}
 }

\def\alato(#1){}
\def\aolado(#1){}

\def\veroparagrafo{\number\numsec}
\def\veraformula{\number\numfor}
\def\verotheo{\number\numtheo}

\def\Eq(#1){\eqno{\etichetta(#1)\alato(#1)}}
\def\eq(#1){\etichetta(#1)\alato(#1)}
\def\teq(#1){\tag{\aolado(#1)\tetichetta(#1)\alato(#1)}}
\def\Eqa(#1){\eqno{\etichettaa(#1)\alato(#1)}}
\def\eqa(#1){\etichettaa(#1)\alato(#1)}
\def\eqv(#1){\senondefinito{fu#1}$\clubsuit$#1
\write16{#1 non e' (ancora) definito}%
\else\csname fu#1\endcsname\fi}
\def\equ(#1){\senondefinito{e#1}\eqv(#1)\else\csname e#1\endcsname\fi}

\def\Lemma(#1){\aolado(#1)Lemma \letichetta(#1)}%
\def\Theorem(#1){{\aolado(#1)Theorem \letichetta(#1)}}%
\def\Proposition(#1){\aolado(#1){Proposition \letichetta(#1)}}%
\def\Corollary(#1){{\aolado(#1)Corollary \letichetta(#1)}}%
\def\Remark(#1){{\noindent\aolado(#1){\bf Remark \letichetta(#1).}}}%
\def\Definition(#1){{\noindent\aolado(#1){\bf Definition 
\letichetta(#1)$\!\!$\hskip-1.6truemm}}}
\def\Example(#1){\aolado(#1) Example \letichetta(#1)$\!\!$\hskip-1.6truemm}

\def\include#1{
\openin13=#1.aux \ifeof13 \relax \else
\input #1.aux \closein13 \fi}

\openin14=\jobname.aux \ifeof14 \relax \else
\input \jobname.aux \closein14 \fi
\openout15=\jobname.aux


{\count255=\time\divide\count255 by 60 \xdef\hourmin{\number\count255}
        \multiply\count255 by-60\advance\count255 by\time
   \xdef\hourmin{\hourmin:\ifnum\count255<10 0\fi\the\count255}}

\def\oramin{\hourmin }

\def\data{\number\day/\ifcase\month\or january \or february \or march \or april
\or may \or june \or july \or august \or september
\or october \or november \or december \fi/\number\year;\ \oramin}

\newcount\pgn \pgn=1
\def\foglio{\number\numsec:\number\pgn
\global\advance\pgn by 1}
\def\foglioa{A\number\numsec:\number\pgn
\global\advance\pgn by 1}


\def\TIPIO{
\font\setterm=amr7 
\def \settepunti{\def\rm{\fam0\setterm}
\textfont0=\setterm   
\normalbaselineskip=9pt\normalbaselines\rm }\let\nota=\settepunti}


\def\TIPITOT{
\font\twelverm=cmr12
\font\twelvei=cmmi12
\font\twelvesy=cmsy10 scaled\magstep1
\font\twelveex=cmex10 scaled\magstep1
\font\twelveit=cmti12
\font\twelvett=cmtt12
\font\twelvebf=cmbx12
\font\twelvesl=cmsl12
\font\ninerm=cmr9
\font\ninesy=cmsy9
\font\eightrm=cmr8
\font\eighti=cmmi8
\font\eightsy=cmsy8
\font\eightbf=cmbx8
\font\eighttt=cmtt8
\font\eightsl=cmsl8
\font\eightit=cmti8
\font\sixrm=cmr6
\font\sixbf=cmbx6
\font\sixi=cmmi6
\font\sixsy=cmsy6
\font\twelvetruecmr=cmr10 scaled\magstep1
\font\twelvetruecmsy=cmsy10 scaled\magstep1
\font\tentruecmr=cmr10
\font\tentruecmsy=cmsy10
\font\eighttruecmr=cmr8
\font\eighttruecmsy=cmsy8
\font\seventruecmr=cmr7
\font\seventruecmsy=cmsy7
\font\sixtruecmr=cmr6
\font\sixtruecmsy=cmsy6
\font\fivetruecmr=cmr5
\font\fivetruecmsy=cmsy5
\textfont\truecmr=\tentruecmr
\scriptfont\truecmr=\seventruecmr
\scriptscriptfont\truecmr=\fivetruecmr
\textfont\truecmsy=\tentruecmsy
\scriptfont\truecmsy=\seventruecmsy
\scriptscriptfont\truecmr=\fivetruecmr
\scriptscriptfont\truecmsy=\fivetruecmsy
\textfont\truecmr=\tentruecmr
\scriptfont\truecmr=\seventruecmr
\scriptscriptfont\truecmr=\fivetruecmr
\textfont\truecmsy=\tentruecmsy
\scriptfont\truecmsy=\seventruecmsy
\scriptscriptfont\truecmr=\fivetruecmr
\scriptscriptfont\truecmsy=\fivetruecmsy
\def \eightpoint{\def\rm{\fam0\eightrm}
\textfont0=\eightrm \scriptfont0=\sixrm \scriptscriptfont0=\fiverm
\textfont1=\eighti \scriptfont1=\sixi   \scriptscriptfont1=\fivei
\textfont2=\eightsy \scriptfont2=\sixsy   \scriptscriptfont2=\fivesy
\textfont3=\tenex \scriptfont3=\tenex   \scriptscriptfont3=\tenex
\textfont\itfam=\eightit  \def\it{\fam\itfam\eightit}%
\textfont\slfam=\eightsl  \def\sl{\fam\slfam\eightsl}%
\textfont\ttfam=\eighttt  \def\tt{\fam\ttfam\eighttt}%
\textfont\bffam=\eightbf  \scriptfont\bffam=\sixbf
\scriptscriptfont\bffam=\fivebf  \def\bf{\fam\bffam\eightbf}%
\tt \ttglue=.5em plus.25em minus.15em
\setbox\strutbox=\hbox{\vrule height7pt depth2pt width0pt}%
\normalbaselineskip=9pt
\let\sc=\sixrm  \let\big=\eightbig  \normalbaselines\rm
\textfont\truecmr=\eighttruecmr
\scriptfont\truecmr=\sixtruecmr
\scriptscriptfont\truecmr=\fivetruecmr
\textfont\truecmsy=\eighttruecmsy
\scriptfont\truecmsy=\sixtruecmsy }\let\nota=\eightpoint}

\newfam\msbfam   
\newfam\truecmr  
\newfam\truecmsy 
\newskip\ttglue
\ifnum\tipi=0\TIPIO \else\ifnum\tipi=1 \TIPI\else \TIPITOT\fi\fi


\def\sqr#1#2{{\vcenter{\vbox{\hrule height.#2pt
     \hbox{\vrule width.#2pt height#1pt \kern#1pt
   \vrule width.#2pt}\hrule height.#2pt}}}}
\def\qed{\hfill $\mathchoice\sqr64\sqr64\sqr{2.1}3\sqr{1.5}3$}



\newcount\foot
\foot=1
\def\note#1{\footnote{${}^{\number\foot}$}{\ftn #1}\advance\foot by 1}
\def\tag #1{\eqno{\hbox{\rm(#1)}}}
\def\frac#1#2{{#1\over #2}}

\def\text#1{\quad{\hbox{#1}}\quad}

\def\proof{{ \noindent {\bf  Proof.} }}

\def\thanks{\noindent{\bf Aknowledgements: }}



\font\ch=cmbx12
\font\ftn=cmr8

\font\it=cmti10
\font\bf=cmbx10
\font\sm=cmr7

%
\catcode`\X=12\catcode`\@=11
\def\n@wcount{\alloc@0\count\countdef\insc@unt}
\def\n@wwrite{\alloc@7\write\chardef\sixt@@n}
\def\n@wread{\alloc@6\read\chardef\sixt@@n}
\def\crossrefs#1{\ifx\alltgs#1\let\tr@ce=\alltgs\else\def\tr@ce{#1,}\fi
   \n@wwrite\cit@tionsout\openout\cit@tionsout=\jobname.cit 
   \write\cit@tionsout{\tr@ce}\expandafter\setfl@gs\tr@ce,}
\def\setfl@gs#1,{\def\@{#1}\ifx\@\empty\let\next=\relax
   \else\let\next=\setfl@gs\expandafter\xdef
   \csname#1tr@cetrue\endcsname{}\fi\next}
\newcount\sectno\sectno=0\newcount\subsectno\subsectno=0\def\r@s@t{\relax}
\def\resetall{\global\advance\sectno by 1\subsectno=0
  \gdef\firstpart{\number\sectno}\r@s@t}
\def\resetsub{\global\advance\subsectno by 1
   \gdef\firstpart{\number\sectno.\number\subsectno}\r@s@t}
\def\v@idline{\par}\def\firstpart{\number\sectno}
\def\l@c@l#1X{\firstpart.#1}\def\gl@b@l#1X{#1}\def\t@d@l#1X{{}}
\def\m@ketag#1#2{\expandafter\n@wcount\csname#2tagno\endcsname
     \csname#2tagno\endcsname=0\let\tail=\alltgs\xdef\alltgs{\tail#2,}%
  \ifx#1\l@c@l\let\tail=\r@s@t\xdef\r@s@t{\csname#2tagno\endcsname=0\tail}\fi
   \expandafter\gdef\csname#2cite\endcsname##1{\expandafter
     \ifx\csname#2tag##1\endcsname\relax?\else{\rm\csname#2tag##1\endcsname}\fi
    \expandafter\ifx\csname#2tr@cetrue\endcsname\relax\else
     \write\cit@tionsout{#2tag ##1 cited on page \folio.}\fi}%
   \expandafter\gdef\csname#2page\endcsname##1{\expandafter
     \ifx\csname#2page##1\endcsname\relax?\else\csname#2page##1\endcsname\fi
     \expandafter\ifx\csname#2tr@cetrue\endcsname\relax\else
     \write\cit@tionsout{#2tag ##1 cited on page \folio.}\fi}%
   \expandafter\gdef\csname#2tag\endcsname##1{\global\advance
     \csname#2tagno\endcsname by 1%
   \expandafter\ifx\csname#2check##1\endcsname\relax\else%
\fi
   \expandafter\xdef\csname#2check##1\endcsname{}%
   \expandafter\xdef\csname#2tag##1\endcsname
     {#1\number\csname#2tagno\endcsnameX}%
   \write\t@gsout{#2tag ##1 assigned number \csname#2tag##1\endcsname\space
      on page \number\count0.}%
   \csname#2tag##1\endcsname}}%
\def\m@kecs #1tag #2 assigned number #3 on page #4.%
   {\expandafter\gdef\csname#1tag#2\endcsname{#3}
   \expandafter\gdef\csname#1page#2\endcsname{#4}}
\def\re@der{\ifeof\t@gsin\let\next=\relax\else
    \read\t@gsin to\t@gline\ifx\t@gline\v@idline\else
    \expandafter\m@kecs \t@gline\fi\let \next=\re@der\fi\next}
\def\t@gs#1{\def\alltgs{}\m@ketag#1e\m@ketag#1s\m@ketag\t@d@l p
    \m@ketag\gl@b@l r \n@wread\t@gsin\openin\t@gsin=\jobname.tgs \re@der
    \closein\t@gsin\n@wwrite\t@gsout\openout\t@gsout=\jobname.tgs }
\outer\def\localtags{\t@gs\l@c@l}
\outer\def\globaltags{\t@gs\gl@b@l}
\outer\def\newlocaltag#1{\m@ketag\l@c@l{#1}}
\outer\def\newglobaltag#1{\m@ketag\gl@b@l{#1}}

\def\t@gsoff#1,{\def\@{#1}\ifx\@\empty\let\next=\relax\else\let\next=\t@gsoff
   \expandafter\gdef\csname#1cite\endcsname{\relax}
   \expandafter\gdef\csname#1page\endcsname##1{?}
   \expandafter\gdef\csname#1tag\endcsname{\relax}\fi\next}
\def\verbatimtags{\let\ift@gs=\iffalse\ifx\alltgs\relax\else
   \expandafter\t@gsoff\alltgs,\fi}
\catcode`\X=11 \catcode`\@=\active
\localtags
%
\setbox200\hbox{$\scriptscriptstyle \data $}
\global\newcount\numpunt
\magnification=1000
\hoffset=0.cm
\baselineskip=14pt  
\parindent=12pt
\lineskip=4pt\lineskiplimit=0.1pt
\parskip=0.1pt plus1pt

\hyphenation{small}

\def\chap #1#2{\line{\ch #1\hfill}\numsec=#2\numfor=1}


\def\AA{{\cal A}}
\def\BB{{\cal B}}
\def\CC{{\cal C}}
\def\DD{{\cal D}}
\def\EE{{\cal E}}
\def\FF{{\cal F}}

\def\JJ{{\cal J}}

\def\LL{{\cal L}}

\def\PP{{\cal P}}
\def\SS{{\cal S}}

\def\MM{{\cal M}}

\def\VV{{\cal V}}

\def\LL{{\cal L}}

\def\ZZ{{\cal Z}}

\def\a{\alpha}
\def\b{\beta}
\def\d{\delta}
\def\e{\epsilon}
\def\ve{\varepsilon}

\def\g{\gamma}

\def\l{\lambda}
\def\r{\rho}
\def\s{\sigma}

\def\th{\theta}

\def\D{\Delta}
\def\L{\Lambda}
\def\G{\Gamma}
\def\O{\Omega}

\def\E{{I\kern-.25em{E}}}
\def\N{{I\kern-.25em{N}}}
\def\M{{I\kern-.25em{M}}}
\def\R{{I\kern-.25em{R}}}
\def\Z{{Z\kern-.425em{Z}}}
\def\1{{1\kern-.25em\hbox{\rm I}}}
\def\eu{{1\kern-.25em\hbox{\sm I}}}
\def\C{{I\kern-.64em{C}}}
\def\P{{I\kern-.25em{P}}}
\def\T{{\rm I\kern -4.3pt{\rm T}}}
\def\W{{\rm W\kern -9.5pt{\rm W}}}
\def\J{{\rm I\kern -4.3pt{\rm J}}}

\def\X{{\rm X\kern -6.0pt{\rm X}}}
\def\V{{\rm V\kern -6.5pt{\rm V}}}
\def\U{{\rm U\kern -6.0pt{\rm U}}}
\def\Y{{\rm Y\kern -6.0pt{\rm Y}}}
\def\IG{\ {\rm I\kern -5.8pt{\rm G}}}
\def\IIG{{I\kern-.5em{G}}}

\def\wt{\widetilde}

\def\leq{\le}
\def\geq{\ge}



 \def \mus {\overline}
 \def\HH{{\cal H}}
 \def\EE{{\cal E}}
  \def\BB{{\cal B}}
  \def \ac   {\AA \CC}

\hfuzz 16pt
\baselineskip 14pt

\catcode`\@=11  

 \def \mus {\overline}
 \def\HH{{\cal H}}
 \def\EE{{\cal E}}
  \def\BB{{\cal B}}
  \def \ac   {\AA \CC}

\hfuzz 16pt
\baselineskip 14pt

\catcode`\@=11  


\centerline {\bf  Quenched  large deviations for Glauber evolution 
   } \centerline {\bf with  Kac  interaction   and   random field.    \footnote{$^*$}  
{\eightrm work supported by ANR-07-BLAN-0230, ANR-2010-BLAN-0108, GDRE GREFI-MEFI,
   INDAM-CNRS,   Prin07:
20078XYHYS,  PHC Galil\'ee 19762TG,
 Universit\'a di Roma TRE,  Universit\'e  de Rouen, Universit\'e Paris Descartes.  }}
\vskip1cm 

\centerline{  
Olivier Benois,     \footnote{$^1$}      {\eightrm 
Universit\'e de Rouen, LMRS, UMR 6085, Avenue de l'Universit\'e, 
BP. 12,  76801 Saint Etienne du Rouvray, France.} 
\footnote{}      {\eightrm Olivier.Benois@univ-rouen.fr}
\hskip.2cm 
 Mustapha  Mourragui,     \footnote{$^2$}
 {\eightrm Universit\'e de Rouen, LMRS, UMR 6085, Avenue de l'Universit\'e, 
BP. 12,  76801 Saint Etienne du Rouvray, France.}
\footnote{}      {\eightrm Mustapha.Mourragui@univ-rouen.fr}
\hskip.2cm 
Enza Orlandi,   \footnote{$^3$}{\eightrm 
Dipartimento di Matematica, Universit\'a di Roma Tre, L.go S.Murialdo 1,
00146 Roma, Italy.}
\footnote{}      {\eightrm orlandi@mat.uniroma3.it}
\hskip.2cm 
Ellen Saada,    \footnote{$^4$}{\eightrm 
CNRS, MAP5, UMR 8145, Universit\'e Paris Descartes, Sorbonne Paris Cit\'e, 45 rue des Saints-P\`eres, 
75270 Paris cedex 06, France.}
\footnote{}      {\eightrm ellen.saada@mi.parisdescartes.fr}
\hskip.2cm 
Livio Triolo   \footnote{$^5$}{\eightrm 
Dipartimento di Matematica, Universit\'a di Roma ``Tor Vergata'', Via della Ricerca Scientifica,
00133 Roma, Italy.}
\footnote{}      {\eightrm triolo@mat.uniroma2.it}
}
\hskip.2cm 
\footnote{}{\eightrm {\eightit Key Words }: Interacting particle systems; 
random environment; Kac Interaction; large deviations.}

\footnote{}{\eightrm {\eightit 2000 Mathematics Subject Classification} :
60K35, 82C22, 60F10.}
\vskip.5cm

{ \bf Abstract} 

We study a spin-flip model with Kac type interaction, in the presence of a random field given by i.i.d. bounded random variables. The system, spatially inhomogeneous, evolves according to a non conservative (Glauber) dynamics. We show an almost sure (with respect to the random
field) large deviation principle for the empirical  magnetizations of  this   process.  The rate functional associated with the large deviation principle depends on the
statistical properties of the external random field, it is lower semicontinuous with compact level sets.

  \bigskip \bigskip  
  \chap {   1.  Introduction} 1
  \numsec= 1 \numfor= 1 \numtheo=1

We consider interacting spin-flip systems, in dimension $d$,  
with Kac type interaction in the presence of a random field 
given by i.i.d. bounded random variables.  
Kac potentials $J_\g$ are two-body interactions with range
$\g^{-1}$ and strength $\g^d$,  where $\g$ is a dimensionless scaling 
parameter. When $\g\to 0$, i.e. very long range compared with the 
inter particle distance, the strength of the interaction becomes very
weak, but in such a way that the total interaction between one particle 
and all the others is finite. Kac potentials were  introduced in [KUH], 
and then generalized in  [LP], to present a rigorous derivation  of the 
van der Waals theory of a gas-liquid phase transition.   
  There has been in the last decades an increasing interest in them. Indeed they induce the intermediate  scale of interaction $\g^{-1}$ (called mesoscopic) between the  microscopic (lattice) one and a macroscopic one much bigger than the latter. They are suitable to interpolate not only  between short  and long range interactions, but, scaling space and time as functions of $\g$, one can hope to obtain more insights into the physics of the model.  Recently they  have been considered as models to describe social interactions and more general complex social systems, see for example [CDS] and references therein. 

There has been several results on Kac Ising spin systems (without random field) in equilibrium and in non equilibrium statistical mechanics.   
We refer for a survey to the book [P]. {The papers  [C], [CE] were among the first dealing with dynamics issues. They considered} spin systems in a torus evolving according to a reversible and non conservative (Glauber) dynamics, with Kac interactions. In [C] the long time analysis of the spin system is studied,  using large deviations techniques. In [CE] the main results are the infinite particle limits of the non-critical  and critical fluctuation processes. In [DOPT] and [KS] (see also references therein), a complete description of the development and motion of interfaces (long time behaviour) has been derived: it is governed by the law of {\sl motion by mean curvature}.

A natural extension of this analysis is its application to disordered systems. One of the simplest prototype models is obtained by adding a random magnetic field to an interacting spin system. Equilibrium statistical properties of these systems have been widely studied
in the last  decades, see [Bo] for a survey of results in this direction. The case of Kac type interaction has been investigated in $d=1$  by [COP], [COPV] and [OP].   
  
In this paper we  study  a reversible,  nonconservative (Glauber) dynamics of $\pm1$ valued spins,  interacting via a Kac potential and under the influence of an external random field. We  assume the latter given by i.i.d. random  variables taking values $a_i \in \R$ with probability $p_i$, for $i=1,\dots, N$, with $N$ a fixed integer. We do not require the Kac potential to be positive (that is we do not restrict the model to the ferromagnetic case).
 
Our main result is a quenched large deviation principle, almost sure with respect to the random field, for the empirical random magnetizations of this spin-flip process.  
  The rate functional associated with the large deviation principle, which depends on the distribution of the random field, is lower semicontinuous,  positive,  with compact level sets.
{In contrast with the non random case studied in [C], the magnetization $m$ of our spin model is not of mean field type. 
Nevertheless, this difficulty is overcome by coloring the sites according to the random external force, so that the colored magnetizations become a mean field system.} 
The large deviation rate functional is then obtained via a contraction principle from the rate functional associated with the large deviation principle of the empirical  {\it colored} magnetizations $m_i$  (i.e. the magnetization over the sites where the random field takes value $a_i$), $i=1,\dots, N$; {we have $m=\sum_{i=1}^{N} m_i$}.     
 As usual, the rate functional is determined by two distinct types of large deviations of the same order. The first one corresponds to large deviations from the initial state, the second one to the stochasticity of the evolution. Suppose  
 $\AA =\{ \pi^\g (\cdot, dr) \simeq v(\cdot,r) dr , t\in
[0,T]\}$   where  $\pi^\g (\cdot,dr)$ is the local magnetization density, $\simeq$ denotes closeness in some norm and $v$ is a profile  different from the solution of the nonlinear macroscopic equation giving the law of large numbers. 
 We need to modify the measure of the process over the  magnetization profiles so that event $\AA$ becomes typical. One possible choice
is to drive the  spin system by weak, slowly varying, space-time dependent external forces.  
 This is the standard choice for spin systems evolving according  to non conservative (Glauber) evolutions  without randomness involved, or to  conserved (Kawasaki) evolutions with gradient type interactions.  For conserved non gradient systems, the force must be configuration dependent  (see [Sp],
p. 248), to take into account that for these systems the response in the current to an external force field is partially delayed. 
Namely, when an external random field is added to the  Hamiltonian,  in the conservative, non gradient case (see [MO]), one needs to take the external force weakly dependent on  the field randomness. In the non conservative case, it turns out that  the external force  strongly depends on the field randomness.  In other words, in dynamics with a conserved quantity, there is less freedom in choosing the class of perturbations than in non conservative dynamics.  

  We distinguish between sites where the random field takes different values;  on each of them  we take a deterministic space-time dependent external force. 
This allows to write the rate functional associated to the large deviation principle in a closed form with respect to the local {\it colored}  magnetization.  
 We  carried out explicitly the computations for a couple ($N=2$); the general case follows. 
  The simplest case to have in mind is $a_1=1$, $a_2=-1$, $p_1=p_2= 1/2$ and $J\ge 0$;  then, when $\beta$ (which is proportional to the inverse temperature) is such that $ \beta \ge  \beta_c= (\int J(r)\,dr)^{-1}$, interesting phenomena  appear when 
studying the long time behaviour of the spin system.  This is  related to the fact that the underlying spin systems at equilibrium undergo to phase transition, even in one dimension in the limit $\g \downarrow 0$. 
In this paper  we will study the dynamics of the spin system for finite time: in this regime, the evolution does not depend crucially on the value of $\beta$. We will then set $\beta=1$. 
  
 The random Curie-Weiss model (RCW), which describes a mean field interaction, has given rise to many results on short and long time dynamics. 
In [DD], short time dynamics has been studied. More precisely the large deviations for the empirical measures in the  product space of magnetization trajectories and realizations of the random field are given. {}From this result one could derive annealed  large deviations   for the RCW  but not quenched ones. In [MP] and [FMP], long time dynamics, convergence to equilibrium when the random field takes only the two values $\pm \e$ are considered. In [BEGK], the RCW model is analyzed when the random field takes finitely many values, as an example of the use of the potential theoretical approach to metastability. Furthermore, in [BBI], the previous results are extended to continuous distributions of the field, and precise asymptotics of metastable characteristics are derived.   
  
There are no available results for short and long time dynamics of the random field Kac model. We make here a first step in addressing this problem. 
 
In Section 2 we present the model, main definitions and results.  
In Section 3 we define the rate functional associated to the large deviation principle, we exhibit different representations for it, and we give its main properties (lower semicontinuity, compactness). There, we follow the scheme of [C], Section III, but {working with the couple $(m_1,m_2)$ induces intricate computations. Since the spins have value  $\pm 1$, the local and colored magnetizations are always between $+1$ and $-1$. A consequence of the randomness is that the functional becomes infinite for the colored particle system at the boundary
$p_i\le |m_i|\le 1$ ($i=1,2)$ of the coupled magnetization. Thus these boundaries are not rare enough in the large deviations regime, and we have to deal with this lack of regularity.
This is different from} the non random case [C], where the boundary is reduced to the two values $\pm 1$ of the magnetization.
A preliminary step to derive the large deviation principle (LDP) is the hydrodynamic behavior for the colored particle process, sketched in Section 4.
The class of time dependent, random perturbations needed to derive the LDP lower bound is introduced in Section 5, 
where the perturbed process is studied. In Section 6 we derive the upper bound and in Section 7 the lower bound of the LDP.  
The lower bound is obtained first for trajectories that are smooth is space and time,
and outside the boundaries $p_i\le |m_i|\le 1$ ($i=1,2)$. Then it is extended to a larger class $\PP$ of paths,
by smoothing by successive steps trajectories with a finite 
rate functional, using techniques introduced in [QRV]. In this context, $\PP$ consists in trajectories absolutely
continuous with respect to the Lebesgue measure and absolutely continuous in time.
Then, in order for the usual martingale technique to be effective to obtain the upper bound,
we need to show that the process concentrates on $\PP$.   
To this aim, we introduce an energy functional via an 
exponential martingale which  excludes the trajectories not in $\PP$ (in the spirit of [QRV], [MO], [FLM]).
{The  appendices (Sections 8 and 9) gather the most technical proofs. }

\bigskip
\chap{2. The model and the main  results}2
\numsec= 2
\numfor=2
\numtheo=2
\medskip

\noindent{\bf The  space of configurations:} 
Let $\L$ be the $d$-dimensional torus of  diameter 1. For $0<\g<1$
such that $\g^{-1} \in \N $, $\L_\g =\Z^d/\g^{-1}\Z^d $ is the $d$-dimensional discrete
torus of diameter $\g^{-1}$. 
We denote by $\SS_\g \equiv  \{-1,+1\}^{\L_\g}$ the configuration space  and  by  $\s =(\s(x))_{x\in\L_\g}$ a spin configuration, where  for each $x\in \L_\g$, $\s(x) \in\{-1,1\}$.  
\medskip

\noindent{\bf The disorder:} It is described by a collection of i.i.d.  random variables $\a=\{\a(x), 
x\in \Z^d\}$ taking two values,  i.e. $\a(x) \in  \{a_1, a_2\}$. 
The corresponding product measure on  $\O= \{a_1, a_2\}^{\Z^d}$ is denoted by $\P$ (and $\E$ is the expectation with respect to $\P$),   
$$\P\{\a(x)=a_i\}=p_i, \quad i=1,2 . \Eq(g-alpha)$$
For $\g^{-1}$ an odd integer, $\a$ induces in a natural way a random field on $\L_\g$,   also  denoted by $\a$. 
  
\medskip
\noindent{\bf The Kac potential:} 
 We consider a pair interaction among particles given by a Kac potential of the form
$$J_\g(x,y)\equiv \g^d J(\g (x-y)), \quad (x,y) \in \L_\g \times \L_\g, \Eq(Jgamma) $$  
where    $J: \L  \to \R $ is a symmetric function, that is $J(r)=J(-r)$, such that  $\int J(r){d}r=1$ (normalization). The interaction  $J$ might have any sign.  
Denote by $\CC(\L)$ (resp. $\CC^1(\L)$, $\CC^2(\L)$) the space of
continuous (resp. continuously differentiable, twice continuously
differentiable) real functions on $\L$. We assume $J \in \CC^1 (\L)$.

\medskip
\noindent{\bf The  Energy:}
Given a realization  $\a$ of the magnetic field,  define  for all $\g$, $  \th>0 $, $\s\in\SS_\g$,   the  Hamiltonian 
$$H^{\g,\a}(\s)= - \frac 12\sum_{(x,y) \in \L_\g \times \L_\g} J_\g(x,y) \s(x) \s(y) -\th\sum_{x\in \L_\g} \a(x)\s(x),\Eq(g1)$$
and  the  Gibbs measure  $\mu^{\g,\a,\beta}$  associated to 
$H^{\g,\a}$ at inverse temperature $\beta$, with normalization constant  $Z^{\g,\a,\beta}$: 
$$
\mu^{\g,\a,\beta}(\s) = \frac{1}{Z^{\g,\a,\beta}}\exp\big[ -\beta H^{\g,\a}(\s) \big]\;.
$$ 
\medskip
\noindent{\bf The Glauber dynamics:}
Denote  by  $\s^x $ the  configuration obtained from $\s$  by flipping the spin at site $x$:
  $$ 
\s^x (z) = \left \{ \eqalign { & -\s(x)\qquad \hbox { if } \qquad  z=x,
\cr & \s(z) \qquad \hbox { otherwise, }  }\right. $$
so that the energy difference  resulting from a spin flip at $x$ is 
 $$ 
H^{\g,\a}(\s^x)- H^{\g,\a}(\s)
=2\s(x)\left[ (J_\g \star \s)(x) +\th \a(x)\right],\Eq (g210)
$$
where without loss of generality we have assumed  
$J(0)=0$, 
and we define the discrete convolution $\star$ between function $J_\g$ and a configuration $\s$
 by 
$$
(J_\g \star \s)(x) = \g^d \sum_{y \in \L_\g} J(\g(x-y))  \s(y).\Eq (g51)
$$
 We consider  
a Markovian evolution on $\SS_\g$, whose generator $ \LL^{\g,\a}$ acts on cylinder functions $f$ as 
$$
\LL^{\g,\a} f (\s)=\sum_{x \in \L_\g}  c_x^{\g,\a}(\s) [f(\s^x)- f(\s)], \Eq (g3)   
$$
where, for $x\in \L_\g$,  
$$ 
c_x^{\g,\a}(\s) = 
\frac {\exp[ -(\b/2) (H^{\g,\a}(\s^x)-H^{\g,\a}(\s))]}{2\cosh[(\b/2)
(H^{\g,\a}(\s^x)-H^{\g,\a}(\s))]  }.
\Eq (g2)$$
Then $\LL^{\g,\a} $  viewed as an operator on $L^2(\mu^{\g,\a,\beta}) $ is self-adjoint.
 Since temperature is kept fixed in all the paper and does not play any role we set for simplicity $\b=1$. 
 We fix a time $T>0$, and we will study the process $(\s_t)_{t \in[0,T]}$   with infinitesimal generator given in \equ(g3). 
\medskip
 \noindent 
   { \bf The  measure spaces: }
 Let $\MM_1 $ be the set of signed Borel measures $\mu$ on the Borel $\sigma$-field of  $\L$ with total variation norm bounded by 1. 
We equip   $\MM_1 $ with the weak  $\tau^*$ topology induced
by $\CC(\L)$ via $<\mu,G>=\int G {d}\mu$ (for $G\in\CC(\L)$).
We denote by $\rho(\cdot,\cdot)$  the distance which makes $ (\MM_1, \tau^*)$ a metrizable compact space, see [Bill]: 
that is, given   $(H_k)_{k \in \N}$   a dense subset in the unit ball of $\CC(\L)$
 for $ \mu_i \in \MM_1$,  $i=1,2$, 
$$
\rho(\mu_1,\mu_2)=\sum_{k\ge 0}2^{-k}|<\mu_1-\mu_2,H_k>|. \Eq (par8)
$$
  Let $0<q\le1$, and 
   $$  \MM_q^{ac} = \left \{ \mu \in \MM_1: \mu << \l \quad  \hbox {and} \quad   \left|\frac {d \mu } {d \l }\right|  \le q \quad \l\,-a.s. \right  \}, \Eq (Ac.1)$$ 
  where $ \l$ is  the  Lebesgue measure  on $\L$. 
 We identify $\mu \in \MM_q^{ac}$  with its Radon-Nikodym derivative $\displaystyle{\frac {d\mu} {d\l}}$, and, by an abuse of notation, we write $ < \mu,G>= <\displaystyle{\frac {d\mu} {d\l}},G>$. Since  $\MM_q^{ac}$ is  a closed ball of $\MM_1$, it is    $\tau^*$ compact.   

  If $\s \in \SS_\g$ we define the  {\sl empirical measure}  $\pi^\g (\s)\in\MM_1$  by   
$$
\pi^\g(\s)(d r )=\g^d\sum_{x\in\L_\g} \s (x)\d_{\g x}(d r), \Eq (g20)
$$
where $\d_{ \g x} $ is the Dirac  measure concentrated on point $\g
x$.
 Remark that if we denote by $\mu*G$ the
convolution of a measure $\mu$ and a function $G$ over $\L$, namely
$(\mu*G)(r')=\int_\L G(r'-r)\,\mu(dr)$, then we can rewrite 
$$
(J_\g \star \s)(x)=(\pi^\g(\s)*J)(\g x).
\Eq(g5)
$$
  We denote by  $D([0,T],\MM_1)$ (resp. $D([0,T], \SS_\g)$) the space of functions from $[0,T]$ to $\MM_1$ (resp. to $\SS_\g$) that are right continuous  with left limits, endowed with the  Skorohod topology, see [Bill]. 
  \medskip
\noindent{\bf The initial condition:}       
Let  $ (\s^\g)_\g$  be a sequence of configurations such that
$\pi^\g (\s^\g)$ converges when $\g \to 0$ in the weak topology to the measure $m_0\l$, for a continuous function $m_0:\L\to[-1,1]$.
This means that
$$ \lim_{\g \to 0}  \rho ( \pi^\g (\s^\g), m_0\l)  =0. \Eq(hyp-G1a)$$ 

We denote by $P^{\g,\a}_{\s^\g}$ the law (and by $E^{\g,\a}_{\s^\g}$ the expectation) of the process $(\s_t)_{t \in[0,T]}$ on $D([0,T], \SS_\g)$  starting at time $t=0$ from the  
deterministic initial configuration $\s^\g$, and by $Q^{\g,\a}_{\s^\g}$ the law on  $D([0,T],\MM_1)$   of
 the corresponding empirical measure  process $(\pi_t^\g)_{t\in[0,T]}$, where $\pi_t^\g$ stands for $\pi^\g(\s_t)$.    
 \medskip \noindent
We first obtain the ``law of large numbers". 
 \medskip
\noindent {\bf \Theorem (th-g)} 
{\it  Assume $ (\s^\g)_\g, m_0$ satisfy    \eqv (hyp-G1a). For all $ t\ge 0$, 
$$ \lim_{\g \to 0}  \rho ( \pi^\g_t , m(t,\cdot) \l)  =0, \quad   \text { $\P$-a.s.,} \Eq(may1)$$  
where $m (\cdot,\cdot)$ is the unique weak solution of  
$$ 
\left \{ \eqalign {& \partial_t m(t,r)= -m(t,r)+\sum_{i=1,2}
p_i\tanh\left [  (J* m(t,\cdot)) (r)+a_i\th \right ]
\cr & m(0,\cdot)=m_0(\cdot).  }\right. 
\Eq(G1a)$$
Furthermore, for all $ G \in\CC^{0,1}([0,T] \times \L) $ (that is, continuous in its first variable, and continuously differentiable
in its second variable), $\d>0$,
$$
\lim_{\g\to 0} P^{\g,\a}_{\s^\g}\left [\sup_{t\in[0,T]}\left|
< \pi^\g_t,G(t, \cdot)>-<  m(t,\cdot), G (t, \cdot)>
\right| \geq\d\right ] =0. 
\Eq(unif0)
$$ }

\noindent
 By an abuse of notation we write from now on $(J* m)(t,r)$ instead of $(J* m(t,\cdot)) (r)$.
 
\noindent
\Remark (R-th-g)    The   Cauchy
problem  \eqv (G1a) in this setup is well posed with a unique
global solution, because the right hand side of
 \eqv (G1a)  is uniformly Lipschitz, and
because  the set $\{  m \in  L^\infty (\L): \|m\|_{\infty}   \le1 \}$ is left invariant,
since $|\tanh z| \le  1$ for all $z$.   Furthermore the solution is differentiable in time.  
  \smallskip 
\noindent Next we state the {\it quenched} large deviation principle for   $Q^{\g,\a}_{\s^\g}$. Different choices of initial conditions could  be treated as well.
  The only difference would be an extra term to add to the rate functional associated with the large deviation principle $\wt I_{m_0}(\cdot)$, taking into account the  deviation from the initial profile at time $t=0$. The functional  $\wt I_{m_0}(\cdot)$ depends on the distribution of the random field but 
not on its realization; it is obtained through a contraction principle, as explained in the introduction. Its explicit formulation relies on  several intermediate steps.
 Let 
$$ \DD (\wt I_{m_0})=  \{\pi  \in D([0,T], \MM_1):   \wt I_{m_0}( \pi) <\infty \}. \Eq (IHP1)  $$ 
\medskip 
\noindent {\bf\Theorem   (main-ldp)} {\it   
 Assume $ (\s^\g)_\g, m_0$  satisfy  \eqv (hyp-G1a). 
  For all closed subsets  $\FF \subset D([0,T], \MM_1)$ and open subsets $\AA\subset D([0,T], \MM_1)$, we have  
$$ \lim\sup_{\g \to 0} \g^d \log   Q^{\g, \a}_{\s^\g} (
\FF) \le
- \inf_{ \pi  \in \FF} \wt I_{m_0} ( \pi), \quad \P-a.s., \Eq (2.6) $$ 
$$ \lim\inf_{\g \to 0} \g^d \log  Q^{\g, \a}_{\s^\g} (
\AA) \ge - \inf_{ \pi   \in \AA} \wt I_{m_0} ( \pi ),  \quad \P-a.s. \Eq
(2.7) $$    
The functional $\wt I_{m_0}(\cdot)$, defined in \eqv (LL.1) below, is non-negative for   $\pi \in D([0,T], \MM_1)$, lower semicontinuous with compact  level sets and, see Definition \eqv (D10) later on,   
  $$\DD (\wt I_{m_0}) \subset  \{m  \in \CC([0,T], \MM_1^{ac}): m(t,.) \quad \hbox {absolutely  continuous for}\quad t \in [0,T] \}. $$ 
 }
 
\medskip
\noindent
{\bf The colored particle system: }
 To derive the rate functional associated with the large deviation principle we  introduce {\sl random} empirical measures
 $\overline \pi^\g= \big( \pi_1^\g,   \pi_2^\g)$. 
For $\a\in\O,x\in\L_\g$,  $i=1,2$, set
$$ \a_i(x)= \1_{\{\a(x)=a_i\}},  \quad 
\Eq (1-alpha) $$ 
$$
\pi_i^\g(\s)( {d} r )=\g^d\sum_{x\in\L_\g}\a_i(x) \s (x)\d_{\g x
}( {d} r).  
 \Eq (g200)$$
Though we do not write   it explicitly,  $\pi_i^\g(\s) \in\MM_1$ 
depends on the randomness.
Moreover the knowledge of $\pi_i^\g(\s)$ for $i=1,2$  determines    $\pi^\g(\s)=\pi^\g_1(\s)+\pi^\g_2(\s)$.
 We denote   by 
$\mus Q^{\g,\a}_{\s^\g}$  the law on  $D([0,T],\MM_1 \times \MM_1)$   of the   empirical measure  process $(\overline\pi_t^\g)_{t\in[0,T]}=(\pi_{1,t}^\g, \pi_{2,t}^\g)_{t\in[0,T]}$ 
under $P^{\g,\a}_{\s^\g}$.
We denote,   for  $\overline G = (G_1, G_2 )    \in \left(\CC(\L)\right)^2$,
$$<\overline \pi_t^\g, \overline G >=  \sum_{i=1,2}  \g^d\sum_{x\in\L_\g} G_i (\g
x) \a_i(x)  \s_t (x) 
  \Eq(g21)$$
and,  for  $\overline m= (m_1, m_2)\in (L^\infty (\L))^2$,
 by an abuse of notation,
$$ < \overline m,\overline G >= < (m_1\l, m_2\l),\overline G >=\sum_{i=1}^2  \int_\L  G_i (r)    m_i(r)\, dr.  \Eq(g21a)$$

\medskip
\noindent {\bf \Theorem (2th-g)} 
{\it  Assume    $ (\s^\g)_\g, m_0$    satisfy    \eqv (hyp-G1a). 
 For all $t \in [0,T]$,  $\d>0$ and  $\overline
 G\in\big(\CC^1 (\L)\big)^2  $,
 $$
 \lim_{\g\to 0} P^{\g,\a}_{\s^\g}\left [\left|<
 {\mus \pi}_t^\g,\mus G >- 
 < \mus m(t,\cdot),\mus G> \right| \geq\d\right ] =0  \quad
 \P-\hbox{a.s.}\; ,
 $$
where  $\mus  m= (m_1,m_2)$  is the unique weak  solution 
  of  
$$ \left \{ 
  \eqalign {
   &\partial_t m_i(t,r)= -m_i(t,r)+  p_i \tanh \left [
  \b((J* m)(t,r)+ a_i\th)  \right ],  
      \cr & 
 m= m_1+m_2;\quad m_i(0,\cdot)=p_i m_0(\cdot),\;\quad     i=1,2. }\right.\; 
\Eq(G2b)$$ 
}
 \bigskip

 \noindent
\Remark (R-2th-g) 
     Similarly to Remark \eqv (R-th-g), the   Cauchy
problem  \eqv (G2b) in this setup is well posed with a unique global solution;
here, the set $\{ \mus m \in  (L^\infty (\L))^2 :\|m_i\|_{\infty}   \le p_i,  i=1,2\}$ is left invariant.
The    solution  is differentiable in time. 
The case  $J\ge0$, $a_1=1$, $a_2=-1$, $p_1=p_2= 1/2$  is  analyzed in  [COP4]. \medskip \noindent
  To derive  Theorem \eqv (3A) below, we need a stronger type of convergence: 
\medskip \noindent
  {\bf \Corollary (ch-g)} {\it 
  For all $\mus G \in\left (\CC^{0,1}([0,T] \times \L) \right)^2$, $\d>0$,
$$
\lim_{\g\to 0} P^{\g,\a}_{\s^\g}\left [\sup_{t\in[0,T]}\left|
<\mus \pi^\g_t, \mus G(t,\cdot)>-<\mus m(t,\cdot),\mus G(t,\cdot)>
\right| \geq\d\right ] =0. 
$$}  

\medskip 
\Remark (R-3th-g)   Theorem \equ(2th-g) and Corollary  \eqv (ch-g) imply Theorem \equ(th-g)  since 
if $ \mus G= (G,G)$,
$$<\overline \pi_t^\g, \overline G >= < \pi_{1,t}^\g, G >
+< \pi_{2,t}^\g, G >=< \pi_{t}^\g, G >.$$  
\medskip 
\noindent
Next theorem states the large deviation principle for the colored particle system. Theorem  \eqv (main-ldp) is based on this important intermediate result, interesting for itself.
   
 \medskip 
\noindent {\bf \Theorem   (3A)}  {\it  Assume $ (\s^\g)_\g, m_0$ satisfy \eqv (hyp-G1a). We have, for all  
open  subset  $\mus \AA$ and  closed subset $ \mus \FF$ in $D([0,T], \MM_1 \times \MM_1)$, 
$$ \lim\inf_{\g \to 0} \g^d \log  \mus Q^{\g, \a}_{\s^\g} ( \mus \AA) \ge -
\inf_{  \mus \pi   \in \mus \AA} I_{m_0} (  \mus \pi ), \quad \P-a.s. \Eq (2.3) $$
$$ \lim\sup_{\g \to 0} \g^d \log  \mus  Q^{\g, \a}_{\s^\g} ( \mus \FF) \le -
\inf_{ \mus \pi  \in \mus \FF} I_{m_0} (    \mus\pi), \quad \P-a.s. \Eq (2.5) $$
where 
$$ I_{m_0} (  \mus\pi) = \left \{ \eqalign { & I_0( \mus\pi) \quad \hbox {if} \quad  \pi_i(0, \cdot) =  p_i m_0(\cdot) \l,   \quad  i=1,2,\cr & +\infty  \quad \hbox {otherwise,} } \right. \Eq (2.4)$$
and $ I_0(\cdot)$, defined in \eqv (funct1) below, is lower semicontinuous with compact level sets. 
} 
\medskip
 
\noindent Define, for a path $\pi\in D([0,T], \MM_1)$,  
$$
\wt I_{m_0}(\pi)=\inf\Bigl(
I_{m_0} ( \mus \pi),\  \mus \pi=(\pi_1,\pi_2),\pi_i\in D([0,T], \MM_1),i=1,2,\, \pi_1+\pi_2=\pi
\Bigr).  \Eq (LL.1)
$$
Since the map $ (\pi_1,\pi_2) \mapsto \pi_1+\pi_2$ is continuous  in  $D([0,T], \MM_1\times \MM_1)$, by the contraction principle, see [V], [DZ], Theorem \eqv (3A) proves  Theorem  \eqv  (main-ldp). Therefore in the following  sections we will focus on the colored particle system.

\bigskip
\chap{3. Rate functional }3
\numsec=3
\numfor= 1
\numtheo=1

In this section we define the rate functional $ I_0 (\cdot)$ of the colored particle system  and state its main properties. Proofs of the latter, quite technical, are carried out in Section 8.  Heuristics to define $ I_0 (\cdot)$ consists in finding, for any path $ \mus \phi $ on $[0,T]$  smooth enough, an exponential change of probability under which the  process  $ (\sigma_t)_{t \in [0,T]}$   
is uniformly close to $ \mus \phi $ on $[0,T]$. When there exists some potential $\mus V(t,\cdot)$, $ t \in [0,T]$ smooth enough for $ \mus \phi$ to be the solution of {\it a perturbed equation} (obtained by the law of large numbers from the process $ (\sigma_t)_{t \in [0,T]}$,  see \eqv (Eq.1) later on), then  $ I_0 (\cdot)$ is related to the Radon-Nikodym derivative of the distribution of $ (\sigma_t)_{t \in [0,T]}$ with respect to the distribution of the original process,    
 see Theorem \eqv (RN). In the general case, i.e. when there is no  such $\mus V(t,\cdot)$, we are still able to provide an explicit  representation of $ I_0 (\cdot)$  (this is similar to the results of [C]).  We will then show that this representation of $ I_0 (\cdot)$ is  equivalent to the usual definition of the rate functional, given through the macroscopic functional associated to the Radon-Nikodym derivative, see \eqv (2.2).
We start by specifying the functional spaces on which we will define  $ I_0 (\cdot)$.
  For $(p_1,p_2)\in[0,1]^2$,  we identify the set
$$
B_{p_1,p_2}=\{ \mus u=(u_1,u_2): u_i \in L^\infty ( \L),\,   \|u_i\|_\infty\le p_i,\,i=1,2 \} \Eq
(Aco.2)$$
with $  \MM_{p_1}^{ac} \times   \MM_{p_2}^{ac} $, see \eqv     (Ac.1),  and extend   the distance $\rho$ (see
\eqv (par8)) to elements of $  \MM_{p_1}^{ac} \times    \MM_{p_2}^{ac} $  by
  $$   \r(\mus \mu, \mus  \nu)=  \sum_{i=1,2} \r(\mu_i, \nu_i).  \Eq (metric3a)  $$

\smallskip
\noindent  \Definition (D10)
 { \it\  Let  $\ac([0,T], B_{1,1})\subset\CC([0,T], B_{1,1})$      be  the  subset  of absolutely continuous functions  $\mus \phi= (\phi_1,\phi_2)$, that is, for $j=1,2$:
 for all $t'\in [0,T],\, t\in [t',T]$,  
 there exists  $\dot  \phi_j \in L^1([0,T] \times \L)$ such that}
$$\phi_j(t)(r)- \phi_j(t')(r)= \int_{t'}^t \dot \phi_j(s,r )\, ds\, , \quad \l - a.s. $$
 By an abuse of notation, from now on we write $\phi_j(t,r)$ instead of $\phi_j(t)(r)$.

\medskip   
 
To write $ I_0 (\cdot)$,
we start by defining, for each $t\in[0,T]$, the following functionals, in which time is kept fixed,  therefore we omit to write   it. 
 For $ \mus  \pi=(\pi_1,\pi_2)\in\MM_1\times\MM_1$ (we write $\pi=\pi_1+\pi_2$), $\mus\mu=(\mu_1,\mu_2) \in
   \MM_1\times\MM_1$  and $\mus V=(V_1,V_2) \in (L^\infty (\L))^2$ denote 
$$ \eqalign { F_{\mus V}(\mus\mu,\mus\pi)  &= 
  \sum_{i=1,2}<\mu_i, \tanh (  \pi*J +a_i\th) \sinh (2V_i) + \cosh (2V_i) -1>  \cr &
 -  \sum_{i=1,2} <\pi_i, \tanh (  \pi*J + a_i\th)[ \cosh (2V_i)-1] + \sinh (2V_i) >, 
 } \Eq (rate1) $$ 
 and for $\mus u = (u_1,u_2) \in B_{1,1} $,  $ \mus g \in \left(L^1(\L)\right)^2$, 
 $$\G_{\mus V}(\mus   u ) = F_{\mus
 V}\bigl((p_1\l,p_2\l),(u_1\l,u_2\l)\bigr),
 \Eq (2.1) $$
 $$  \HH^* ( \mus u, \mus g) =  \sup_{\mus V \in (L^\infty ( \L))^2}  [ <\mus V ,\mus g > -    \frac 12 \G_{\mus V}(\mus u )].  \Eq (2.2a) $$  
The function  $\mus g \to \HH^* ( \mus u,   \mus g)$  is convex. Next lemma ensures that $\HH^*( \mus u,\cdot)$ is the Fenchel-Legendre transform of $\G_{(\cdot)}(\mus u)$ when  $\mus u\in B_{p_1,p_2}$, and we will derive in that case an explicit formula for $\HH^* ( \mus u,\mus g)$.
   
\medskip 
\noindent {\bf \Lemma  (2C1)}  {\it 
As a function of $\mus V\in  (L^\infty ( \L))^2$, $ \G_{\mus V}(\mus u) $ is convex  differentiable  for   $\mus u \in B_{p_1,p_2}$. }
  
 \smallskip\medskip 
 
\noindent  \Definition  (D1)   
 {\it The dynamical rate functional  $I_0:  D([0,T], \MM_1  \times \MM_1 )\to \R \cup \{\infty\} $ is given by  
  $$I_0( \mus \pi ) =\left \{ \eqalign { &  I_0(\mus \phi)= \int_0^T  \HH^* (\mus \phi(s,\cdot), \dot {\mus \phi} (s,\cdot) )  ds,    \quad
 \hbox{for } \mus\pi = (\phi_1 \l, \phi_2 \l), \mus \phi= (\phi_1, \phi_2)  \in   \ac ([0,T], B_{1,1}),   \cr &     \infty  \quad \hbox {otherwise.} } \right.  \Eq (funct1) $$
}
  
\medskip \noindent 
To derive properties of the rate functional associated with the large deviation principle it is convenient to have different representations of $I_0$. To this aim let  
 $ \mus V= (V_1, V_2) \in \left(L^\infty   ( [0,T] \times \L)\right)^2 $.  We  define, for   $\mus \pi \in  D([0,T], \MM_1  \times \MM_1 )$ (cf. \eqv (2.1)), 
$$   K_{\mus V}(\mus \pi)=\left \{   \eqalign { &   K_{\mus V}(\mus \phi), 
 \ \hbox{ for }   \mus\pi = (\phi_1 \l, \phi_2\l),\ \mus \phi= (\phi_1, \phi_2)  \in   \ac ([0,T], B_{1,1}), 
  \cr & \infty    \quad
\hbox  {otherwise}, } \right . \Eq (DC1)$$
where 
$$K_{\mus V}(\mus \phi)=
\int_0^T < \mus V(s,\cdot), \dot {\mus \phi} (s,\cdot)> ds -
  \frac 12  \int_0^T \G_{\mus V(s,\cdot)}( \mus \phi (s,\cdot))  ds, $$
    $$J_0(  \mus \pi )= \sup_{\mus V \in (L^\infty   ( [0,T] \times
  \L))^2}   K_{\mus V} (  \mus \pi ).  \Eq (2.2) $$  
  $$J_1(\mus \pi)= \left \{ \eqalign { & J_1(\mus \phi)=  \int_0^T
  \int_{\L} \HH(\mus  \phi (t,r),   \dot {\mus  \phi} (t,r))\,dr\,dt,\    
  \hbox{ for } \mus\pi = (\phi_1 \l, \phi_2 \l),\ \mus \phi= (\phi_1, \phi_2)  \in   \ac ([0,T], B_{1,1}), 
 \cr &     \infty  \quad \hbox{ otherwise}, } \right.  \Eq (funct2) $$
    where for $\mus u=(u_1,u_2)\in B_{1,1}$, $\mus g=(g_1,g_2)\in (L^1([0,T]\times\L))^2$, $(t,r)\in[0,T]\times\L$, 
$$ \HH (\mus u,  \mus g )(t,r)=\HH (\mus u (t,r),  \mus g (t,r))= \sum_{i=1}^2 H_i (\mus u,   g_i)  (t,r),
\Eq (Le2)$$
$$  H_i (\mus u,   g_i)(t,r) =  \sup_{v_i\in \R }  \left \{    g_i (t,r)  v_i- \frac 12  B_i(\mus u(t,r), v_i) \right \}, \quad  i=1,2,  \Eq (Le2c)$$  
$$  \eqalign{B_i(\mus u (t,r),  v_i)=  
 &( p_i-u_i(t,r))  \frac { e^{  (J* u)(t,r)   + a_i\th}  } { 2\cosh  [ (J* u)(t,r)  + a_i\th] }  \left [ e^{ 2v_i  }-1 \right ]\cr & + ( p_i+u_i(t,r))  \frac { e^{- [  (J* u)(t,r)    + a_i\th]}  } { 2\cosh   [  (J*u)(t,r)  + a_i\th] }  \left [ e^{ -2v_i  }-1 \right ].} \Eq (2.1n) $$ 
When $\mus u (t,r)\in[-p_1,p_1]\times[-p_2,p_2]$,  
$\sum_{i=1}^2 B_i(\mus u(t,r),\cdot)$ is convex so that $\HH (\mus u (t,r), \cdot  )$ is its
Fenchel-Legendre transform. 
 We now give an explicit representation of $  \HH (\cdot, \cdot)$.
  To simplify notations denote    
 $$\eqalign { &  A_i= A_i(u,\th) (t,r)= (J * u)(t,r)+a_i\th, \cr &
    R_i=R_i(\mus u, g_i,\th)(t,r)=\sqrt { \left ( g_i(t,r) \cosh [A_i(u,\th) (t,r)]\right ) ^2+   p_i^2- u^2_i (t,r)}   , \cr & 
D_i=  D_i(\mus u, g_i,\th) (t,r)=   g_i(t,r)  \cosh  [A_i(u,\th) (t,r)]+ R_i(\mus u, g_i,\th) (t,r). 
}  \Eq (DD1a) $$ 
Note that $D_i (\mus u,g_i,\th)(t,r) \ge 0$   regardless  of the sign of $g_i(t,r)$. 
When $ (t,r) $ is kept fixed  we  omit to  write it.
The function $\hbox{\rm {{\rm sgn}}}:\R\to \R$ is given by
$$ 
\hbox{\rm {{\rm sgn}}}(x)= 
 \left \{ \eqalign { & 
\frac {x}{|x|}\qquad \hbox{ if }\ \ x \not=0 \, ,\cr
&0  \qquad \hbox{ if }\ \ x=0  \, .
}\right. \Eq (sng)
$$ 
 
\medskip
\noindent {\bf \Proposition  (3A22)}  {\it 

\noindent
(a) If  $|u_1| >p_1$ or   $|u_2| >p_2$,   then 
$ {\HH} (\mus u,  \mus g)=+\infty$. 

\noindent
(b)  For $i=1,2$, when  $|u_i|< p_i $,  then  
        $$  H_i (\mus u,   g_i)  =           
     \frac {g_i}2 \left(\log \frac { D_i } { p_i- u_i  }-  A_i \right)  
     +\frac {p_i}2 -\frac {u_i}2 \tanh  A_i   - \frac { R_i}{2\cosh  A_i }.   \Eq (DD1c) 
$$ 
 (c) For $i=1,2$, when  either  ($u_i=p_i $    and    $g_i\le 0$) or  ($u_i=-p_i$ and  $ g_i\ge 0$), then 
  $$ 
       H_i(\mus u,   g_i) =  \1_{\{g_i\neq 0\}} \frac {| g_i|} 2 \left( \log \left \{ \frac { |g_i| \cosh   A_i }   {p_i  e^{-{\rm sgn}(u_i) A_i}  }\right \} 
-1 \right)     +  p_i \frac { e^{-{\rm sgn}(u_i) A_i}  } { 2\cosh   A_i } . 
\Eq (DD1b)$$
(d) For $i=1,2$, when  either ($ u_i=p_i$ and $g_i>0$)  or ($ u_i=-p_i$ and $g_i<0$), then
$H_i (\mus u,g_i)=+\infty$.
}
 \medskip
  The following proposition shows that the order of supremum and the integrals can be reversed.  In particular we can compute the supremum for each point $(t,r) \in [0,T] \times \L$.  
    
\medskip
\noindent {\bf \Proposition  (3A1)}  {\it
For  $\mus\pi = (\phi_1 \l, \phi_2 \l)$, $\mus \phi= (\phi_1, \phi_2) \in  \ac ([0,T], B_{1,1})$, we have $I_0(\mus \pi) = J_0(\mus\pi ) =  J_1(\mus \pi).$
Furthermore
if $\mus \phi \in \ac ([0,T], B_{1,1}) \setminus   \ac ([0,T],  B_{p_1,p_2})$, then  
  $I_0(\mus \phi) =+\infty.$  }

\medskip
 \noindent {\bf  Proof.}  
 This follows and extends the proof in [C], p. 171, Properties III(a). By their respective Definitions \eqv(funct1), \eqv(2.2), \eqv(funct2) (see also \eqv(2.2a), \eqv(DC1), \eqv(Le2)), we have
$ J_0(\mus \pi) \le  I_0(\mus \pi) \le J_1(\mus \pi)$. We now prove that we have equalities.
In all cases, for $i=1,2$, we denote by $\vartheta_i$  the value of $v_i$ 
that realizes the extremum of $H_i (\mus u, g_i)$. {}From Proposition \eqv(3A22), 
$\vartheta_i$ belongs to $\R\cup\{+\infty,-\infty\}$. Let
 $\vartheta_i^m= {\rm sgn}(\vartheta_i)\times[|\vartheta_i|\wedge m]$
 and $b_i^m$ be the corresponding (finite) value of $H_i (\mus u, g_i)$. Then as $m\to\infty$, $\vartheta_i^m\to \vartheta_i$  and $b_i^m\to H_i (\mus u,    g_i)\in\R^+\cup\{\infty\}$. According to the case we consider, either $a_i$ and/or $b_i$ are finite, and there is no problem, or $b_i=+\infty$ thus $b_i^m>0$ for $m$ large enough, or, when $\mus u\in B_{p_1,p_2}$, $b_i^m$
is non-negative because $\vartheta_i^m$ is between 0 and $\vartheta_i$, and 
 $v_i\mapsto B_i(\mus u, v_i)$ is a convex function.  
 Therefore in all cases we apply  Fatou's Lemma to get $$J_1(\mus \pi)\leq \liminf_{m\to\infty}\int_0^T \int_\L\sum_{i=1,2}b_i^m(t,r)\,drdt $$
which is smaller than $J_0(\mus \pi)$, whence the result. Notice that this implies that $I_0(\mus \pi)$ is infinite when $\mus u\notin B_{p_1,p_2}$.\qed 

\medskip
  Next   we  characterize  the finite energy trajectories. 

\smallskip
\noindent {\bf \Proposition  (3A3)} {\it  Take $(t,r)\in[0,T]\times\L$.

(a)
 Let $\mus u$ be such that for $i=1,2$, $|u_i | < p_i$.  There exist positive constants $K_1$,$K_2$ and $C$ such that }
$$\eqalign { &  \HH (\mus u,  \mus g) (t,r) \cr &  \le \sum_{i=1,2} \frac {|g_i|} 2 \left [  \left(\log |g_i|\right)^+  + \1_{ \{g_i>0\}}  \left(\log \frac 1 {p_i-u_i}\right)^+  +  \1_{ \{g_i<0\}}  \left(\log \frac 1 {p_i+u_i}\right)^+   + K_i\right ](t,r)   +C  
 }  \Eq  (2.1d) $$ 
$$  \HH (\mus u,  \mus g) (t,r)  \ge \sum_{i=1,2}  \frac {|g_i|} 2 \left [ \log |g_i|   - K_i\right ](t,r)   
  -  C. \Eq  (2.1e) $$
 {\it  (b)  $I_0(\mus \phi) < \infty $ if and only if for $i=1,2$, 
$\dot \phi_i \log |\dot \phi_i|$,   
$ \displaystyle{\dot \phi_i  \log \frac 1 {p_i-\phi_i}  \1_{\{\dot \phi_i>0\}}}$, $\displaystyle{\dot \phi_i \log \frac 1 {p_i+\phi_i}  \1_{\{\dot  \phi_i<0\}}}$ belong to $L^1 ( [0,T]\times \L)$}.

  \medskip
\noindent {\bf \Proposition  (2A)}  {\it   
(1) The functional $I_0 (\cdot) $ is lower semicontinuous on
$D([0,T], B_{1,1})$.

(2) The set $ D_{L_0}= \left \{ \mus  \pi;   I_0( \mus  \pi) \le L_0\right \} $ is compact in  $D([0,T], B_{1,1}
)$ for all  $L_0>0$.  

(3)  $I_0( \mus  \phi) \ge 0$,  and $I_0( \mus  \phi  ) = 0$ if and only if $ \mus  \phi $ is the solution of equation  \eqv (G2b).}

\bigskip
\chap{4. Hydrodynamic behavior for the colored system}4
\numsec=4
\numfor=1
\numtheo=1

In this section, we prove Theorem \equ(2th-g), through a by now standard scheme. 
Nevertheless, we detail it since many of its parts will also appear in the following sections.
\smallskip
We first highlight that throughout the paper, one of the key ingredients to deal with the
 randomness of the interaction will be the    
following  applications of the ergodic theorem and strong law of large numbers. 
For all function $h$ on $\L_\g$, integer $l$, we denote by
$h^{(l)}$  the averaged function
$$
h^{(l)}(x)=\frac 1{(2l+1)^d}\sum_{y\in\L_\g,|y-x|\le l}h(y)\; ,\quad x\in\L_\g.
\Eq (microav)
$$ 
\medskip 
\noindent {\bf \Lemma (erg)} (ergodic theorem for local functions)
{\it Let $\Theta (\a)$ be a bounded measurable cylinder function on
$\O$ and $G\in \CC(\L)$. Then, for almost any disorder configuration $\a$,
$$
\lim_{\g \to 0} \g^d \sum_{x\in \L_\g}G(\g x) \tau_x \Theta (\a)=
\E \big[ \Theta \big]\int_\L G(r) dr \, .
$$}

\smallskip
\proof 
Write 
$$ \g^d \sum_{x\in \L_\g}G(\g x) \tau_x \Theta (\a) =  \g^d \sum_{x\in \L_\g}G(\g x) \big[ \tau_x  \Theta (\a) - \E [ \Theta]  \big] +   \E  [ \Theta ] \g^d \sum_{x\in \L_\g}G(\g x). $$
For any $l \in \N$, by the regularity of $G$,
$$
\big |\g^d \sum_{x\in \L_\g}G(\g x) \big [ \tau_x  \Theta (\a) - \E [ \Theta \big] ] \big | \le  \|G\|_{\infty} \g^d \sum_{x\in \L_\g}  \big | 
 \left(\tau_{\cdot}   \Theta (\a)\right)^{(l)}(x)  
  - \E \big[ \Theta]\big|  + \e(\g l),$$
where $ \lim_{s \to 0} \e (s)=0$.
Keeping $l$ fixed, by the ergodic theorem,
$$ \lim_{\g \to 0} \g^d \sum_{x\in \L_\g}  \big | \left(\tau_{\cdot}   \Theta(\a)\right)^{(l)}(x)  
 - \E  [ \Theta]\big| = \E \big [\big | \left(\tau_{\cdot}   \Theta  \right)^{(l)}(x)
  -  \E  [\Theta] \big |   \big]. $$ 
The law of large numbers (letting $l\to \infty$)  gives the result. \qed

\medskip\smallskip  \noindent
We introduce (cf. [K]), for $i=1,2$ and $\d>0$, 
$$
 A_{l, \d}(x,i)=\left\{\a\in\O:\left| \a_i^{(l)}(x)-\E(\a_i(x)) \right|\le   \d \right\},\quad x\in\L_\g,\Eq(gAlbis) $$ 
 $$ \EE_i (\d, l,\g,\a)= \g^d  \sum_{x\in\L_\g}    \1_{   A^c_{l,\d}(x,i)}(\a). \Eq (ju1)$$ 

\medskip 
\noindent {\bf \Lemma (erg1)}  { \it For any $\d>0$, for $i=1,2$, 
$ \lim_{l \to \infty} \lim_{\g \to 0} \EE_i (\d, l,\g,\a) = 0,  \quad \P-a.s.$  
 }
 \smallskip
  
\proof  
Applying     Lemma \eqv (erg)   to the function $\displaystyle{\Theta=\1_{ A^c _{l, \d} (0,i) }}$ gives   $
 \lim_{\g\to 0}\EE_i (\d, l,\g,\a) = \P( A^c_{l,\d}(0,i) ),  \quad  \P -a.s.
  $ 
Then by  the  strong law of large numbers  $\lim_{l \to \infty} \P( A ^c_{l,\d}(0,i) )=0$. \qed 
\medskip
  \noindent In  the following  it is convenient  to define 
  the  random discrete measures $\mus \l^\g(\a)=(\l^\g_1(\a),\l^\g_2(\a))$, where 
$$ \l^\g_i(\a)= \g^d \sum_{x \in \L_\g} \a_i(x)\d_{\g x},\,\,\,i=1,2,  \qquad 
\l^\g =\l^\g_1(\a)+\l^\g_2(\a)=
\g^d \sum_{x \in \L_\g} \d_{\g x}.
\Eq(randmeas)
$$
   \medskip \noindent {\bf Proof of  Theorem \equ(2th-g)}  We
follow the general scheme introduced in [KL] chap. 4.  
We have to show:
\item{(i)} For   any $\a$, the sequence $(\mus Q^{\g,\a}_{\s^\g})_{\g}$ is tight.
\item{(ii)}  Any  limit point $\mus Q^{\a}$ of  $(\mus Q^{\g,\a}_{\s^\g})_{\g}$ is   $\P$-a.s.  concentrated on measures $(\mus \pi_t)_{t\in[0,T]} \in \CC([0,T],  \MM^{ac}_{p_1} \times   \MM^{ac}_{p_2} )$. 
\item{(iii)} For $\P$-a.s. $\a$, any limit point $\mus Q^{\a}$ of  $(\mus Q^{\g,\a}_{\s^\g})_{\g}$ 
is concentrated on trajectories $(\mus \pi_t)_{t\in[0,T]}$ such that $\mus \pi_t({d}r)=\mus m(t,r){d}r$, where the density $\mus m$ is a  weak solution  of   \equ(G2b).
\item{(iv)} Equation  \equ(G2b)  has a unique weak solution. 

\bigskip
For    (ii) we use that the spins are finite-valued  (cf. [KL]).   
Namely, fix  $ G  \in \CC(\L)$, 
$$  \sup_{ 0 \le t \le T }  \left |    <\pi_{i,t}^{\g},G   >  \right | \le   \g^d   \sum_{x\in \L_\g}       |G (\g x)| \a_i(x) ,\quad i=1,2,$$
because there is at the most one spin per site and   $ \a_i(x) \ge 0$.   As in the case without random field, the application  $(\pi_{i,t})_{t\in[0,T]}\mapsto\sup_{t\in [0,T]}  <  
  \pi_{i,t}, G>$ is continuous in the weak topology. Thus by weak convergence and Lemma \eqv   (erg) (by the independence of the r.v. $\a$'s)
  all limits points are concentrated on trajectories $(\pi_{i,t})_{t\in[0,T]}$  such that 
 $$    \left | <   \pi_{i,t},G>  \right | \le      \int_{\L }        |G(r)|  p_i  dr, \quad \P- a.s.$$
 Point (iv) is derived similarly to the proof of the Cauchy-Lipschitz
theorem. For Points (i) and (iii), let  $ \mus G= (G_1,G_2)\in  \left(\CC^{1,0}([0,T]\times\L)\right)^2$. For $\mus\pi\in D([0,T],\MM_1 \times \MM_1)$, let 
$$
\ell_t(\mus\pi,\mus G)=<\mus\pi_t,\mus G(t,\cdot)>-<\mus\pi_0,\mus
G(0,\cdot)>-\int_0^t  < \mus \pi_s,\partial_s \mus G(s,\cdot) >\,ds. 
\Eq (linpart)
$$ 
We have, for $x\in\L_\g$,
$$\eqalign {
\LL^{\g,\a}(\s(x))&= -\s(x)+\s(x)(1-2c_x^{\g,\a}(\s)) \cr &=
-\s(x)+\tanh[ (J_\g \star\s) (x)+\th \a(x) ]. }\Eq(g34)  $$ 
The
$P^{\g,\a}_{\s^\g}$-martingale 
 $\overline N^G_\g \equiv (\mus
N^{\mus G}_\g(t))_{t\in[0,T]} $  with respect to the natural
filtration associated to  $(\s_t)_{t\in[0,T]}$ (cf.  \equ(1-alpha)) given by
 $$\eqalign {
\mus N^{\mus G}_\g (t)&=\ell_t(\mus\pi^\g,\mus G)
-\g^d \int_0^{t} 
 \sum_{i=1,2} \sum_{x\in\L_\g} G_i(s,\g x)\a_i(x) \LL^{\g,\a}(\s_s(x))
  \,ds \cr 
&=\ell_t(\mus\pi^\g,\mus G)
+\int_0^{t}< \overline{\pi}_s^\g,\mus G(s,\cdot)>ds
-\sum_{i=1,2} \int_0^t <\l_i^\g(\a),G_i(s,.)\tanh[\pi^\g_s*J+a_i\th]>  ds\; ,
}\Eq (g24) $$
  has quadratic variation 
 $$ 
 < \mus N^{\mus G}_\g,\mus N^{\mus G}_\g > (t) 
 = - 2\g^{2d} \sum_{i=1,2} \sum_{x\in\L_\g}\a_i(x) G^2_i(s,\g x)  \int_0^t\big\{-1+\s_s(x)\tanh[(J_\g \star\s_s)(x)+a_i\th]\big\}ds.   \Eq (g24qv) $$
Hence,  for any $\a \in \O$,  since $\tanh$ is a smooth function and $(J_\g \star\s)
(x)+a_i\th $ is uniformly bounded in $x$, $\s$, by Doob's inequality,
 $$  \lim_{\g \to 0}
 P^{\g,\a}_{\s^\g} \Bigl( \sup_{t \in [0,T]}\bigl |\mus
 N^{\mus G}_\g (t)\bigr | >\d \Bigr ) 
\le \lim_{\g \to 0}\frac1 {\d^2} E^{\g,\a}_{\s^\g} \Bigl ( \bigl |\mus
 N^{\mus G}_\g (T)\bigr |^2\Bigr ) 
\le \lim_{\g \to 0}\frac 1 {\d^2} C(\mus G,T) \g^{d}=0. 
 \Eq(g27) 
 $$
 Bound \equ(g27) yields Point (i), by Prohorov's
criterion. Point (iii)  will consist in identification of the limit. 
To obtain a closed form for the limiting equation, we only
need to average over the disorder,   that is to replace in the limit
$\g\to 0$ the random
discrete measures $\l_1^\g(\a)$ and $\l_2^\g(\a)$ by their expectations
$p_1\l$ and $p_2\l$ with respect to the environment.
Denote
$$
{\widetilde \ell}_t (\mus \pi^\g ,\mus G)\, =\,
\int_0^t \Big\{\big< {\mus\pi}^\g_s, {\mus G}(s,\cdot)\big>  -
\sum_{i=1}^2\big< p_i \lambda^\g,G_i (s,\cdot)   \tanh (\pi_s^\g * J+a_i\theta\big)\big>
\Big\}\,  ds.  \Eq (elltilde1)
$$
Putting together  \eqv (g24), \equ(g27), 
 and applying Lemma \eqv (unifbound)  below,  we get that for all subsequences    $$\eqalign{& 
\liminf_{k\to+\infty}\,  Q^{\g_k,\a}_{\s^{\g_k}}\Bigl(\sup_{t\leq T}\,\,\Bigl|{    \ell_t (\mus \pi^{\g_k} ,\mus G) + \widetilde \ell}_t (\mus \pi^{\g_k},\mus G)  
 \Bigr|>\frac{\d}{2}\Bigr)=0.
}\Eq (g8bis) $$
 Denoting 
  $  m= m_1+m_2$,  for $ \mus m = (m_1,m_2)$, see \eqv (G2b),   this gives,  for almost any $\a$,
   $$   \eqalign { &
  Q^{\a}\Bigl(\sup_{t\leq T}\,\,\Bigl|\sum_{i=1}^2 \int_\L  \left [ \left\{G_i(t,r) m_i(t,r)- G_i(0,r) m_i(0,r)\right\}   -
    \int_0^{t} \partial_sG_i(s,r) m_i(s,r)\,ds \right ] dr \cr &
   +  \sum_{i=1,2}  \int_0^{t}\int_\L  \left [ 
G_i(s,r) m_i(s,r)  - \E(\a_i(0))  \tanh[(J * m)  (s,r)+a_i\th]  \right ]\, drds
  \Bigr|>\frac{\d}{2}\Bigr)\cr &= 
  Q^{\a}\Bigl(\sup_{t\leq T}\,\,\Bigl|   \ell_t (\mus m,\mus G) + {\widetilde \ell}_t (\mus m,\mus G)   \Bigr|>\frac{\d}{2}\Bigr)
  =0,  } 
\Eq (g88) $$
where we set by an abuse of notation 
$$
\ell_t(\mus m,\mus G)=< \mus m(t,\cdot) ,\mus G(t,\cdot)>-< \mus m(0,\cdot) ,\mus
G(0,\cdot)>-\int_0^t  <  \mus m(s,\cdot) ,\partial_s \mus G(s,\cdot) >d s, 
\Eq (linparta)
$$ 
$$
{\widetilde \ell}_t (\mus m,\mus G)\, =\,
\int_0^t \Big\{\big<  \mus m(s,\cdot) , {\mus G}(s,\cdot)\big>  -
\sum_{i=1}^2\big< p_i \lambda,G_i (s,\cdot)   \tanh ( J*m(s, \cdot)+a_i\theta\big)\big>
\Big\}  ds.  \Eq (elltilde1a)
$$
 This leads to identification of the limit (iii), that is to equation  \eqv (G2b).
\qed

 \medskip

\noindent{\bf \Lemma (unifbound)} {\it 
For $i=1,2$, $  G_i \in \CC^{1,0}([0,T]\times\L)$,   there exists a positive
function $\e$ on $\R_+$ with $\lim_{s\to 0}\e(s)=0$ 
such that for all $l\in\N\setminus\{0\}$, $\d>0$, and  $\EE_i(\d, l,\g,\a)$ defined in \eqv  (ju1), the quantity
$$
\Delta^\g_i(\a,\s,T)=\int_0^T\Bigl|<\l_i^\g(\a)-p_i \l,G_i(s,.)\tanh[\pi_s^\g*J+a_i\th]>\Bigr|\,ds\Eq (def:Delta)
$$
satisfies
$$
\Delta_i^\g(\a,\s,T)\le      \frac{\d}2  T \|  G_i (s,\cdot)\|_1  +
2  T \|G_i(s,\cdot) \|_\infty   
   \EE_i (\d, l,\g,\a)
+\e(\g l) T +\e(\g) T,  \Eq (PA1)
$$
$$ \lim_{\g \to 0}   
\Delta_i^\g(\a,\s,T) =0 , \qquad \P-a.s.  \Eq (par15) $$}
\smallskip

\proof  We introduce averages over large
microscopic boxes of size $l$ but small w.r.t. the range $\g^{-1}$ of
the interaction ($l$ will go to infinity but after the limit $\g\to 0$). 
To keep notation readable, the function $\e$ may
vary from one line to another but keeping the same property that $\lim_{s\to 0}\e(s)=0$.
Since $J$ and  $\tanh$ are uniformly Lipschitz, and $G_i$ is uniformly continuous (in space), there are
constants $c_1>0$, $c_2>0$ such that, see \eqv (microav), for $x\in\L_\g$,  
$$
\sup_{\s\in\SS_\g}\left|(\pi^\g(\s)* J)(\g x)
-\bigl((\pi^\g(\s)* J)(\g .)\bigr)^{(l)}(x)\right|\le c_1\g l,
$$
$$
\sup_{\s\in\SS_\g}\left|\tanh[\beta((\pi^\g(\s)*J)(\g x)+a_i \th)]
-\bigl(\tanh[(\pi^\g(\s)*J)(\g .)+a_i\th]\bigr)^{(l)}(x)\right|\le c_2\g l.
$$
$$
\sup_{0\le s\le T}\left|G_i(s,\g x)
-\bigl(G_i(s,\g .)\bigr)^{(l)}(x)\right|\le \e(\g l).
$$
Recalling   notation 
 \equ(randmeas),  by summation by parts  we get
$$
\int_0^T\!\left|<\l_i^{\g}(\a)-\l_i^{\g}(\a^{(l)}),G_i(s,.)\tanh[\pi^\g_s*J+a_i\th]>\right|\,ds\le \e(\g l) T
$$
 and by    uniform continuity or
Lipschitz condition  
$$
\int_0^T\!\left|<p_i\l^{\g}-p_i\l,G_i(s,.)\tanh[\pi^\g_s*J+a_i\th]>\right|\,ds\le \e(\g)T.
$$
Therefore, we have 
$$
\Delta_i^\g(\a,\s,T)\le \e(\g
l) T+\e(\g)T+\int_0^T\!\left|<\l_i^\g(\a^{(l)})-p_i\l^\g,
G_i(s,\cdot)\tanh[\pi^\g_s*J+a_i\th]>\right|\,ds.  
$$
To derive \eqv (PA1),  we take into account definitions \eqv (gAlbis), \eqv  (ju1), and  that  $ | \a_i^{(l)}(x)-p_i |\le 2$,
 to write
$$\eqalign{ &| <\l_i^\g(\a^{(l)})-p_i\l^\g,
G_i(s,\cdot)\tanh[\pi^\g_s*J+a_i\th]> | = \cr &
\left|\g^d\sum_{x\in\L_\g}
G_i(s,\g x)
\tanh[(J_\g \star \s_s)(x)+a_i\th]\left[\1_{A_{l,\d}(x,i)}(\a)+\1_{A^c_{l,\d}(x,i)}(\a)\right]
\left( \a_i^{(l)}(x)-\E(\a_i(x))\right )\right|\cr &\le 
\g^d\sum_{x\in\L_\g}\left|G_i(s,\g x)
\tanh[(J_\g \star \s_s)(x)+a_i\th]   \right |\left[ \d+2\1_{ A^c_{l,\d}(x,i)}(\a) \right]\cr 
 &   \le  \d   \|  G_i (s,\cdot)\|_1  +
2  \|G_i(s,\cdot) \|_\infty \EE_i(\d,l,\g,\a).   
}
$$
 Applying   Lemma \equ(erg1)   to \eqv (PA1)  we get \eqv (par15).  \qed 
\medskip\smallskip
  
 {\bf Proof of Corollary \eqv (ch-g):}    First notice that applying Lebesgue dominated convergence Theorem
in the time integral, Theorem \eqv(2th-g) implies that for any $G_i\in \CC^{0,1}([0,T]\times\L) $
we have
$$
\lim_{\g\to 0} P^{\g,\a}_{\s^\g}\left
[\int_0^T\left|<\pi_{i,s}^\g,G_i(s,\cdot)>-<m_i(s,\cdot),G_i(s,\cdot)> \right|\,ds \geq\d\right ] =0.  
\Eq(time_int)
$$
Now remark that  integrating in time \equ(G2b), 
$$\eqalign { & 
<m_i(t,\cdot),G_i(t,\cdot)>=<m_i(0,\cdot),G_i(0,\cdot)>+\int_0^t\!<m_i(s,\cdot), \partial_s G_i (s,\cdot)-G_i(s,\cdot)>\,ds  \cr &
 +p_i\int_0^t\!< \tanh[(J * m)(s,\cdot) +a_i\th], G_i(s,\cdot)>\,ds .}
$$
Introducing the martingale $\mus N^G_\g$, see 
\equ(g24), and using \equ(g27)
 we get
$$
 P^{\g,\a}_{\s^\g}\left [\sup_{t\in[0,T]}\left|
<\pi_{i,t}^\g,G_i(t,\cdot)>)-<m_i(t,\cdot),G_i(t,\cdot)>
\right| \geq\d\right ] \le A_\g+B_\g+C_\g+\e(\g),
$$
with $\lim_{\g\to 0}\e(\g)=0$ and
$$
\eqalign{
A_\g&= P^{\g,\a}_{\s^\g}\left [\int_0^T\!\bigl|
<\pi_{i,s}^\g-m_i(s,\cdot), \partial_s G_i   (s ,\cdot)-G_i(s,\cdot)>
\bigr|\,ds\ge\frac{\d}4\right ], \cr
C_\g&= P^{\g,\a}_{\s^\g}\left [\bigl|
<\pi_{i,0}^\g-m_i(0,\cdot),G_i(0,\cdot)>
\bigr|\ge\frac{\d}4\right ],  } $$
$$ \eqalign {  
B_\g&= P^{\g,\a}_{\s^\g}\left [\int_0^T\!\bigl|
<\l_i^\g(\a), G_i   (s ,\cdot)\tanh[\pi^\g_s*J+a_i\th]>-<p_i\l,G_i   (s ,\cdot)\tanh[(J* m) (s,\cdot)+a_i\th]>
\bigr|\,ds\ge\frac{\d}4\right ] \cr &\le  P^{\g,\a}_{\s^\g}\left [\int_0^T\!\left|
<\l_i^\g(\a) -p_i\l ,G_i   (s ,\cdot) \tanh[\pi^\g_s*J+a_i\th]>\right|\,ds\ge\frac{\d}8\right ] \cr
\ & \ \ +
P^{\g,\a}_{\s^\g}\left [\int_0^T\!\left|
<p_i\l, G_i   (s ,\cdot)
\bigl(\tanh[\pi^\g_s*J+a_i\th] -\tanh[(J* m_i) (s,\cdot)+a_i\th]\bigr) 
 >\right|\,ds\ge\frac{\d}8\right ].
} \Eq (paris15a)
$${}From  \equ(hyp-G1a) and \equ(time_int), $\displaystyle{\lim_{\g \to 0}A_\g =\lim_{\g \to 0}C_\g =0}$. 
{}For $B_\g$, from Lemma \equ(unifbound), the limit when $\g \to 0$ of the first term in the right hand side of \eqv   (paris15a) is equal to zero; the second term vanishes from \equ(time_int) since the function $\tanh$ is Lipschitz continuous.
\qed

\bigskip
\chap{5.  The perturbed dynamics  and Radon-Nikodym  derivative }5
 \numsec= 5
 \numfor= 1
 \numtheo=1

The general strategy to derive the large deviation principle prescribes to find a family of mean one positive martingales that can be expressed as functions of the empirical measures. Following  [DV], the relevant martingales are obtained as Markovian  perturbations of the original process. In this section we define a class of time dependent, random external potentials, {\it the perturbations}, to which we can associate a  trajectory  $(m(t,\cdot))_{t \in [0,T]}$ smooth in time. We  show the law of large numbers for the empirical measures of the dynamics associated to these perturbations and derive the Radon-Nikodym derivative of the perturbed process with respect to the unperturbed one. 
 
Given a realization  $\a$ of the magnetic field, $ \mus V=(V_1,V_2) \in\left( \CC^{1,0} ([0,T]\times \L)\right)^2$,  let
     $$ V(t, \g x,\a(x))
 = \sum_{i=1,2} \a_i(x) V_i(t, \g x) \Eq (n2d)
 $$ 
 be the full external random perturbation for the magnetization trajectories $ \pi^\g (\s)$ (not colored). As pointed out in the introduction this perturbation strongly depends on the randomness. It is therefore convenient to consider a Glauber evolution for the colored particle system, associated to the Hamiltonian obtained by summing up  \eqv (g1) and 
 $$ 
 H^{\mus V,\g,\a}(\s)= - \sum_{x\in \L_\g} \sum_{i=1,2} \a_i(x) V_i(t, \g x) \s(x) .  \Eq(n2) $$
To this aim we define time dependent rates, 
for all $ x \in \L_\g$,  $\s \in \SS_\g$,
 $$ c_x^{\mus V,\g,\a}(\s, t) =  e^{- \s (x) 2 V (t,\g x,\a(x)) } c_x^{\g,\a}   (\s) =  
 \frac { e^{- \s (x)[(J_\g \star \s)(x)+ \th \a(x) + 2 V(t, \g x, \a(x))]}}  {2 \cosh [( J_\g \star \s)(x)+ \th \a(x)] }.  \Eq (P.2a)$$ 
 Assume    $ (\s^\g)_\g, m_0$    satisfy    \eqv (hyp-G1a). We denote by $ P_{\s^\g}^{\mus V,\g, \a}$ the law (and by $ E_{\s^\g}^{\mus V,\g, \a}$ the expectation)  of the corresponding inhomogeneous Glauber process $(\s_t)_{t\in[0,T]}$ on  $ \SS_\g$, that is the unique probability measure on $D([0,T], \SS_\g)$ with initial condition $\s^\g$ under which 
 $\displaystyle{ f (\s_t) - f (\s_0) -\int_0^t {\cal L}^{\mus V,\g,\a}_{s} (f) (\s_s) ds}$
 is a martingale w.r.t. the canonical filtration, for all cylinder function $f$, where
  $${\cal L}^{\mus V,\g,\a}_{s} (f) (\s_s) = \sum_{x\in\L_\g} c_x^{\mus V,\g,\a}(\s_s, s) [ f(\s_s^x)-f(\s_s)] .\Eq(L-pert1)$$
Let $ \mus Q^{\mus V,\g,\a}_{\s^\g}$ be the law of the corresponding empirical measures.  
\medskip

\noindent {\bf \Theorem (pert1)} 
{\it   
 Assume $ (\s^\g)_\g, m_0$ satisfy    \eqv (hyp-G1a).
 For all $t \in[0,T]$, $ \mus G=(G_1,G_2)\in \left(\CC^1 (\L)\right)^2$, and $\d>0$,   
$$
\lim_{\g\to 0}  \mus Q^{\mus V,\g,\a}_{ \s^\g}\left [ \left|  <   \mus \pi_t^\g,\mus G> - <  \mus m^{\mus V} (t,\cdot),   \mus G> \right| \geq\d \right ] =0, \qquad \P-a.s.,
$$
where $ \mus m^{\mus V} = (m^{\mus V}_1,m^{\mus V}_2)$ is the solution of, for $i=1,2$, 
$$ 
\left \{ \eqalign {& \partial_t m_i (t,r)=  \left \{ -m_i(t,r)+ p_i   \tanh \left [(J * m )(t,r)+a_i\th +2V_i(t,r) \right] \right\} \frac {\cosh  \left[( J * m)(t,r)+a_i \th +2V_i(t,r) \right]} {\cosh  \left[ (J * m)  (t,r)+a_i\th \right] }, 
 \cr & m_i(0,\cdot)=p_i m_0(\cdot),\quad m= m_1+m_2. }\right. \Eq (Eq.1)
$$
 }

\Remark (R-pert1) For existence and uniqueness of the solution $\mus m^{\mus V}\in\left( \CC([0,T], L^\infty (\L))\right)^2$, we refer to Remark \eqv (R-2th-g). Notice that 
the set $\{ \mus m\in(L^\infty (\L))^2:\|m_i\|_{\infty}\le p_i, i=1,2 \}$ is still left invariant.

\smallskip
 \proof  
We proceed as for Theorem \eqv (2th-g). We use 
 $$ 
 1 = \frac {(1-\s_s(x))} 2 + \frac {(1+\s_s(x))} 2 = \1_{\{\s_s(x)=-1\}} + \1_{\{\s_s(x)=1\}}.\Eq (1-sigma) 
 $$
 For  $i\in\{1,2\}$  we have
 $$\eqalign {&
{\cal L}^{\mus V,\g,\a}_{s} (\a_i(x)\s(x)) = -2\a_i(x)\s(x)c_x^{\mus V,\g,\a}(\s, s)  \cr &= - \a_i(x)   e^{- \s (x) 2 V_i (s,\g x) }2\s(x)  
 \frac { e^{- \s (x)[(J_\g \star \s)(x)+ a_i\th ]}}  {2 \cosh [ J_\g \star \s)(x)+a_i \th ] }
  \cr & 
   = - \a_i(x)   \left[ \frac{(\s(x)+1)}2   
 \frac { e^{- [2 V_i (s,\g x)+(J_\g \star \s)(x)+ a_i\th ]}}  { \cosh [ (J_\g \star \s)(x)+a_i \th ] }
 +\frac{(\s(x)-1)}2 
 \frac { e^{[ 2 V_i (s,\g x)+(J_\g \star \s)(x)+ a_i\th ]}}  { \cosh [ (J_\g \star \s)(x)+a_i \th ] }\right]
  \cr &
= - \a_i(x)\s(x) \frac { \cosh[ 2 V_i (s,\g x)+(J_\g \star \s)(x)+ a_i\th ]}  { \cosh [ (J_\g \star \s)(x)
+a_i \th ] }  
- \a_i(x)    
 \frac { \sinh [2 V_i (s,\g x)+(J_\g \star \s)(x)+a_i \th ]}  { \cosh [ (J_\g \star \s)(x)+ a_i\th ] }.}\Eq(g34-pert1)$$
\qed
\bigskip
We have the analogous result to Corollary \eqv (ch-g):
\medskip \noindent
 {\bf \Corollary (ch-g1)} {\it 
  For all $\mus G=(G_1,G_2)\in
\left(\CC^1(\L)\right)^2$,     and $\d>0$,
$$
\lim_{\g\to 0} \mus Q^{\mus V,\g,\a}_{\s^\g}\left [\sup_{t\in[0,T]}\left|
<\mus \pi^\g_t, \mus G>-<\mus m^{\mus V}(t,\cdot),\mus G>
\right| \geq\d\right ] =0. 
$$}

\medskip
\noindent {\bf\Theorem   (RN)}  {\it  Let   $ \mus V=(V_1,V_2) \in \left(\CC^{1,0} ([0,T] \times  \L)\right)^2$.  
 The Radon-Nikodym derivative is given by 
 $$
  \frac{ d P_{\s^\g}^{\mus V,\g,\a} } {d  P_{\s^\g}^{\g,\a}} ( 
 \s_{[0,T]}) =
\exp\Bigl\{ \g^{-d}\Bigl(\ell_T(\mus \pi^\g (\s) ,\mus V)
 -   \frac 12  \int_0^T   F_{\mus V(s)}(\mus\l^\g(\a),\mus \pi^\g_s) \, ds\Bigr)   \Bigr\},  
 \Eq (enza2b) 
 $$
where  $\ell_T$  was defined  in \equ(linpart),  
  $\mus\l^\g(\a)$ in \equ(randmeas),  $F_{\mus V(s)}(\cdot,\cdot)$ in \eqv (rate1), and we have abbreviated 
  $\s_{[0,T]}=(\s_t)_{t\in[0,T]}$.}
  
\smallskip
 \proof 
The Radon-Nikodym derivative associated with rates \equ(P.2a) is given by
 (see [HS] or  [KL], Appendix 1, Proposition 7.3)  
 $$\eqalign {
 \frac{ d P_{\s^\g}^{\mus V,\g,\a} } {d  P_{\s^\g}^{\g,\a}} ( \s_{[0,t]} ) &=\exp\Big\{ - H^{\mus V,\g,\a}(\s_t) + H^{\mus V,\g,\a}(\s_0)\cr &-\int_0^t  \exp\Big\{ H^{\mus V,\g,\a}(\s_s) \Big\} \big( \partial_s +\LL_\g\big)\exp\Big\{ -  H^{\mus V,\g,\a}(\s_s) \Big\}\, ds\Big\}\cr & =\exp\Big\{ \ell_t(\mus \pi^\g(\s),\mus V)  
- \int_0^t\sum_{x\in\L_\g} c_x^{\g,\a}(\s_s) \left [ e^{-2 \s_s (x)V(s,\g x,\a(x)) } -1 \right ] 
\,ds\Big\}, 
 } 
 $$
 because of \equ(linpart), \equ(n2). To get \equ(enza2b), 
 we use trigonometric formulas
 to write   (remember  \equ(g2), \equ(1-alpha), \equ(g34-pert1)) 
 $$ \eqalign {   
 2\g^d\sum_{x\in\L_\g} & c_x^{\g,\a}(\s_s) \left [ e^{-2 \s_s (x)V(s,\g x,\a(x)) } -1 \right ] \cr= 
 & 2\g^d\sum_{x\in\L_\g} \left \{ \frac {(1-\s_s(x))} 2
 + \frac {(1+\s_s(x))} 2\right \}\left \{ \sum_{i=1,2} \a_i(x) \right \}c_x^{\g,\a}(\s_s) \left [ e^{-2 \s_s (x)V(s,\g x,\a(x)) } -1 \right ] \cr 
 = & \g^d \sum_{i=1,2}  \sum_{x\in\L_\g}  \frac {(1-\s_s(x))}2
 \a_i(x) \frac{\exp[ (J_\g \star \s_s)(x) +a_i\th ]}{ \cosh [ (J_\g \star \s_s)(x)+a_i\th]} \left [ e^{2 V_i(s,\g x) } -1 \right ] \cr 
  & +\g^d \sum_{i=1,2}   \sum_{x\in\L_\g}  \frac {(1+\s_s(x))} 2 \a_i(x)\frac{\exp[-(J_\g \star \s_s)(x) -a_i \th ]}{ \cosh [ (J_\g \star \s_s)(x)+a_i\th]} \left [ e^{-2 V_i(s,\g x) } -1 \right ] 
\cr = & \g^d \sum_{i=1,2}  \sum_{x\in\L_\g} \a_i(x)\left \{ \cosh[2 V_i(s,\g x)]-1 +\tanh[ ( J_\g \star \s_s)(x)+a_i\th]\sinh[ 2V_i(s,\g x)]\right \}\cr 
     &- \g^d \sum_{i=1,2}  \sum_{x\in\L_\g}  \a_i(x)\s_s(x)\left \{\tanh[ ( J_\g \star \s_s)(x)+a_i\th](\cosh[ 2 V_i(s,\g x)]-1)+\sinh[ 2 V_i(s,\g x)]\right \}\cr 
 =& F_{\mus V(s)}(\mus\l^\g(\a),\mus \pi^\g_s).
  }
$$  
  \qed 
\medskip\smallskip
Note that the Radon-Nikodym derivative depends on the randomness through $\mus \pi^\g$ and $\mus\l^\g(\a)$. By next proposition, which is proved in Appendix B, we can
 replace   $\mus\l^\g(\a)$ in  $ F_{\mus V(s)}(\cdot,\mus \pi^\g_s)$  with $(p_1\l^\g, p_2\l^\g)$, 
   making an error which goes uniformly (for all $\s\in \SS_\g$  and $\P$-a.s.) to zero as $ \g\to 0$. 
  
\medskip
\noindent {\bf \Proposition    (RN1)}  {\it  Let   $ \mus V=(V_1,V_2) \in \left(\CC^{1,0} ([0,T] \times  \L)\right)^2$. 
  There exists a positive
function $\e$ on $\R_+$ with $\lim_{s\to 0}\e(s)=0$ such that for any $\d>0$, $l\in\N\setminus\{0\}$,  we have 
 $$
 \left|  \int_0^T \Bigl[
F_{\mus V(s)}(\mus\l^\g(\a),\mus\pi^\g_s)
-F_{\mus V(s)}( (p_1\l^\g, p_2\l^\g), \mus\pi^\g_s)\Bigr]
\, ds\right|  
\le  \e(\g l) T +
   T  C(V_1,V_2)[\d +  \sum_{i=1,2}\EE_i(\d, l,\g,\a)]
$$
where the positive constant $C(V_1,V_2)$ depends on the $L^\infty$ norm of $(V_1,V_2)$.
    }

 \bigskip
 \chap{ 6. Upper Bound  }6
\numsec=6
\numfor= 1
\numtheo=1

In this section we investigate the upper bound of the large deviation principle for compact sets and then closed sets
of the topological space $D([0,T),\MM_1 \times \MM_1)$. Notice that in [C] the result was stated for closed sets
in $\CC([0,T),\MM_1 \times \MM_1)$.
We first prove exponential tightness, so that it is enough to derive the upper bound  of the large
deviation principle for compact subsets. The strategy then follows the martingale approach  introduced by [DV]:
we need to show that trajectories which are not absolutely
continuous with respect to the Lebesgue measure and not absolutely continuous in time can be neglected in the large deviations regime. 
To exclude these ``bad" paths, as in [FLM], we introduce an energy functional via an exponential martingale.
With this we prove an upper bound with an auxiliary rate functional which is infinite on the set of bad trajectories.

  \medskip
\noindent {\bf \Proposition (exptight)}  {\it  For any $ \ell \ge 1 $, there exists  
 a compact subset  $\mus  K_\ell\subset D ([0,T], \MM_1 \times \MM_1) $  such that for any $\s_\g \in \SS_\g $,
  $$
 \limsup_{\g \to 0} \g^d \log   \mus  Q^{\g, \a}_{\s^\g} (  \mus K_\ell ^c ) \le  -\ell .  
$$
} 

The proof is standard, however the main lines are recalled in Appendix B.
\bigskip
 
\noindent 
For $\mus \pi\in D\big( [0,T],\MM_1\times \MM_1\big)$, $\mus G=(G_1,G_2) \in (\CC^{1,0}([0,T]\times \L))^2$
denote 
$$\eqalign{
\J_{\mus G}(\mus \pi)&\, =\, {\ell}_T(\mus \pi ,\mus G)\, +\,{\widetilde \ell}_T(\mus \pi ,\mus G)\cr
\ &\ \
\, - 2 \sum_{i=1}^2\int_0^T \Big\{ \big< p_i\lambda, G_i^2(s,\cdot)\big> 
    -\big<\pi_{i,s},G_i^2(s,\cdot)\tanh (\pi_s * J+a_i\theta\big)\big>\Big\} \, ds\, ,
}\Eq(j-functa)
$$
for $\ell_T$, $\tilde \ell_T$ given  in \equ(linparta),   \eqv (elltilde1a).
 We define  the auxiliary rate functional $\JJ : D\big( [0,T],\MM_1\times \MM_1\big)\to \mus\R$  as 
$$
\JJ(\mus \pi)=\left \{ \eqalign { 
\sup_{\mus G\in (\CC^{1,0}([0,T]\times \L))^2}\Big( \J_{\mus G}(\mus \pi)\Big)&\qquad  \hbox{ if }\ \mus\pi \in D([0,T], \MM_1^{ac}\times\MM_1^{ac} )\, ,\cr
+\infty \qquad &\qquad \hbox{ otherwise }.
}\right.
\Eq(j-funct)
$$

\medskip
\noindent {\bf \Lemma (up3)} {\it For all $\mus\pi \in D\big( [0,T],\MM_1\times \MM_1\big)$, if
$\JJ(\mus \pi )<\infty$, then $\mus\pi\in \CC\big( [0,T],\MM_1^{ac}\times \MM_1^{ac}\big)$.
}

\smallskip
\proof
Fix $\mus\pi \in D\big( [0,T],\MM_1\times \MM_1\big)$ such that $\JJ(\mus \pi )<\infty$.  By  definition of $\JJ(\cdot)$,
$\mus\pi \in D\big( [0,T],\MM_1^{ac}\times \MM_1^{ac}\big)$.
Let $\mus g=(g_1,g_2)\in ( \CC(\L))^2$ and $0 \leq s < t < T$. For each $\delta>0$, let $\psi^{\delta}_{s,t}:[0,T]\to\R$ be the
function given by
$$
\psi^{\delta}_{s,t} (\tau) = 
\left \{ \eqalign { & 
0 \qquad \hbox{ if }\qquad  0\leq \tau\leq s \;\hbox{ or }\; t \leq \tau\leq T\, , \cr
&\frac{\tau-s}{\delta}  \qquad \hbox{ if }\qquad s\leq \tau \leq s+\delta\, , \cr
&1 \qquad \hbox{ if } \qquad s+\delta\leq \tau\leq t-\delta\, , \cr
&\frac{t-\tau}{\delta} \qquad \hbox{ if } t-\delta\leq \tau\leq t \, .
}\right. \Eq (psi1)
$$
Denote $\mus G^{\delta}(\tau,r) = \psi^{\delta}_{s,t}(\tau)\mus g(r)$. 
Since $ \mus G^{\delta} $ can be approximated by functions in $(\CC^{1,0} ( [0,T]\times \L))^2$,
considering $\displaystyle{\frac{\mus G^{\delta}}{\sqrt{t-s}}}$ as a test function and performing the limit $\delta \to 0$, we obtain
$$\eqalign{
&\sqrt{t-s} \; \lim_{\delta\to 0} \J_{\frac{\mus G^{\delta}}{\sqrt{t-s}}}(\mus \pi) \, = \, 
< \mus\pi_t,\mus g>-< \mus\pi_s,\mus g>\cr
&\qquad\qquad\qquad
\, +\,  \int_s^t \Big\{< {\mus\pi}_{\tau}, {\mus g}>  -
\sum_{i=1}^2< p_i \lambda,g_i  \tanh (\pi_\tau * J+a_i\theta\big)>
\Big\}   d\tau\, \cr
&\qquad\qquad\qquad   - 2
\frac{1}{\sqrt{t-s}}\sum_{i=1}^2\int_s^t \Big\{< p_i\lambda, g_i^2 >
    -< \pi_{i, \tau},g_i^2 \tanh (\pi_\tau* J+a_i\theta\big)>\Big\} \, d\tau\,.
}\Eq(stbound)
$$
Since 
$\displaystyle{\sqrt{t-s} \; \lim_{\delta\to 0} \J_{\frac{\mus G^{\delta}}{\sqrt{t-s}}}(\mus \pi) \le \sqrt{t-s}\; \JJ(\mus \pi)}$, we get
$$\eqalign{
\Big|< \mus\pi_t,\mus g >-< \mus\pi_s,\mus g > \Big|&
\le  \;C_0 (t-s) \sum_{i=1}^2\Big\{\| g_i\|_1 +\frac{1}{\sqrt{t-s}}\| g_i\|_2^2\Big\} \,+\, \sqrt{t-s}\JJ(\mus \pi)  \cr
\ &\ \ =C_0\, (t-s) \sum_{i=1}^2\| g_i\|_1 \, +\    \sqrt{t-s} \Big\{C_0  \sum_{i=1}^2\|g_i\|_2^2 \, +\,\JJ(\mus \pi) \Big\} \,,
}
$$
for some positive constant $C_0$.
This  implies that $\mus \pi \in \CC\big( [0,T],\MM_1^{ac}\times \MM_1^{ac}\big)$.
 \qed
 \medskip
  \smallskip \noindent   To prove next Lemma, we   will use  the   following    characterization of absolutely continuous functions, see [DS].
\medskip
\noindent {\bf \Proposition (P100)}
 { \it 
A function  $ \mus \phi$ belongs to $\ac([0,T], B_{1,1})$ if and only if:
 for all $\e>0$, there exists $\D>0$ such that for all integer $k>0$, rectangles $A_1, \dots, A_k$ of $\L$ and
$\{(s_i,t_i) ,\, 1\le i\le k\}$ nonempty
disjoint intervals  of $[0,T]$,
$$
\sum_{i=1}^k |t_i-s_{i}|  \l (A_i)<\D \; \Rightarrow \; \sum_{i=1}^k\left| \int_{A_i} \left (\phi_j (t_i,r)-\phi_j  (s_{i},r)\right ) \,dr  \right|<\e\, , \ \ j=1,2\, .
$$
}
\medskip   

\noindent {\bf \Lemma (up5)} {\it Let $  \mus\pi = (\phi_1(s,r)dr, \phi_2(s,r)dr) \in D\big( [0,T],\MM_1^{ac}\times \MM_1^{ac}\big)$ such that $\JJ(\mus \pi )<\infty$,  then \par
\noindent{(a)}
for $i=1,2$, 
$F_i(s,r):=\big[p_i -\phi_i(s,r)\tanh (\pi_s(r) * J+a_i\theta\big)\big]\ge 0$ for almost all $(s,r)\in [0,T]\times \L$,\par
\noindent{(b)} $(\phi_1,\phi_2)\in \AA\CC\big( [0,T],B_{1,1}\big)$.
 }

\smallskip
\proof  
\noindent{(a)}
Taking as a test function $AG(\cdot,\cdot)$,  for all 
$\mus G=(G_1,G_2) \in (\CC^{1,0}([0,T]\times \L))^2$ and $A>0$, we obtain from Definition \equ(j-funct) of the rate function $\JJ$,
$$  
-\sum_{i=1}^{2}\int_0^T\int_\L  G_i^2(s,r)\, F_i(s,r) \, dr \,ds  \,\le\,
-\frac{1}{2A} \Big\{\ell_T (\mus \pi ,\mus G)\,+\, {\widetilde \ell}_T (\mus \pi ,\mus G)\Big\}\, +\, \frac{1}{2A^2} \JJ(\mus \pi )\; . \Eq(mparis2)
$$
Letting $A\uparrow \infty$, we get
$$
\sum_{i=1}^{2}\int_0^T\int_\L  G_i^2(s,r)\, F_i(s,r) dr ds  \,\ge\, 0\,.
$$
Since $\mus G$ is arbitrary, we conclude that $F_i(s,r) \ge 0$ for $1\le i\le 2$, almost everywhere.

\noindent{(b)}
We show the absolute continuity in time   for $\phi_1$.  The proof for $\phi_2 $ is similar. 
We apply the  characterization of $\AA\CC\big( [0,T],B_{1,1}\big)$ given in Proposition \eqv (P100). 
For all  positive integer $k$ let  $\{A_i\ ,\  1\le i\le k\} $  be rectangles of $\L$  and 
$\{(s_i,t_i)\, ,\ \ 1\le i\le k\}$  be nonempty
disjoint intervals  of $[0,T]$. For $i=1, \dots, k$,  denote
$$
\eta_i:=\hbox{\rm {{\rm sgn}}}\Big(\int_{A_i} \left [ \phi_1(t_i, r)- \phi_1(s_i, r)\right ]\, dr\Big)\, ,
$$
For each $1\le i\le k$, $\displaystyle 0< \delta < \frac{1}{4}\min_{1\le i\le k}\big( t_i-s_i\big)$, we set, see \eqv (psi1),
$$
V_1 (t,r)=  \sum_{i=1}^{k} \eta_i\times \psi^{\delta}_{s_i,t_i}(t)\times \1_{A_i} (r), \quad  V_2 (t,r)=0, \quad \mus V= (V_1,V_2).
    \Eq (EE.1) 
$$
Since $V_1$ can be approximated by functions in $\CC^{1,0} ( [0,T]\times \L)$, proceeding as in \equ(stbound), we obtain for any $b>0$, see \eqv (mparis2), 
$$\eqalign{
\sum_{i=1}^{k} \eta_i \Big\{\int_{A_i} \left [ \phi_1(t_i, r)- \phi_1(s_i, r)\right ] \,dr\Big\}&
\, \le \, - \sum_{i=1}^{k}\int_{s_i}^{t_i} < \eta_i\, \1_{A_i},\big[\phi_1(s,\cdot) -p_1\tanh (\pi_s *J+a_1\theta\big)\big]>  \, ds \cr
\ &\ \ 
+\,  2 b\, \sum_{i=1}^{k}\int_{s_i}^{t_i} \int_\L |\eta_i|\, \1_{A_i}(r) \, F_1(s,r) \, dr ds\, +\, \frac{\JJ(\mus \pi)}{b}\, .
}$$
Minimizing over $b$ yields
$$\eqalign{
\sum_{i=1}^{k}  \Big|\int_{A_i} \left [ \phi_1(t_i, r)- \phi_1(s_i, r)\right ]\,dr\Big|&
\, \le\, 2 \sum_{i=1}^{k} (t_i -s_i)\lambda\big(A_i\big)\cr
\ &\ \ \, +\, 2\sqrt{2 \JJ(\mus \pi)}\Big( \sum_{i=1}^{k}\int_{s_i}^{t_i} \int_\L |\eta_i|\, \1_{A_i}(r) \, F_1(s,r) \, dr ds \Big)^{1/2}\cr
& \ \ \ \le\,  2 \sum_{i=1}^{k} (t_i -s_i)\lambda\big(A_i\big)\, +\,
4\sqrt{\JJ(\mus \pi)}\sqrt{\sum_{i=1}^{k} (t_i -s_i)\lambda\big(A_i\big)} \, .
}\Eq(EE.12)$$
For all $\ve>0$ denote $\Delta = \min\big(\ve/4\, ,\, \ve^2/(64\JJ(\mus \pi) )\big)$. It follows from \equ(EE.12) that $\sum_{i=1}^{k} (t_i -s_i)\lambda\big(A_i\big)\le \Delta$
implies $\sum_{i=1}^{k}  \Big|\int_{A_i} \left [ \phi_1(t_i, r)- \phi_1(s_i, r)\right ]\,dr\Big|\le \ve$. This  concludes the proof.
\qed

\medskip
For  $\mus \pi \in D\big ( [0,T],\MM_1\times \MM_1\big)$,  $\ell_T$  as in \eqv (linparta) and 
 $ F_{\mus V (s,.)}$ defined   in \eqv (rate1) let
 $$
{\widehat J}_{\mus V}( \mus \pi)\, =\, \ell_T( \mus \pi,\mus V )
-\frac12 \int_0^T F_{\mus V (s,.)}\big( (p_1\lambda,p_2\lambda),\mus \pi_s\big)\,ds,  \Eq (mparis4)
$$  
$$
{\widehat J}( \mus \pi) =\sup_{ \mus V \in  (\CC^{1,0}([0,T]\times \L))^2} \Big\{ {\widehat J}_{\mus V}( \mus \pi)\Big\}\,. \Eq (mparis1) 
$$
 Remark that when $\mus\pi = (\phi_1 \l, \phi_2 \l)$, with $ \mus \phi= (\phi_1, \phi_2)  \in   \ac ([0,T], B_{1,1})$,   $\widehat J$ coincides  with the functional $I_0=J_0=J_1$ (cf. Proposition \equ(3A1)). 
 The proof of the upper bound of the large deviation principle relies on the following proposition.

\medskip\noindent {\bf \Proposition (up1)}
 { \it 
Let  $\mus K$ be a compact set of $D\big( [0,T],\MM_1\times \MM_1\big)$.  For any $0< b< 1$,
$$\limsup_{\g\to 0}
\g^d\log \mus  Q^{\g, \a}_{\s^\g} ( \mus K)  \le  -\frac{1}{1+b}\inf_{\mus\pi\in\mus  K} \left [  {\widehat J}(\mus \pi)  \; +\; b \JJ (\mus \pi)\right ]. 
$$
}

\smallskip
\proof
For $\ve>0$,  $\mu\in \MM_1$,  $g\in \CC(\L)$   denote by $\iota_\ve$ the approximation of the identity
$$
\iota_\ve (r)=\frac{1}{(2\ve)^d} \1_{\big\{\big[-\ve,\ve \big]^d \big\}}(r), \quad r \in \L
$$
and  by $\mu * \iota_\ve$ the measure defined   by 
$\big< \mu * \iota_\ve,g\big>=\big< \mu,g * \iota_\ve\big> $. 
It is  absolutely continuous with respect to the Lebesgue measure with density
$$
\frac{d(\mu * \iota_\ve)}{d\lambda}(r)\, =\, \left< \mu,\iota_\ve(r-\cdot)\right>\, ,\qquad r\in\Lambda \, .
$$
In general, we can only bound this density by $\| \iota_\ve \|_\infty$ which is of order $\ve^{-d}$. Nevertheless, in the case of the empirical measure, we have
$$
\left|\left< \pi^\g_s,\iota_\ve (r-\cdot)\right> \right|\, =\,
\left| \frac{\g^d}{(2\ve)^d}\sum_{x\, :\, \g x\in [r-\ve,r+\ve]}\s_s (x)\right|\, \le\, 1\, ,\qquad \text{for\ almost\ all\ } 0\le s\le T\, ,
$$
which means that $\pi^\g*\iota_\ve\in \MM_1^{ac}$, when $0<\g<\ve$.
Furthermore for  any ${\mus \pi}\in D\big( [0,T],\MM_1\times \MM_1\big)$, 
denote by $\mus \pi_s *\iota_\ve := \big(\pi_{1,s} *\iota_\ve\, ,\,\pi_{2,s} *\iota_\ve \big)$, $0\le s\le T$ the trajectory in 
$D\big( [0,T],\MM_1^{ac}\times \MM_1^{ac}\big)$.

Fix a function $\mus G \in (\CC^{1,0}([0,T]\times \L))^2 $. Consider the mean one   exponential martingale
$\big({\mus\ZZ}_t^{\mus G,\g}\big)_{t\ge 0}$   
$$ {\mus \ZZ}_t^{\mus G,\g}= \exp \left\{ \g^{-d}  {\mus N}_\g^{\mus G}(t) - 
\frac{\g^{-2d}}{2} \big< {\mus N}_\g^{\mus G},{\mus N}_\g^{\mus G}\big>(t) \right\}\, ,  
$$
where the martingale $\big({\mus N}_\g^{\mus G}(t)\big)_{t\ge 0}$ and its quadratic variation 
$\big(\big< {\mus N}_\g^{\mus G},{\mus N}_\g^{\mus G}\big>(t) \big)_{t\ge 0}$ were given in \equ(g24) and \equ(g24qv). 
 Using the same arguments as in Proposition \equ(RN1), by smoothness of $\mus G$ and $\pi^\g* J$, a spatial summation by parts and Taylor expansion permit
to rewrite the martingale ${\mus\ZZ}_t^{\mus G,\g}$ as  
$$
{\mus\ZZ}_t^{\mus G,\g}=\exp \Big\{\g^{-d} \J_{\mus G} (\mus \pi^\g * \iota_\ve) \,+\,\g^{-d} r(\mus G,\g,\ve,l,\d,\a) \Big\}\, ,
\Eq(martexp)$$
where $0<\g,\ve,\d<1$, $l$ is a positive integer. Here and in the sequel, $r(\mus G,\g,\ve,l,\d,\a)$  (resp. $r(\mus G,\mus V,\g,\ve,l,\d,\a)$ later on)
stands for some random variable satisfying
$$
\limsup_{\d\to 0}\limsup_{l\to \infty}\limsup_{\ve\to 0}\limsup_{\g\to 0} r(\mus G,\g,\ve,l,\d,\a) =0,\quad \P-\hbox{\rm {a-e}}\, .
\Eq(rest-upp)
$$
Let
$\mus K$ be a compact set of $D\big( [0,T],\MM_1\times \MM_1\big)$. By H\"older inequality,
$$\eqalign{
\g^{d}\log \mus  Q^{\g, \a}_{\s^\g} ( \mus K)&
=\g^{d}\log E^{Q^{\g, \a}_{\s^\g}}\Big[ \1_{\mus K}(\mus\pi^\g)
\big({\mus \ZZ}_t^{\mus G,\g}\big)^{\frac{b}{1+b}}\times \big({\mus \ZZ}_t^{\mus G,\g}\big)^{\frac{-b}{1+b}}\Big]\cr
\ & \le 
\frac{b}{1+b}\g^{d}\log E^{Q^{\g, \a}_{\s^\g}}\Big[ \1_{\mus K}(\mus\pi^\g)
{\mus \ZZ}_t^{\mus G,\g}\Big]\,+\, 
\frac{1}{1+b}\g^{d}\log E^{Q^{\g, \a}_{\s^\g}}\Big[ \1_{\mus K}(\mus\pi^\g)
\big({\mus \ZZ}_t^{\mus G,\g}\big)^{-b}\Big]\cr
\ &\ \ \le \frac{1}{1+b} \g^{d}\log E^{Q^{\g, \a}_{\s^\g}}\Big[ \1_{\mus K}(\mus\pi^\g)
\big({\mus \ZZ}_t^{\mus G,\g}\big)^{-b}\Big]\, .
}
\Eq(ub2)$$

We now exclude paths whose densities are not absolutely continuous with respect to the Lebesgue measure.
 Fix a sequence $\{F_k: k\geq 1\}$ of smooth nonnegative functions
dense in $\CC (\Lambda)$ for the uniform topology.  For
$k\ge 1$, $\varrho>0$  and $\delta>0$, let
$$
D_{k,\varrho} = \Big\{\mus \pi\in D([0,T],\MM_1\times \MM_1): \, 0\le <|\pi_{i,t}| , F_k>
\le \int_{\L} F_k(x) \, dx  \,+\, C_k \varrho
\;,\, 0\le t\le T, i=1,2 \Big\} \, ,
$$
where $C_k = C(\Vert \nabla F_k \Vert_\infty)$ is a constant depending on the
gradient $\nabla F_k$ of $F_k$. The sets $D_{k,\varrho}$, $k \ge 1$   are   closed subsets of $D([0,T],\MM_1\times \MM_1)$, as well as 
$$
E_{m,\varrho} \;=\; \bigcap_{k=1}^m D_{k,\varrho}\;, \quad m\ge 1\, .
$$
Note that the empirical measure $\mus \pi^\g$  belongs to $E_{m,\varrho}$ for $\g$ sufficiently small.
We have that $$\displaystyle{D([0,T],\MM_1^{ac} \times \MM_1^{ac}) = \cap_{n\ge 1} \cap_{m\ge 1}
E_{m,1/n}}. \Eq (mparis5) $$  
 Fix $0<b<1$. For  $\mus G,\mus V \in (\CC^{1,0}([0,T]\times\Lambda))^2$,  
$\varepsilon>0$ and   $m, n \in\Z_+$, let
${\widehat \J}_{\mus V,\mus G,\varepsilon}^{b,m,n}:D([0,T],\MM_1 \times \MM_1)\to\R\cup\{\infty\}$ be
the functional given by
$$ {\widehat \J}_{\mus V,\mus G,\varepsilon}^{b,m,n}(\mus \pi)
 =  \left \{ \eqalign {  
  {\widehat J}_{\mus V} (\mus \pi *\iota_\ve) 
+ b\J_G(\mus \pi * \iota_\ve) &\qquad \hbox{ if }\ \  \mus \pi\in E_{m,\frac{1}{n}},\cr
+\infty &\qquad \hbox{ otherwise }.
}\right. \Eq (tranq)
$$
It is lower semicontinuous because so is $\mus \pi \mapsto {\widehat J}_{\mus V} (\mus \pi *\iota_\epsilon) + b\J_G(\mus \pi * \iota_\epsilon)$,
and because $E_{m,1/n}$ is closed.

We now return to inequality \equ(ub2). 
By  Proposition \equ(RN1),  the exponential martingale $\MM^{\mus V,\g}_t$ defined by the Girsanov formula \equ(enza2b)
satisfies
$$
\MM^{\mus V,\g}_T \, :=\, \frac{ d P_{\s^\g}^{\mus V,\g,\a} } {d  P_{\s^\g}^{\g,\a}} ( 
 \s_{[0,T]})
\;=\;  \exp \Big\{ \g^{-d} {\widehat J}_{\mus V} (\mus \pi^\g *\iota_\ve) \,+\,\g^{-d} r(\mus V,\g,\ve,l,\d,\a) \Big\}.
\Eq(ub3)$$
We rewrite   \equ(ub2) as 
$$
\g^{d}\log \mus  Q^{\g, \a}_{\s^\g} ( \mus K)\le
 \frac{1}{1+b} \g^{d}\log E^{Q^{\g, \a}_{\s^\g}}\Big[ \1_{\mus K}(\mus\pi^\g)\MM^{\mus V,\g}_T\times \big(\MM^{\mus V,\g}_T\big)^{-1}
\times \big({\mus \ZZ}_t^{\mus G,\g}\big)^{-b}\Big]\, . 
$$
Since $\MM^{\mus V,\g}_T$ is a mean one positive martingale, taking into account \equ(martexp) and \equ(ub3) and
optimizing over $\mus \pi$ in $\mus K$, we obtain, for all positive integers $m,n$,
$$\eqalign{
\limsup_{\g\to 0} \g^d\log
\mus  Q^{\g, \a}_{\s^\g} ( \mus K)  &\le
\frac1{1+b}\sup _{\mus\pi\in \mus K\cap E_{m,\frac{1}{n}}} \Big\{ - {\widehat J}_{\mus V} (\mus \pi *\iota_\epsilon) 
\; -\ b\J_G(\mus \pi * \iota_\epsilon)\Big\}\cr
\ &\ \;+\;\limsup_{\g\to 0} r(\mus V, \mus G,\g,\ve,l,\d,\a)\cr
\ &= \frac1{1+b}\sup _{\mus\pi\in \mus K} \Big\{ - {\widehat \J}_{\mus V,\mus G,\varepsilon}^{b,m,n}(\mus \pi)\Big\}
\,+\,\limsup_{\g\to 0}  r(\mus V, \mus G,\g,\ve,l,\d,\a) \; .
}
$$
Optimizing the previous expression with respect to $\mus V$, $\mus G$, $\ve,l,\delta,m,n$, taking into account \eqv (rest-upp),  we get
$$
\limsup_{\g\to 0} \g^d\log
\mus  Q^{\g, \a}_{\s^\g} ( \mus K)  \le
\inf_{ \mus V, \mus G, \ve,l,\delta,m,n}\Big\{ \frac1{1+b}\sup _{\mus\pi\in \mus K} \Big\{ - {\widehat \J}_{\mus V,\mus G,\varepsilon}^{b,m,n}(\mus \pi)\Big\}  \Big\}\; .
\Eq(ub55)
$$
Since  $\mus K$ is compact and  $\displaystyle{\mus \pi \mapsto \frac1{1+b}\sup _{\mus\pi\in \mus K} \Big\{ - {\widehat \J}_{\mus V,\mus G,\varepsilon}^{b,m,n}(\mus \pi)\Big\}}$
is lower semi-continuous for all $\mus V$, $\mus G$ and $\ve,l,\delta,m,n$, we
may apply the arguments presented in [V], Lemma 11.3 to exchange
the supremum with the infimum.  In this way we obtain that the right hand side of \equ(ub55) is bounded above by
$$
\sup _{\mus\pi\in \mus K}\inf_{ \mus V, \mus G, \ve,l,\delta,m,n} \Big\{-\frac1{1+b} {\widehat \J}_{\mus V,\mus G,\varepsilon}^{b,m,n}(\mus \pi)
\Big\}\; .
$$
By \eqv (mparis5)  we have
$$
\limsup_{\ve\to 0}\limsup_{m\to \infty}\limsup_{n\to \infty} {\widehat \J}_{\mus V,\mus G,\varepsilon}^{b,m,n}(\mus \pi)
\,:=\, {\widehat \J}_{\mus V,\mus G}^{b}(\mus \pi)
 =  \left \{ \eqalign {  
  {\widehat J}_{\mus V} (\mus \pi) 
+b\J_G(\mus \pi) &\qquad \hbox{ if }\ \  \mus \pi\in  D([0,T],\MM_1^{ac}\times \MM_1^{ac}), \cr
+\infty &\qquad \hbox{ otherwise } \, .
}\right. 
$$
By \eqv (j-funct) and \eqv (mparis1)  we have that 
$\sup_{ \mus V, \mus G} \Big\{ {\widehat \J}_{\mus V,\mus G}^{b}(\mus \pi)\Big\} = \widehat J(\mus \pi) +b\JJ(\mus \pi)$. 
Therefore,
$$\eqalign{
\limsup_{\g\to 0} \g^d\log
\mus  Q^{\g, \a}_{\s^\g} ( \mus K) & \le
\sup _{\mus\pi\in \mus K}\inf_{ \mus V, \mus G} \Big\{-\frac1{1+b} {\widehat \J}_{\mus V,\mus G}^{b}(\mus \pi)\Big\}
  \cr
\ &\ \ =\, -\frac1{1+b}\inf _{\mus\pi\in \mus K}\sup_{ \mus V, \mus G} \Big\{ {\widehat \J}_{\mus V,\mus G}^{b}(\mus \pi)\Big\} 
 \ =\, -\frac1{1+b}\inf _{\mus\pi\in \mus K}\Big\{\widehat J(\mus \pi) +b\JJ(\mus \pi) \Big\} \, .
}
\Eq(ub56)$$

\qed
\medskip
\noindent {\bf Proof of the upper bound.}  Let  $\mus K$ be a compact set of $D\big( [0,T],\MM_1\times \MM_1\big)$.
By Proposition \eqv (up1),  if $ \JJ \equiv +\infty$ on $\mus K$, then the upper bound of the large deviation principle is satisfied.
Otherwise, there exists $\mus \pi\in\mus K$  such that $\JJ(\mus \pi)<\infty$.  By semicontinuity of the
functional $\mus \pi\mapsto \JJ(\mus \pi)$,  we obtain from \equ(ub56) for any $0<b<1$,
$$
\limsup_{\g\to 0} \g^d\log
\mus  Q^{\g, \a}_{\s^\g} ( \mus K) \, \le\,
-\frac1{1+b} \inf _{\mus\pi\in \mus K\, ,\JJ(\mus \pi)<\infty \; }  \widehat J(\mus \pi)  
 -\frac{b}{1+b}\inf _{\mus\pi\in \mus K} \JJ(\mus \pi)   \; .
$$
Letting $b\to 0$, we get  $$ 
\limsup_{\g\to 0} \g^d\log
\mus  Q^{\g, \a}_{\s^\g} ( \mus K)  \le 
-\inf_{\mus\pi\in \mus K\, ,\JJ(\mus \pi)<\infty \; } \big\{\widehat J(\mus \pi) \big\} 
  \le   -\inf_{\mus\pi\in \mus K} \big\{I_0(\mus \pi) \big\} .
$$
For the last inequality we used   Lemma \equ(up5). 
By  Proposition \eqv (exptight) the proof of the upper bound of the large deviation principle is completed. \qed

\bigskip
\chap{ 7. Lower Bound  }7
\numsec=7
\numfor= 1
\numtheo=1
  
We first get  in Lemma \eqv (A11)  a lower estimate for the probability of a neighborhood of suitable trajectories.   
 We perform the computation  with  the uniform metric on the time interval $[0,T]$ defined as following:
   for $\mu$ and $\nu$ in  $D([0,T],\MM_1)$ and $\rho (\cdot, \cdot)$  defined  in \eqv (par8),
  $$ \r_{[0,T]}(\mu,\nu)= \sup_{t\in [0,T]}\r(\mu_t, \nu_t)  \quad \hbox {and} \quad   \r_{[0,T]}(\mus \mu, \mus  \nu)=   \sum_{i=1,2}\r_{[0,T]}(\mu_i, \nu_i).  
 \Eq  (metric2) $$
 Taking into account that if $d^S_{[0,T]} (\mu,\nu)$  denotes the Skorohod distance, then $$ d^S_{[0,T]} (\mu,\nu) \le  \r_{[0,T]}(\mu,\nu) \Eq (metric3),$$  
 the result holds for the Skorohod topology as well, see [Bill].    
 
 \medskip
To conclude the proof of the lower bound of the large deviation principle in Theorem \eqv (3A), it will remain to show that 
all $\mus \pi$'s  such that $I_0(\mus \pi)<\infty$   can be approximated by a sequence  $(\mus \pi_n)_n$ of smooth 
trajectories, for which Lemma \eqv (A11) holds with  
$\lim_{n \to \infty}  I_0(\mus \pi_n)=  I_0(\mus \pi)$. For this, in Lemma \eqv (3E) we prove that any trajectory  $ \mus m$ smooth enough and far away from the boundaries $(\pm p_1,\pm p_2)$ is associated to a  function  $\mus V (\cdot, \cdot)$. 
\smallskip \noindent 
 Then, given $ \mus m_0 \in B_{p_1,p_2}$,    
denote by  $ \mus R (t,\cdot)$,   $ t \in [0,T]$  the solution of \eqv (G2b) with $ \mus R (0,\cdot)= \mus m_0(\cdot) $: for $i=1,2$,
$$  R_i(t,\cdot) = e^{-t}m_i(0,\cdot)  + p_i \int_0^t    e^{-(t-s)}\tanh [(J* R) (s,\cdot) +a_i \theta]\,ds,  
\Eq(Ri=)$$ 
where $R= R_1+R_2$. 
  It is continuously differentiable in time, actually it is $\CC^\infty $ in time for $ t  \ge t_0 >0$,
  and  there exists  $\delta_i$ which depends on $T$ such that  $ |R_i(t,\cdot)| \le p_i- \delta_i$  for    $ t \in [t_0,  T ]$. 
  Namely,
  since  $|\tanh  z| \le 1-d  $,  for $ |z|   \le K (\beta, \theta) $  with $1>d=d(\beta,\theta)>0$,   we have, for $t \in [0, T]$, 
$$\eqalign {    |R_i(t,\cdot)| &\le e^{-t}m_i(0,\cdot)  + p_i  (1-d)   \int_0^t e^{-(t-s)}\,ds 
\le     p_i [ 1 - d(1-e^{-t}) ].} \Eq (par1) $$
 Recall  that $I_0( \mus R) =0 $, see  (3) of Proposition \eqv (2A). Define the  sets:
  $$ \CC_0 = \CC_0(m_0)=   \{ \mus \phi \in \ac([0,T], B_{p_1,p_2}):  \mus \phi(0)=\mus m_0,   I_0(\mus \phi) < \infty \}, \Eq (25.1)$$  
  $$ \CC_1 = \{ \mus \phi \in \CC_0: \exists\, 0<\eta<T,    \mus \phi (t)= \mus R(t),  t \in [0, \eta ]    \}, \Eq (25.2) $$
 $$ \CC_2 = \{ \mus \phi  \in \CC_1:  \forall   \eta \in (0,T],   \exists \delta_i=\delta_i(\mus\phi)>0,  i=1,2:     \|\phi_i(t)\|_\infty \le p_i - \delta_i, \quad  t \in [\eta, T] \},  \Eq (25.3)$$
    $$ \CC_3 = \{ \mus \phi  \in \CC_2:  \phi_i \in \CC^{2} ( ( 0,T] ,B_{p_1,p_2}), i =1,2 ,\, \phi_i(t)\in\CC(\L),  \forall   t\in (0,T] \}. \Eq (25.4) $$
   By construction
    $ \CC_3 \subset \CC_2 \subset \CC_1 \subset \CC_0 $.  By Lemma \eqv (3E) below we can associate a function 
                    $\mus V$ to  $\mus \phi \in\CC_3$.   
 To extend the lower bound, we show that for $i\in \{ 1,2,3\}$, $ \CC_i $  {\it is $(\rho_{[0,T]},I_0)$-dense } in $\CC_{i-1} $, that is,   for all $\mus  \phi \in \CC_{i-1} $
     there exists a sequence $( \mus  \phi_n)_n \subset   \CC_i$ such that $$\lim_{n \to \infty} \rho_{[0,T]}  (\mus  \phi_n,\mus \phi) =0, \quad  \lim_{n \to \infty} I_0(\mus  \phi_n) = I_0( \mus \phi). \Eq (dense1)$$ 
  This method has been inspired by a similar strategy in [QRV]. 

\medskip
\noindent {\bf \Lemma   (A11)}  {\it  
 Assume   $ (\s^\g)_\g, m_0$ satisfy   \eqv (hyp-G1a). 
  Let $ \d>0$ and  $\mus \mu =\mus  m^{\mus V}  \l$,
where $\mus m^{\mus V}$ is the solution of \eqv (Eq.1) for  $ \mus V = (V_1,V_2) \in \left(\CC^{1,0}  ( [0,T] \times \L)\right)^2$ and $ m^{\mus V}_i (0, \cdot)= p_i m_0 (\cdot)$ for $i= 1,2$. Then we have, for $  \VV_\d (\mus \mu) = \{ \mus \mu' \in D ([0,T], \MM_1 \times \MM_1):   \r_{[0,T]} (\mus \mu,\mus  \mu') < \d \}$, and $ I_{m_0}$ given in \eqv (2.4),
   $$   \liminf_{\g \to 0} \g^d \log  \mus Q^{\g, \a}_{\s^\g} ( \VV_\d(\mus \mu)) \ge   -    I_{m_0} (\mus \mu)  ,   \quad \P-a.s.  
$$ 
 }
\proof
We introduce the perturbed process. By Jensen inequality  we get 
$$
\eqalign{
  \log   \mus Q^{\g, \a}_{\s^\g} (\VV_\d(\mus \mu))
&\ge
E_{\s^\g}^{\mus V,\g,\a}
\left[ \1_{\VV_\d(\mus \mu)}\bigl(\mus \pi^\g_{[0,T]} \bigr) \log \frac{d P_{\s^\g}^{\g, \a} }{d  P_{\s^\g}^{\mus V,\g, \a}} ( \s_{[0,T]}
) \right]
\Bigl(\mus Q_{\s^\g}^{\mus V,\g,\a}(\VV_\d(\mus \mu))\Bigr)^{-1}
 +  \log \mus Q_{\s^\g}^{\mus V,\g,\a}(\VV_\d(\mus \mu)).
}
$$ 
By   Corollary \eqv (ch-g1), $
\lim_{\g\to 0}\mus Q_{\s^\g}^{\mus V,\g,\a}(\VV_\d(\mus \mu))=1$.  
 By  Lebesgue
dominated convergence Theorem,
$$
\liminf_{\g\to 0}\g^d \log   \mus Q^{\g, \a}_{\s^\g} (\VV_\d(\mus \mu))
\ge \liminf_{\g\to 0}E_{\s^\g}^{\mus V,\g,\a}
\left[\g^d \log \frac{d P_{\s^\g}^{\g, \a} }{d  P_{\s^\g}^{\mus V,\g, \a}} ( \s_{[0,T]}
)\right].
$$
By  Radon-Nikodym formula, see Theorem  \equ(RN),  and Proposition \equ(RN1)
  we have 
$$\eqalign { & 
\g^d \log \frac{d P_{\s^\g}^{\g, \a} }{d  P_{\s^\g}^{\mus V,\g, \a}} (
\s_{[0,T]})
\ge -\ell_T\Bigl(\mus \pi^\g(\s),\mus V\Bigr)
+\frac 12\int_0^T\!\G_{\mus V(s,\cdot)}(\mus \pi^\g_s)\,ds \cr & -    \e(\g l) T -   T  C(V_1,V_2)[\d +  \sum_{i=1,2}\EE_i(\d, l,\g,\a)]. }
 $$ 
 {}From Theorem \eqv(pert1), recalling the definition of   $K_{\mus V}(\cdot)$  given  in \eqv(2.2),     we get that for any $l$,
$$\eqalign { & 
\liminf_{\g\to 0}\g^d \log   \mus Q^{\g, \a}_{\s^\g} (\VV_\d(\mus \mu))
\ge -K_{\mus V}(\mus \mu ) -  T  C(V_1,V_2)[\d +  \lim_{\g \to 0} \sum_{i=1,2}\EE_i(\d, l,\g,\a)], }
$$ 
which yields the result letting $l\to\infty$  by Lemma \eqv (erg1) and Proposition \eqv (3A1). \qed

\medskip

\noindent {\bf{\Lemma   (3E)} } {\it  Given $\mus m = (m_1,m_2)\in \big(\CC^{2,0} ( [0,T] \times \L)\big)^2$, with, for $i=1,2$,  
$|m_i (t, r)| <p_i$,
  for all $ t \in [0,T]$, $r \in\L$,  there exists $\mus  V = (V_1,V_2)\in\left( \CC^{1,0}  ( [0,T] \times \L)\right)^2$ such that $\mus m= \mus m^{\mus V} $ 
  is the solution of \eqv (Eq.1).  For $(t,r)\in (0,T]\times \L$,  
  $$\eqalign {& 2 V_i(t,r) =\cr 
&\log 
\left \{   \partial_t m_i (t,r)    \cosh   \left [ (J * m)  (t,r)+a_i \th \right ]   +  \sqrt { \left (\partial_t m_i (t,r)    \cosh   \left [ (J * m)  (t,r)+a_i \th \right ]\right ) ^2+   p_i^2- m^2_i(t,r)}\right \}  \cr & -    [(J *  m) (t,r)+a_i  \th] -\log { \{p_i- m_i(t,r)\} }, }  \Eq (L1)
$$
and for $t=0$ we set  $  \lim_{t \to 0} V_i(t,r)=  V_i(0,r)$.   } 

\smallskip
\proof
   By   \eqv (Eq.1),
     for $t\in(0,T]$,  we determine $\mus  V(t, \cdot) = (V_1(t, \cdot), V_2(t, \cdot))$ with 
$ V_i \in \CC^{1,0}  ( (0,T] \times \L) $ for $i=1,2$, such that $\mus m = \mus m^{\mus V}$.   
 Namely, for $(t,r) \in (0,T] \times \L$, denoting  
 $A _i= ( J * m )(t,r)+ a_i \th $,  $Y_i=     (\cosh  A_i ) \partial_t m_i (t,r) $,  $Z_i= -m_i(t,r) $,   \eqv (Eq.1) is written as
   $$ 
   Y_i   =  
  Z_i \cosh   \left [ A_i+2V_i (t,r)\right ]  + p_i   \sinh  
\left [ A_i +2V_i (t,r)\right ].  $$ 
  We   
  multiply both sides by $X_i= e^{2V_i(t,r)}$  and obtain 
 $$ e^{A_i}(Z_i+p_i)X_i^2-   2Y_i X_i+e^{-A_i}(Z_i-p_i)=0.$$
Its positive solution is 
$$X_i=    \frac { Y_i + \sqrt { Y_i^2-  (Z_i^2-p_i^2)}}{ e^{A_i}(Z_i+p_i) },$$ 
 which gives \eqv  (L1). 
  Note that   $\mus V=(V_1,V_2) $ has the same spatial regularity as  $\mus m$, namely the argument of the square root is always strictly positive.  
   \qed  

\medskip
\noindent {\bf {\Corollary  (26J)} }   {\it If $\mus m$ is solution of \eqv (G2b) then $\mus V=0$ in \eqv (L1). 
   }
   \medskip\smallskip
    \noindent 
   \Remark  (26JJ) Lemma \eqv (3E) could have been stated requiring $ \mus m \in   \ac ([0,T],B_{1,1})$.
 In this case one would get $ \mus V\in \big(L^1([0,T], \CC(\Lambda))\big)^2$.  
 We prefer to obtain more regularity in time for $\mus V$, so that uniformity and other technical needs become straightforward. 
  
 \medskip\noindent {\bf \Lemma (72)}        {\it  
  $   \CC_1 $   is  $(\rho_{[0,T]},I_0)$-dense   in $\CC_{0} $.}

\smallskip
\proof  Fix $ \mus  m \in \CC_{0} $. Let $\mus  R(t,\cdot)$, $t\in[0,T]$,
be the solution of \eqv (G2b)  with initial datum $ \mus R(0, \cdot)= \mus m_0(\cdot)$.
For any $\eta \in (0,T)$, define
$$ \mus m^{\eta} (t, \cdot )= \left \{ \eqalign { & \mus R
 (t, \cdot)  \qquad \hbox {for } \qquad  t \in  [0,\eta], \cr &
\mus R (2\eta - t, \cdot)  \qquad \hbox {for } \qquad  t \in (\eta, 2\eta], \cr &
\mus m (t- 2\eta, \cdot)  \qquad \hbox {for } \qquad  t \in ( 2\eta, T].  }  \right .  $$
We  have  $\mus m^{\eta} \in \CC_{1}$ for any $0<\eta<T$ and
$\lim_{\eta \to 0} \rho_{[0,T]} (\mus m^{\eta}, \mus m)=0$.  Since $I_0$ is lower semicontinuous it remains to show
$$ \lim_{\eta\to 0}  I_0(\mus m^\eta)  \le I_0 (\mus m ).  \Eq (7.2b)  $$ 
  We  split  $[0,T]$ into  $[0,2\eta]$ and  $[2\eta,T]$ in the integration. We have that
$$  
    \int_{2\eta}^T \int_\L  \HH ( \mus m^\eta , \dot {\mus m^\eta}) (t,r)\,drdt   =
    \int_0^{T-2 \eta}  \int_\L   \HH ( \mus m, \dot {\mus m})  (t,r)\,drdt   \le
I_0(\mus m).  
$$
 Next we show that 
$$ \lim_{\eta \to 0} \int_{0}^{2\eta}\int_\L  \HH ( \mus m^\eta , \dot {\mus m^\eta}) (t,r)\,drdt    = 0. $$
Since  $\mus m^\eta=  \mus R $ for $ t \in [0,\eta]$  solves  \eqv (G2b), by (3) of Proposition \eqv (2A), 
$$ \int_{0}^{ \eta}  \int_\L \HH ( \mus R, \dot {\mus R}) (t,r)\,drdt    = 0. \Eq (7.29)$$ 
Since the profile $\mus m^\eta$ in $(\eta, 2\eta]$ is the profile in $(0,\eta]$ backwards in time, we have 
$$ \int_{\eta}^{ 2\eta}\int_\L  \HH ( \mus m^\eta , \dot {\mus m^\eta}) (t,r)\,drdt= 
   \int_{0}^{ \eta}\int_\L \HH ( \mus R, -\dot {\mus R}) (t,r)\,drdt. 
$$ 
  Because  $ \mus R$ solves  \eqv (G2b) and  for $ t>0$, $ |R_i(t, \cdot)| \le p_i- \delta_i$, for $i=1,2$,  $ \HH ( \mus R, \dot {\mus R})$  belongs to  $ L^1 ([0,T] \times  \L)$,    as well as $ \HH ( \mus R, -\dot {\mus R})$, see  explicit formula \eqv(DD1c).
By dominated convergence,  
    $$\lim_{\eta \to 0}  \int_{\eta}^{ 2\eta} \int_\L  \HH ( \mus m^\eta , \dot {\mus m^\eta}) (t,r) \,drdt    =0.$$
    In this way we prove \eqv (7.2b).   
 \qed
 
\medskip
\noindent {\bf \Lemma  (28J)}   {\it  $\CC_2$ is $(\rho_{[0,T]},I_0)$-dense in $ \CC_1$.}

\smallskip
\proof   Let   $ \mus m  \in   \CC_1$ and $\eta \in (0,T)$ so that $ \mus m (t,\cdot)= \mus R(t,\cdot)$ for $t \in [0,\eta]$.  By  \eqv (par1), 
$ \|m_i ( \eta, \cdot )\|_\infty \le p_i- \d_i$ for some $ \d_i>0$  and  $i=1,2$. 
    Define    
$$
m_i^n (t,r)= \left \{  \eqalign { &m_i ( t, r )  \quad\hbox{ for } t \in [0, \eta], \cr 
     &m_i(  \eta,r)+  \left(1 -\frac 1 n\right) (m_i (t,r)-  m_i(  \eta,r))
            \quad\hbox{ for } t \in  ( \eta, T]. } \right.  
\Eq (V.1)
$$
By construction and from \equ(7.29),
$\displaystyle{I_0 (\mus m^n)=\int_\eta^T\int_\Lambda \HH (\mus m^n (t,r),  \frac {\partial   \mus m^n } {\partial t} (t,r)  )\,drdt}$.
Moreover, since  $I_0 (\mus m) < \infty$, by  Proposition \eqv(3A22) 
we have  $ \|m_i (t)\|_\infty \le p_i$ for $t \in [\eta, T]$, then  $$ \|m_i^n (t)\|_\infty \le  p_i - \frac {\d_i} n,
\quad \forall t\in[\eta,T].  
\Eq (V.2a) $$
Hence $\mus m^n \in \CC_2$ for all $n$. Furthermore  $\displaystyle{\lim_{n \to \infty} m_i^n (t,r) =  m_i  (t,r)}$ and 
$ \displaystyle{ \frac {\partial  m_i^n } {\partial t} (t,r) =  (1 -\frac 1 n)    \frac {\partial  m_i } {\partial t} \to   
\frac {\partial  m_i  } {\partial t} (t,r)}$  for almost all  $(t,r)\in [\eta,T]\times\Lambda$.
Then, by Proposition \eqv(3A22), $ \displaystyle{\HH (\mus m^n (t,r),  \frac {\partial\mus m^n } {\partial t} (t,r)  ) }$  is given by \equ(DD1c), while $\displaystyle{\HH (\mus m  (t,r),   
\frac {\partial   \mus m } {\partial t}    (t,r) )}$ is given either by  \equ(DD1c) when $|m_i(t,r)|<p_i$,
or, when $|m_i(t,r)|=p_i$, by \equ(DD1b), or is infinite. We hence check that pointwise 
$$ \lim_{n \to \infty} \HH (\mus m^n (t,r),  \frac {\partial   \mus m^n } {\partial t} (t,r))       
=  \HH (\mus m  (t,r),   \frac {\partial   \mus m } {\partial t}    (t,r) ). 
$$
To apply the Lebesgue dominated convergence Theorem we give an upper bound, uniformly with respect to $n$, of $\displaystyle{| \HH (\mus m^n,  \frac {\partial   \mus m^n } {\partial t})(t,r) |}$ 
(see also [C] p. 174). For that we combine \eqv (2.1d) with the facts that, 
$$ \{ (t,r):  \frac {\partial  m_i^n } {\partial t} >0 \}=  \{ (t,r):  \frac {\partial  m_i } {\partial t} >0\}, $$
    and on the set  $ \{ (t,r):     m_i^n (t,r)   \ge  p_i- \delta_i \}$ we have 
    $ m_i  (t,r)-  m_i (0,r)\ge   0$ and $p_i- m^n_i (t,r) \ge  p_i-
    m_i (t,r) $.  To get shorter notation, we 
    denote for $ \mus \phi, \mus \psi \in  \ac ([0,T],  B_{p_1,p_2})$ 
    $$ \Upsilon (\phi_i, \psi_i)=  \1_{ \{ \dot \phi_i >0;\,  p_i-\delta_i\le \psi_i \} }\left( \log \frac 1 {p_i-\psi_i} \right)^++  \1_{ \{ \dot \phi_i<0;\,   -p_i+\delta_i\ge   \psi_i \} }\left( \log \frac 1 {p_i+\psi_i}\right)^+.  $$
   We have
        $$\eqalign { &  2\HH (\mus m^n,  \frac {\partial   \mus m^n } {\partial t})(t,r)   \cr & \le  \sum_{i=1,2}
      | \dot m_i| \left [ \left(\log | \dot m_i|\right)^+ +  \Upsilon ( \dot m_i, m_i^n)   +K_i \right ](t,r)   
       \cr &+  \sum_{i=1,2} | \dot m_i| \left [   \1_{ \{\dot m_i>0;\,  m_i^n<p_i-\delta_i\} }\left( \log \frac 1 {p_i-m_i^n}\right)^+ +  \1_{\{ \dot m_i<0; \,   m_i^n>-p_i+\delta_i\}} \left(\log \frac 1 {p_i+m_i^n}\right)^+ \right ](t,r) + C 
      \cr & \le \sum_{i=1,2}  | \dot m_i| \left [ \left(\log | \dot m_i|\right)^+ +\Upsilon ( \dot m_i, m_i) +K_i \right ](t,r)    \cr &
     +   \sum_{i=1,2} | \dot m_i| \left [   \1_{ \{\dot m_i>0;\, m_i^n<p_i-\delta_i\}} \log \frac 1 {\delta_i} +  \1_{\{ \dot m_i<0;\,    m_i^n>-p_i+\delta_i\}} \log \frac 1 {\d_i} \right ](t,r)  + C 
 \cr &    
 \le   \sum_{i=1,2}
 | \dot m_i| \left [ \left(\log | \dot m_i|\right)^+ + \1_{ \{\dot m_i>0\}}  \left(\log \frac 1 {p_i-m_i}\right)^+ +  \1_{ \{\dot m_i<0\}}   \left(\log \frac 1 {p_i+m_i}\right)^+ +K_i +\log \frac 1 {\d_i} \right ](t,r)      +C.  
 }  
$$  
  Since by assumption  $ I_0 (\mus m) < \infty$, by Proposition \eqv (3A3),  part (b), the above  upper bound is integrable. By Lebesgue dominated convergence Theorem   we then have 
    $$ \lim_{n \to\infty}  I_0 (\mus m^n)=  I_0 (\mus m). \Eq (S.1)$$
    Obviously $   \mus m^n \to    \mus m $ in the metric \eqv (metric3a). 
   \qed
\medskip

\noindent {\bf \Lemma  (J4)} { \it  $\CC_3 $  is $(\rho_{[0,T]},I_0)$-dense  in $\CC_2$. }  

\smallskip
 \noindent \proof
Take $ \mus \psi \in \CC_2$.  To get more regularity we    convolve   with a smooth kernel the function both  in time and space.
 To  perform the convolution in time 
 we extend the definition of $\mus \psi $ to $[T,T+1]$   by setting, for each $s\in [0,1]$, 
 if $\mus u=({\mus u}_1,{\mus u}_2)$ is the solution of equation 
\equ(G2b) with initial condition $\mus \psi (T, \cdot)$,
$$\mus \psi (T+s,r)=\mus u (s,r). \Eq (S.2z)$$
Since $ \mus \psi \in \CC_2$ there exist $  \delta_i$, $i=1,2$, such that    $ |\psi_i  (T,r)| \le p_i-   \delta_i $. It  follows from
 \equ(Ri=) that $\mus \psi_i (T+s,r)\le p_i-{\tilde \delta}_i$   for all $s \in [0,1]$,  for some ${\tilde \delta}_i$  smaller than $ \delta_i$.  
In the following we will denote it  always by $ \delta_i$.  
 Denote by $\theta_s \mus  \psi $ the time translation
of $\mus \psi$,  $  (\theta_s\mus \psi) (t,r)= \mus \psi(t+s,r)$ for $(t,r) \in [0,T] \times  \L$. 
Let $\Phi_{ \e_1} $ be a smooth non-negative kernel,   $ \Phi_{ \e_1 }  \in  \CC^\infty ( \L) $ with support in a ball of radius $\e_1$ and integral one
which we use as spatial mollifier. For $\e_0>0$, 
let $\Psi_{\e_0} $ be the $\CC^\infty (\R) $
non-negative temporal mollifier with support $[0,\e_0]$ and integral one. 
Set $\e \equiv (\e_0,   \e_1 )$,   $ \e \downarrow 0$  stands for  $  \e_0 \downarrow 0$ and $ \e_1 \downarrow 0$.    
Let $\eta>0$ be such that $\mus \psi(t,\cdot) = \mus R(t,\cdot)$ for $t\in[0, 3\eta]$.  Let  $ \chi_1 (t) $,  $ \chi_2 (t) $  be a  $\CC^2$ partition of the unity enjoying the properties:
$$  \left \{ \eqalign { & \chi_1 (t) = 1\qquad \hbox {for} \qquad t\in [0,   \eta ], \qquad    \chi_1 (t) = 0 \qquad  \hbox {for} \qquad t\in [2 \eta, T], \cr & 
 \chi_2 (t) = 0 \qquad\hbox {for} \qquad t\in [0,  \eta ], \qquad    \chi_2 (t) = 1 \qquad \hbox {for} \qquad t\in [ 2\eta, T],  \cr & 
  \chi_1 (t)+\chi_2 (t) = 1, \qquad  \forall t \in [0,T]. } \right. $$
Let 
$$ \psi^{\e}_i (t,\cdot) = \chi_1 (t) \psi_i  (t, \cdot) +  \chi_2 (t) \int_{ \R } \Psi_{\e_0}(s)   ( \Phi_{ \e_1} * \theta_s \psi_i) (t, \cdot)\,ds. \Eq (7.5) 
 $$ 
By construction $\psi_i (\cdot,\cdot)$, $i=1,2$, satisfies all the regularity requirements to be in  $\CC_3$. Furthermore, since  $ |  \psi_{i }  (t,r) | \le    p_i-  \d_i, \delta_i>0$, for all $\e>0$ and $t \in  [0,T]$, we still have that   
  $$  |   \psi^\e_{i }  (t,r) | \le    p_i-  \d_i, \qquad i=1,2,  
 \Eq (Sa.6) $$
therefore $\mus\psi^{\e}\in\CC_3$. 
 Moreover
   $$ \lim_{\e \to 0} \,  \r_{[0,T]} ( \mus \psi^{\e}, \mus \psi)=0.$$     
Since $I_0$ is  lower semicontinuous, see Proposition \eqv (2A), (1), it is enough to prove 
$$ \lim_{\e \to 0}  I_0( \mus \psi^{\e}) \le I_0( \mus \psi ).  \Eq   (EJu4a) $$ 
    By using the expression \eqv(funct2) of $I_0$,
see  Proposition \eqv(3A1),   
  we have
  $$   I_0 ( \mus \psi^{\e  }) - I_0 ( \mus \psi) =
  \int_0^T \int_{\L} \left[\HH (\mus \psi^\e,  \frac {\partial   \mus \psi^\e} {\partial t}) (t,x)-
   \HH (\mus \psi,  \frac {\partial   \mus \psi} {\partial t}) (t,x)\right]\,dxdt. 
\Eq (Sa.3) $$
We split the time integral into 3 pieces: (i) a first integral on $[0,\eta]$, which is equal to 0 
  by  definition \eqv (7.5) of $\mus \psi^\e$; (ii) a second one on $[\eta,2\eta]$,  treated 
   in Lemma \eqv (EA10) below; (iii) a third one on $[2\eta,T]$, that we now analyze.
   Notice that for $t \ge 2 \eta$, see \eqv (7.5), $\chi_1(t)=0$ and  $\chi_2(t)=1$, therefore  $\psi^{\e}_i(t,\cdot)$  reduces to  a convex combination,   and we exploit that 
$ \HH (\mus m,  \mus a )$ is convex with respect  to $ \mus a $.   Then, for $ t \ge 2 \eta $,  by   Jensen inequality we obtain 
$$ \HH (\mus \psi^\e,  \frac {\partial   \mus \psi ^\e} {\partial t}) (t,x)   \le
   \int_{\R} \Psi_{\e_0} (s)  \int_{\L} \Phi_ {\e_1} (y)   \HH (\mus \psi^\e  (t,x),  \frac {\partial   \mus \psi } {\partial t}  (t+s ,x-y) )\,dyds. \Eq (D.5) $$   
     For all $ s \in [0,1]$, $s <T$, 
we have
 $$\eqalign { &  \int_{2\eta+s}^T  \int_\L \HH (\mus \psi, \frac
 {\partial   \mus \psi } {\partial t}) (t,x)\,dxdt   =
  \int_{2\eta}^T  \int_\L  \HH (\mus \psi, \frac {\partial   \mus \psi } {\partial t}) (t+s,x)\,dxdt  
  -  \int_{T-s}^T  \int_\L   \HH (\mus \psi, \frac {\partial   \mus \psi } {\partial t}) (t+s,x)\,dxdt \cr &
  =
  \int_{2\eta}^T  \int_\L  \HH (\mus \psi, \frac {\partial   \mus \psi } {\partial t}) (t+s,x)\,dxdt   
\cr &  = 
    \int_{\R} \Psi_{\e_0} (s)   \int_{2\eta}^T {d}t  \int_\L \HH (\mus \psi, \frac {\partial   \mus \psi } {\partial t}) (t+s,x)\,  dxds    \cr & =
       \int_{\R} \Psi_{\e_0} (s)   \int_{2\eta}^T\int_\L dy  \Phi_ {\e_1} (y)  \int_\L \HH (\mus \psi, \frac {\partial   \mus \psi } {\partial t}) (t+s,x-y)\,dxdtds, }  \Eq (EJ4)$$
where the first equality comes from a change of variables, the second one from 
the definition of $\mus \psi$ in $[T,T+1]$ (see  \eqv  (S.2z)), the third one from
$\int_\L dy  \Phi_ {\e_1} (y)=1$ and $ \int_\L  {d}x \HH (\mus \psi, \frac {\partial   \mus \psi } {\partial t}) (t,x-y) =   \int_\L  {d}x \HH (\mus\psi, \frac {\partial   \mus \psi } {\partial t}) (t,x) $, and the last one from
$\int_{\R}   {d}s \Psi_{\e_0} (s)   =1$.   
  Therefore 
$$ \eqalign {&  \int_{2\eta}^T \int_{\L} \left[\HH (\mus \psi^\e,  \frac {\partial   \mus \psi^\e} {\partial t}) (t,x)-
   \HH (\mus \psi,  \frac {\partial   \mus \psi} {\partial t}) (t,x)\right]\,dxdt\cr &
   = \int_{2\eta}^T \int_{\L} \HH (\mus \psi^\e,  \frac {\partial   \mus \psi^\e} {\partial t}) (t,x)\,dxdt -
 \int_{2\eta}^{2\eta +s}  \int_\L  \HH (\mus \psi,  \frac {\partial   \mus \psi} {\partial t}) (t,x)\,dxdt  
 - \int_{2\eta +s} ^T  \int_\L\HH (\mus \psi,  \frac {\partial   \mus \psi} {\partial t}) (t,x)\,dxdt \cr &
  \le
  \int_{\R} \Psi_{\e_0} (s)      \int_{2\eta}^T \int_{\L}   \Phi_
   {\e_1} (y)
  \int_{\L} \HH (\mus \psi^\e  (t,x),  \frac {\partial   \mus \psi } {\partial t}  (t+s ,x-y) )\,dxdydtds\cr &-
     \int_{\R}   \Psi_{\e_0} (s)   \int_{2\eta}^T \int_\L \Phi_ {\e_1} (y)  \int_\L  \HH (\mus \psi, \frac {\partial   \mus \psi } {\partial t})(t+s,x-y)\,dxdydtds .     
    } \Eq (Sa.30) $$
       The 
     inequality holds by  \eqv (D.5), and because
     $ \HH (\mus \psi, \frac {\partial   \mus \psi } {\partial t}) (t,x)  \ge 0$.  Finally we use \eqv (EJ4).    To estimate the last difference  in \eqv  (Sa.30) we   
  add   and subtract    to  it  the term 
 $$ 
   \int_{\R} {d}s \Psi_{\e_0} (s)  \int_{2\eta}^T {d}t  \int_{\L} \Phi_ {\e_1} (y)  dy  \int_{\L} {d}x   \HH (\mus \psi(t,x),  \frac {\partial \mus \psi }{\partial t}  (t+s ,x-y) ),  $$ 
     which gives
     $$  \int_{2\eta}^T  \int_{\L} \left[\HH (\mus \psi^\e,  \frac {\partial   \mus \psi^\e} {\partial t}) (t,x)-
   \HH (\mus \psi,  \frac {\partial   \mus \psi} {\partial t}) (t,x)\right]\,dxdt \le W_1 +W_2,$$
 where 
   $$ \eqalign {W_1 &=  
    \int_{\R}   \Psi_{\e_0} (s)      \int_{2\eta}^T   \int_{\L}   \Phi_ {\e_1} (y)   \int_{\L}   \left[\HH(\mus \psi^\e  (t,x),  \frac {\partial   \mus \psi } {\partial t}  (t+s ,x-y) )-  \HH (\mus \psi     (t,x),  \frac {\partial   \mus \psi } {\partial t}  (t+s ,x-y) )\right]\,dxdydtds, \cr  
  W_2 &=  \int_{\R}   \Psi_{\e_0} (s)      \int_{2\eta}^T    \int_{\L}   \Phi_ {\e_1} (y)   \int_{\L} \left[ \HH(\mus \psi     (t,x),  \frac {\partial   \mus \psi } {\partial t}  (t+s ,x-y) ) -
  \HH (\mus \psi , \frac {\partial   \mus \psi } {\partial t})(t+s,x-y) \right]\,dxdydtds.   } \Eq (Sa.4)$$
   Taking into account  Lemma   \eqv  (J40) below  we get the result.\qed
 \medskip 
 The proofs of the next two Lemmas are postponed to  Appendix B.
 \medskip

 \noindent {\bf \Lemma  (J40)} { \it  
 $$ \lim_{\e \to 0}  |W_i| = 0, \qquad i=1,2.  $$} 

\medskip 
 
 \noindent {\bf \Lemma   (EA10)}  {\it  
 $$  \lim _{\e\to 0} \int_{ \eta }^{2\eta}\int_{\L} \left[\HH (\mus \psi^\e,  \frac {\partial   \mus \psi^\e} {\partial t}) (t,x)-
   \HH (\mus \psi,  \frac {\partial   \mus \psi} {\partial t}) (t,x)\right]\,dxdt= 0. $$ 
  }
   
  \bigskip
\chap{8. Appendix A}8
\numsec= 8
\numfor= 8
\numtheo=8
  
In this Appendix we give the proofs of the properties of the rate functional stated in Section 3.

 \noindent {\bf  Proof of Lemma \eqv (2C1)}  
  The differentiability  of  $\G_{\mus V}(\mus u)$
    in $(L^\infty (\L))^2$ is easily verified. For the convexity  we
    compute first  the Hessian of  $ \G_{\mus V}(\mus u) $ with respect to $V_1$ and $V_2$. Since the Hessian is a diagonal matrix, it is enough to study separately the convexity with respect to $V_1$ and $V_2$, we do it for $V_1$.    For $r\in\L$, we set $V_1(r)=x,u_1(r)=m$ and denote by   
 $$f_1(x) = (p_1  \tanh \vartheta  -  m)  \sinh (2x ) + (p_1 
  - m  \tanh  \vartheta )(\cosh (2x) -1)    $$
  the integrand term in $ \G_{\mus V}(\mus u) $ which depends only on $\mus V$  
   with   $\vartheta$ varying in some bounded interval of $\R$,   $x \in \R$, $ |m| \le p_1 $.  We then   study the sign of the second derivative of  $f_1$.  
    $$\eqalign {\frac 14  f_1'' (x)
  = \cosh (2x)    [ p_1-m  \tanh \vartheta] +  \sinh (2x ) [ p_1 \tanh \vartheta  -m]}.$$
  Notice that $ p_1-m  \tanh \vartheta \ge 0$, 
  and $ p_1-m \tanh \vartheta \ge p_1 \tanh \vartheta  -m  \ge -(p_1-m \tanh \vartheta ). $
   Since $\cosh(2x) > |\sinh(2x)| $  and $| \tanh \vartheta| <1$ when $\vartheta$ varies in a bounded interval  we obtain that $f_1$ is convex. \qed 
   \bigskip
   \noindent {\bf  Proof of Proposition  \eqv (3A22)}
Recalling \equ (2.1n), for $i=1,2$,
denote 
$$ \eqalign { & F_i(v_i)= g_i v_i- \frac 12  B_i(\mus u, v_i) \cr & =  g_iv_i-(p_i-u_i)\frac { e^{  A_i}  } { 4\cosh  A_i }[e^{2v_i}-1]-(p_i+u_i)\frac { e^{- A_i}  } { 4\cosh    A_i }[e^{-2v_i}-1].}\Eq (F1bis)
$$
Hence 
 $$ \frac{ \partial {F_i}}{\partial v_i}= g_i- (p_i-u_i)
   \frac { e^{ A_i}  } {2\cosh  A_i }    e^{ 2v_i } 
 +( p_i+u_i)  \frac { e^{- A_i}  } { 2\cosh A_i }  e^{ -2v_i }.  \Eq (JP1) $$
 First assume that $u_i>p_i$.  By \eqv (F1bis), since $-(p_i-u_i)>0$, we have $\lim_{v_i\to 
+\infty}F_i(v_i)=+\infty$. In the same way, we get $\lim_{v_i\to
-\infty}F_i(v_i)=+\infty$ if $u_i<-p_i$ because then $-(p_i+u_i)>0$. 
Therefore, (a) holds. 

For the remaining cases, we exploit that for $\mus u\in B_{p_1,p_2}$, the function $v_i \mapsto B_i(\mus u, v_i)$ is convex  differentiable  on $\R$.   

\noindent (b) To compute the Legendre transform  of $ B_i(\mus u, v_i)$, 
when 
  $|u_i|< p_i $,   by  \eqv (JP1), 
the maximum  in
\eqv (Le2c) is obtained for  (remember \equ (DD1a))  
  $$e^{ 2v_i }= e^{- A_i}\frac { D_i}  {p_i- u_i }, \qquad 
 {\rm hence\ } v_i =  \frac 12 \left(  \log\frac{ D_i}{p_i- u_i}- A_i\right). 
  \Eq (L1aA) $$ 
 Inserting \eqv   (L1aA) in \eqv (F1bis)  
   we have 
 $$\eqalign {   H_i (\mus u, g_i) & =   
   \frac {g_i}2 \left [\log \frac{D_i}{p_i- u_i}  -A_i \right ] 
+ p_i \frac { e^{ A_i}+e^{ -A_i}  } { 4\cosh  A_i }
-u_i\frac { e^{ A_i}-e^{ -A_i}  } { 4\cosh  A_i }
   - \frac 1{ 4\cosh  A_i} \left [
  D_i + 
  \frac {p_i^2- u_i^2 } { D_i}\right ],  }  
  $$
  which yields \eqv (DD1c) since, using \equ (DD1a), we write
  $$ \eqalign {     D_i + 
  \frac {p_i^2- u_i^2 } { D_i}  & =
        g_i \cosh A_i +  R_i+ 
  \frac { ( p_i^2- u_i^2)  \left ( g_i   \cosh A_i - R_i \right ) }
  {    ( g_i  \cosh A_i)^2  - \left ( g_i \cosh A_i\right )^2 - p^2_i+ u^2_i  } 
 =
    2R_i.
   }   
   $$  
 
 \noindent 
 (c) When $ u_i =p_i$   (resp. $ u_i =-p_i$), 
$$ \frac{ \partial {F_i}}{\partial v_i}= g_i
 +p_i \frac { e^{-{\rm sgn} (u_i) A_i}  } { \cosh   A_i }    e^{ -2 {\rm sgn} (u_i) v_i} $$
 and  to solve  $\displaystyle{\frac{ \partial {F_i}}{\partial v_i}=0}$ (that is to find a finite extremum)  we need $g_i<0$ (resp. $g_i>0$), 
 namely 
$$g_i= -  p_i  \frac {e^{-{\rm sgn} (u_i)A_i}} { \cosh  A_i} e^{ -2  {\rm sgn} (u_i) v_i } .  
$$
Inserting this value  in   \eqv (F1bis)     we get \eqv (DD1b) when $g_i\neq 0$.
 
 When $ u_i =p_i$  (resp. $ u_i =-p_i$) and $g_i=0$, \eqv (F1bis) becomes
$$F_i(v_i)=p_i\frac { e^{-{\rm sgn} (u_i) A_i}  } { 2 \cosh   A_i }[1-e^{-2  {\rm sgn} (u_i) v_i}].
$$
It is an increasing  (resp. decreasing) function with a finite maximal limit:
$$\lim_{v_i\to 
+\infty}F_i(v_i)= p_i\frac { e^{- {\rm sgn} (u_i) A_i}} {2 \cosh A_i }=H_i (\mus u,g_i),  \quad  \hbox {resp.} \quad    \lim_{v_i\to 
-\infty}F_i(v_i)=p_i\frac { e^{- {\rm sgn} (u_i) A_i} } { 2\cosh   A_i }=H_i (\mus u,g_i).  $$

\noindent (d)  When $ u_i=p_i$ and $g_i>0$  (resp. $ u_i=-p_i$ and $g_i<0$),   \eqv (F1bis) becomes
$$F_i(v_i)=g_iv_i+p_i\frac { e^{-{\rm sgn} (u_i) A_i}  } { 2\cosh   A_i }[1-e^{-2 {\rm sgn} (u_i) v_i}].
$$
Hence
$$\lim_{v_i\to +\infty}F_i(v_i)=+\infty=H_i (\mus u,g_i), \qquad \hbox {resp.} \quad \lim_{v_i\to -\infty}F_i(v_i)=+\infty=H_i (\mus u,g_i). $$
 \qed 
 
  \smallskip\medskip 
     \smallskip\medskip  
  \noindent {\bf  Proof of Proposition  \eqv (3A3)} 
  We use the explicit representation of $ \HH(\cdot, \cdot)$ given in  Proposition \eqv(3A22).  
 \medskip (a)  We give an upper bound  of expression \eqv (DD1c). 
 The  difficulty comes from the term 
$$ F(\mus u, g_i,\th)= 
  g_i \log \frac {D_i (\mus u,g_i,\th)} { p_i- u_i},   \Eq (EJ1) $$ 
   where $D_i (\mus u,g_i,\th)$  is  defined in \eqv (DD1a).  
 Let $- \mus u= (-u_1, -u_2)$.  We have
 $$\eqalign {   F(- \mus u, -g_i,\th)&=  
g_i \log \frac  {(p_i + u_i)\left \{    g_i    \cosh   [ (J* u) -a_i\th]   +  \sqrt { \left (g_i    \cosh  [  (J* u) -a_i\th ]\right ) ^2+   p_i^2 - u^2_i }\right\} } {
  - \left ( g_i \cosh [ (J* u) -a_i\th ] \right )^2   +   \left (g_i  \cosh [(J* u) -a_i\th ]\right ) ^2+   p_i^2 - u^2_i
 }  \cr & =
g_i \log \frac  {     g_i \cosh [ (J* u) -a_i\th ]   +  \sqrt { \left (g_i \cosh [(J* u) -a_i\th ]\right ) ^2+   p_i^2 - u^2_i } } {
    p_i  - u_i
 } \cr & =  F(\mus u, g_i,-\th). 
 } 
 $$  
We write
 $$   F(\mus u, g_i,\th) 
   = F(\mus u, g_i,\th)\1_{\{g_i\ge 0\}}+  F(-\mus u, -g_i,-\th)\1_{\{g_i< 0\}}. 
  \Eq (DD2ZZ) $$ 
Hence it suffices to estimate  $ F(\mus u, g_i,\th)$ for $g_i>0$ and $\th\in \R$. We get 
 $$  F(\mus u, g_i,\th) \le |g_i| \left\{ \log D_i (\mus u,g_i,\th)  + \left( \log \frac 1 {p_i  - u_i}  \right )^+\1_{\{g_i>0\}} + \left( \log  \frac1{p_i + u_i} \right )^+\1_{\{g_i< 0\}} \right \}.  
 $$
 We  obtain \eqv (2.1d)  by the upper bound  $D_i (\mus u,g_i,\th) \le 2 |g_i| + 1$.
  The lower  bound  \eqv (2.1e) is obtained as in [C], p.  171. We rely on formulas  \eqv (Le2c),  \eqv (2.1n). Since
  $  e^{\b a} \le  e^{\b |a|}$, there exists a  constant $C$ such that
$$  
B_i(\mus u, v_i)
 \le 2C    [ e^{ 2|v_i| }-1 ]  
 $$ 
$$2 p_i \ge 2C \ge   \max \left \{    ( p_i-u_i)  \frac { e^{ A_i}  } { 2\cosh  A_i },   ( p_i+u_i)  \frac { e^{- A_i}  } {2 \cosh A_i } \right\}.   
$$
  Then
 $$  H_i(\mus u, g_i) \ge   \sup_{v_i \in \R} \left \{ g_i v_i-  C  \left [ e^{ 2|v_i| }-1 \right ]   \right \}= \max  \left \{ \frac { |g_i|}2 \left [ \log {\frac {|g_i|}{2C} }- 1\right ] +C, 0 \right \}.  $$
  \medskip (b)
 If \eqv (2.1d) holds then  $I_0(\mus \phi)<\infty$. 
  For the converse,  by \eqv (2.1e), it is necessary to have $\dot \phi_i \log |\dot \phi_i|\in L^1  ( [0,T] \times \L)$.  To conclude,   notice that 
  when  $g_i>0$, uniformly in  $\theta \in \R$,
 $$F(\mus u,g_i,\th)-g_i\log\{g_i\cosh A_i\}  \ge  2 g_i \1_{\{g_i>0\}} \log\frac 1{p_i-u_i}.  $$  
\qed 

   \bigskip
   \noindent {\bf  Proof of Proposition  \eqv (2A)} 
For (1), (2) we refer to the similar proof of  [C], Theorem III.4, p. 148 (indeed, the rate functional in infinite outside
$\CC([0,T], \BB_{1,1})$).
To show the first part of  (3), notice that  for $ \mus V=0$ the r.h.s. of the argument of the $\sup$ in  \eqv (2.2a) is equal to zero.  This implies that for $s \in [0,T]$, 
 $\HH^* (\mus \phi(s,\cdot),   \dot {  \mus \phi }(s,\cdot) )\ge 0 $ in \eqv (funct1),       
  therefore  $I_0( \mus  \pi) \ge 0$.   
  For the second half of (3),  
  we start by proving that if $I_0(\mus \pi)=0$, then  $ \mus  \pi = (\phi_1 \l, \phi_2\l)$  with $\mus \phi= (\phi_1, \phi_2)  \in
\ac ([0,T], B_{1,1})$  is the solution of equation  \eqv (G2b).
{}From Proposition \eqv(3A1),  we know  that
$J_0(\mus\pi)=0$ (see  \eqv(DC1), \eqv(2.2)), that is, for any $\mus V= (V_1, V_2) \in
\left(L^\infty   ( [0,T] \times \L)\right)^2$, we have
$$
\int_0^T < \mus V(s,\cdot),   \dot {\mus \phi} (s,\cdot)>ds
\le  \frac 12  \int_0^T \G_{\mus V(s,\cdot)}( \mus \phi (s,\cdot))\,ds.
$$
Now take $V_2=0$ and $\eta V_1$ instead of $V_1$, where
$\eta>0$.  Denote $\phi=\phi_1+\phi_2$, then recalling
definitions \eqv(rate1)  and  \eqv(2.1), we get
$$
\eqalign{
&2\eta\int_0^T < V_1(s,\cdot),   \dot {\phi_1} (s,\cdot)>ds \cr
&\le p_1  \int_0^T
<\tanh(\phi(s,\cdot)*J+a_1\theta)\sinh(2\eta V_1(s,\cdot))+\cosh(2\eta V_1(s,\cdot))-1 >ds\cr
&\quad-\int_0^T<\phi_1(s,\cdot)\Bigl(\tanh(\phi(s,\cdot)*J+a_1\theta)[\cosh(2\eta
V_1(s,\cdot))-1]+\sinh(2\eta V_1(s,\cdot))\Bigr)>ds.\cr
}
$$
Using Taylor expansion in $\eta$ when $\eta\to 0$,
dividing by $\eta$ and letting $\eta\to 0$, we obtain
$$
\int_0^T \!\!< V_1(s,\cdot),   \dot {\phi_1} (s,\cdot)>ds
\le {p_1}  \int_0^T\!\!\!\!
<\tanh(\phi(s,\cdot)*J+a_1\theta)V_1(s,\cdot)>ds
-\int_0^T\!\!\!\!<\phi_1(s,\cdot)V_1(s,\cdot)>ds.
$$
Since all  terms in the previous expression are linear in
$V_1$, we may change $V_1$ into $-V_1$ to obtain the converse
inequality.  Then, exchanging the roles of indices $1$ and $2$, we have, for $i=1,2$,
$$
\int_0^T \!\!< V_i(s,\cdot),   \dot {\phi_i} (s,\cdot)>ds
= {p_i}  \int_0^T\!\!\!\!
<\tanh(\phi(s,\cdot)*J+a_i\theta)V_i(s,\cdot)>ds
-\int_0^T\!\!\!\!<\phi_i(s,\cdot)V_i(s,\cdot)>ds.
$$
This   means that $\mus\phi$ is the (unique) weak solution
of \eqv(G2b), since by definition \eqv (2.4)    of the rate functional   the initial condition is fulfilled.

For the  reverse, we  prove that if    $\mus \phi  \in
\ac ([0,T], B_{1,1})$  is the solution  of   equation  \eqv (G2b), then  $ \mus  \pi = (\phi_1 \l, \phi_2\l)$ is such that $J_1(\mus \pi)=0$; hence, by  Proposition \eqv (3A1), $I_0(\mus \pi)=0$.  We  insert     equation  \eqv (G2b)
into the explicit representation  \eqv (DD1c).   
Namely if $ \mus  \pi $ solves  \eqv (G2b)  then,  by Corollary \eqv (26J), 
$$ \log \frac { D_i } { p_i- \phi_i  }-  A_i =0, \Eq  (L33)$$  
$$      R_i =
 D_i- (\partial_t \phi_i )    \cosh A_i  
= e^{A_i}  (p_i- \phi_i) +( \phi_i  -p_i \tanh  A_i )   \cosh   A_i.
  $$
  Hence   
  $$   p_i-\phi_i \tanh  A_i   - \frac { R_i}{\cosh  A_i }  =
     (p_i - \phi_i)(1+ \tanh  A_i -\frac {e^{A_i}}{\cosh  A_i })=0.\Eq (L2a) $$
 By \eqv(L33), \eqv(L2a), the right hand side of \eqv (DD1c) is equal to zero, which completes the proof of (3).
\qed

\bigskip
\chap{9. Appendix B}9
\numsec= 9
\numfor= 1
\numtheo=1

This appendix is devoted to proofs postponed from Sections 5, 6 and 7. 
\medskip \noindent {\bf Proof of Proposition \equ(RN1)}
Let $s\in[0,T]$, 
$${C(V_1,V_2)=\sum_{i=1,2}C(V_i)=\sum_{i=1,2}\sup_{s\in [0,T]} \sup_{r\in\L} \left(\sinh^2[ V_i(s,r)] + |\sinh[ 2V_i(s, r)] |\right)}.$$ 
Then  
$$ 
\left|
F_{\mus V(s)}(\mus\l^\g(\a),\mus\pi^\g_s)
-F_{\mus V(s)}( (p_1\l^\g, p_2\l^\g), \mus\pi^\g_s)
\right|  
\le   
 \left | 
 \frac {\g^d}2 \sum_{x\in\L_\g} \left [ \a_i (x)-p_i \right ] \BB_i(x,\s,s) \right |, 
  $$ 
 with
$$  \BB_i(x,\s,s)=  \cosh[2 V_i(s,\g x)]-1 +\tanh[ ( J_\g \star \s_s)(x)+a_i\th]\sinh[ 2V_i(s,\g x)] $$  
  $$     | \BB_i( x,\s,s) |   \le
  \sinh^2[  V_i(s,\g x)]   + |\sinh[ 2V_i(s,\g x)] |   \le  C (V_i).   
$$
Take $ l \in \Z $, $ l \neq 0$. Since $\E[\a_i(x)] =p_i$ for all $x \in \L_\g$, 
$$ \eqalign {   \frac {\g^d}2\sum_{x\in\L_\g}  \left [ \a_i (x)-p_i \right ]    \BB_i( x,\s,s)  &=
\frac {\g^d}2\sum_{x\in\L_\g} \frac 1 {(2l+1)^d} \sum_{|y| \le l}   \left [   \BB_i( x+y,\s,s)-   \BB_i( x,\s,s) \right ]  \a_i(x+y)  \cr &  -  
\frac {\g^d}2\sum_{x\in\L_\g}    \BB_i( x,\s,s)   \left [        \a^{(l)}_i(x) -\E[\a_i(x)]  \right ]. } \Eq (PA.5)
$$
Using  uniform continuity  as in the proof of Lemma
\eqv(unifbound), there exists a positive
function $\e$ on $\R_+$ with $\lim_{s\to 0}\e(s)=0$ (depending only on $T$,
$J$ and $\mus V$)   such that the first term on the r.h.s. of \eqv(PA.5) is bounded uniformly in $\a$ and $ \s$; for the second term, let   $ \d>0$ and  $\EE_i(\d, l,\g,\a)$  defined in \eqv (ju1). We  conclude  by
  $$ 
\left |  \frac {\g^d}2\sum_{x\in\L_\g} \left [ \a_i(x)-p_i \right ]   \BB_i( x,\s,s)  \right | \le \e (\g l) 
+     C(V_i) \left[  \d    +   \EE_i(\d, l,\g,\a)\right ].  $$  
\qed

\medskip \noindent {\bf Proof of Proposition \equ(exptight)}
    Consider a sequence of functions $ \{H_{k}\}_{k\ge 1}$ in $\CC^2(\Lambda)$ dense in  $\CC(\Lambda)$ 
for the uniform topology with  $ \|H_{k}\|_\infty\le 1$.   Denote  for all integers $m \ge 1 $, $\ell\ge 1$,  and $\delta>0$, 
$$ \AA_{m, \delta, \ell}= \{ \mus \pi \in  D ([0,T], \MM_1 \times \MM_1): \inf_{\{t'_i\}} \max_i \sup_{t'_i \le s <t'_{i+1} }  \sum_{j=1}^{m} \frac { | <\mus \pi_t, H_j>- <\mus \pi_s, H_j>|}{ 4^j} \le \frac  {\ell+1} m\},
$$
 where the infimum extends over all positive integers $K$ and all finite sets of points $ \{t'_i,  0 \le i \le K\}$ satisfying $0= t'_0<t'_1 < \dots <t'_K=T $,  $t'_{i+1} - t'_{i}  >  \delta$.
 We first show 
 that for $m\ge 1$ and for $\ell\ge1$ there exists $ \delta (m, \ell) $  and $ \g_0(m, \ell)$  so that for all $\g \le \gamma_0 $,
 $$   \mus  Q^{\g, \a}_{\s^\g} [ \mus \pi \notin  \AA_{m, \delta(m, \ell), \ell} ]  \le    e^{- \frac   { \ell +1} {\g^d} }.  
$$
This is done taking into account that 
 $$ \eqalign { &  \{ \inf_{\{t'_i\}} \max_i \sup_{t'_i \le s <t'_{i+1} }| <\mus \pi_t, H >- <\mus \pi_s, H>|\ge \frac   {\ell+1}  m \} \cr &  \qquad \qquad \subset \cup_{ k=0}^{\frac T \delta}  \{ \sup_{k\delta \le t<(k+1) \delta} | <\mus \pi_t, H >- <\mus \pi_{k\delta}, H >|\ge \frac   {\ell+1}   {4m} \}, } $$
and estimating the right hand side as in [KL] (p. 271, after formula (4.6)). 
Then the construction of the compact $\mus K_\ell$ is obtained by a general
procedure as explained in [Bill] and Section 8 of [QRV].    \qed

\medskip

   Next lemma states some technical results needed in the extension of the lower bound of the large deviation principle. 
  \medskip   
 
\noindent {\bf {\Lemma   (3A4)}}  {\it  
For $ (t,t') \in  [0,T]^2$, $ (x,y) \in \L^2$, $ \d_0>0$,  let $ \tau\in \{t, t'\}$, $ \zeta \in \{x, y\}$. We assume      
$(\mus u (\tau,\cdot), \mus v(\tau,\cdot))\in(B_{p_1-\d_0,p_2-\d_0})^2$,   $g_i(\tau, \zeta) \in \R$ and $h_i(\tau, \zeta) \in \R$, for   $i=1,2$;  we have 
$$| H_i (\mus u (t,x),  g_i(t',y)) - H_i (\mus v (t,x),   g_i(t',y)) |  \le 
   |g_i(t',y)|  \left(1 +  \frac {K} {\d_0 }\right) \left(  \|  v(t,\cdot)-    u(t,\cdot) \|_1    +     |u_i(t,x)-  v_i(t,x)| \right),   \Eq (T.1)
  $$
  $$| H_i(\mus v (t,x),  g_i(t',y)) - H_i (\mus v (t',y),   g_i(t',y)) |  \le 
     |g_i(t',y)|  \left(1 +  \frac {K} {\d_0 }\right) \left(  |x-y|   +     |v_i(t,x)-  v_i(t',y)| \right),  \Eq (T.2)
  $$
   $$| H_i(\mus u (t,x),  g_i(t,x)) - H_i (\mus u (t,x),   h_i(t,x)) |   \le \left(K+\log\frac 1{\d_0} +K |g_i(t,x)|\right) | g_i(t,x)- h_i(t,x)|,
   \Eq (T.2a)
  $$
where the constant $K=K(J,\th)$   may change from one  occurrence to the other.   
} 

\smallskip
\proof  
The assumptions   
  enable  to prove \eqv (T.1)--\eqv(T.2a) by writing  formula \eqv (DD1c) for $H_i$, using \eqv (DD1a) for $A_i,D_i,R_i$. The latter depend on $\mus u$ or $\mus v$, 
$g_i$ or $h_i$, $(t,x)$ or $(t',y)$. In each computation, we stress the dependence on the involved quantities, writing e.g. $A_i(u)$ for \eqv (T.1), $R_i(g_i)$ for \eqv (T.2a). Notice that, unlike in \eqv (DD1c), those functions depend not only on $(t,x)$, but on $(t,x)$ and $(t',y)$; this does not change the expression of $H_i$, since \eqv (DD1c)  was established pointwise in the proof of Proposition \eqv (3A22). In the intermediate computations, we omit to write $(t,x),(t',y)$.

We begin with auxiliary estimates. For \eqv (T.2a), notice that $ |(J * u)(t,x)|\le 1$ (since $\int J(r)\,dr=1$), and $$|A_i(u)|\le K(J,\th).\Eq (E0)$$   When  $|v_i | \le p_i- \d_0 $ and $g_i>0$, 
we have   $p_i^2- v^2_i=(p_i- |v_i|)(p_i+|v_i|)\ge p_i\d_0$, hence
$$  g_i K(J,\th)+1 \ge  R_i (\mus u,g_i) \ge \sqrt{g_i^2+p_i\d_0} \ge \max(g_i,\sqrt{p_i\d_0}),  
 \Eq (E2)$$
 $$ \eqalign{  
   g_i K(J,\th)+1 &\ge    D_i (\mus v,g_i)\ge   g_i+  \max(g_i,\sqrt{p_i\d_0}). }\Eq (BB.10) $$
 For \eqv (T.1), we need  
 $$ |(J * u)(t,x)-(J* v)(t,x)|  \le  \|J\|_\infty  \|u(t, \cdot)- v (t, \cdot)\|_1  \Eq (T.3)$$ 
 and its consequences
$$ \eqalign { & 
 |v_i (t,x)\tanh [A_i(v)(t,x)]   - u_i(t,x) \tanh [A_i(u)(t,x)] |
 \cr & \le |v_i(t,x)- u_i(t,x)| |\tanh  [A_i(v)(t,x)]  |+
 |u_i(t,x)| |\tanh [A_i(v)(t,x)]   -\tanh  [A_i(u)(t,x)] |
  \cr & \le  |v_i(t,x)- u_i(t,x)|+   \|J\|_\infty\|u(t, \cdot)- v (t, \cdot)\|_1. 
 }\Eq (E1)$$
  $$      \left |   \cosh [A_i(u)] -   \cosh  [A_i(v)]   \right |    \le   K(J,\th)  \|  u (t,\cdot)- v (t,\cdot) \|_1.    \Eq (BB.11a)  $$
 Respectively for \eqv (T.2), we need,     since $| v(\cdot,\cdot)| \le1$,  
 $$ \eqalign { & |(J * v)(t,x)-(J * v)(t',y)|  \le
|(J * v)(t,x)-(J * v)(t,y)|+ |(J * v)(t,y)-(J *  v)(t',y)|\cr &   \le \|J' \|_\infty  |x-y|
+ | \int_\L J(y-z)    [v(t',z) - v(t,z)]\,dz |
\cr & \le \|J' \|_\infty  \left(|x-y|+ \|v(t',\cdot) - v(t,\cdot)\|_1\right). 
 }\Eq (ET.4) $$
 as well as its consequence analogous to \eqv (E1).
   The proofs of \eqv (T.1), \eqv (T.2) go along the same scheme. Namely 
all the estimates are done pointwise, they rely respectively on  \eqv (BB.10) to \eqv (ET.4), the other changes being straightforward starting from expressions analogous to  \eqv(EJ2) below. Hence we detail only the proof of \eqv (T.1).  We have 
   $$ \eqalign { &  2\left[   H_i (\mus u (t,x),   g_i(t',y)) - H_i (\mus v (t,x),   g_i(t',y)) \right] \cr &=
    g_i(t',y) [(J *  v)(t,x)- (J * u)(t,x)] +  g_i(t',y)  \left [  \log \frac {
D_i (\mus u) (t,x) } { p_i- u_i(t,x) } - \log \frac {
D_i (\mus v) (t,x) } { p_i- v_i (t,x)} \right ] \cr & 
+ v_i(t,x) \tanh  [A_i(v)(t,x)]  - u_i(t,x) \tanh   [A_i(u)(t,x)]    
 +    \frac { R_i(\mus v)(t,x)} {\cosh  [A_i(v)(t,x)]   } -     \frac { R_i (\mus u)(t,x)} {\cosh   [A_i(u)(t,x)]  }.   }   \Eq (EJ2)
$$

Next we  show
$$\left |  g_i(t',y) \left [  \log \frac {
D_i (\mus u)(t,x)  } { p_i- u_i(t,x) } -  \log \frac {
D_i (\mus v)(t,x)  } { p_i- v_i(t,x) } \right ] \right |  \le  \frac K {\d_0} | u_i(t,x)-v_i(t,x)| (1+ |g_i(t',y)|) + |g_i(t',y)| \|u(t,.)-v(t,.)\|_1.      \Eq (PP.3) $$ 
To this aim,  see  \eqv (EJ1)--\eqv (DD2ZZ), it is enough  to estimate,
when $g_i(t',y)>0$ and uniformly for $\theta \in \R$,
$$\eqalign { & \left| F (\mus u(t,x), g_i(t',y),\th) -F (\mus v(t,x), g_i(t',y),\th)\right|=   \left|g_i(t',y) \left [  \log \frac {
D_i (\mus u)(t,x)  } { p_i- u_i(t,x) } -  \log \frac {
D_i (\mus v)(t,x)  } { p_i- v_i(t,x) } \right ]\right|  \cr &\le g_i\left [\left|\log \frac { p_i- v_i } { p_i- u_i }\right| +
 \left| \log \frac {D_i (\mus u)} {D_i (\mus v)  }\right| \right ] \cr & \le
 g_i \left|\frac {u_i- v_i } { p_i- u_i } \right|+
 g_i \left|\frac {D_i (\mus u)-D_i (\mus v)} {D_i (\mus v)  }\right| 
 },  \Eq (BB.1)$$
because  
$|\log (1+a)|\le\log (1+|a|) \le |a|$.  By \eqv (BB.10),   $  g_i \le   D_i (\mus v) $. 
Using also \eqv (BB.11a) we get 
   $$  \eqalign {  
\left| g_i \frac{D_i(\mus u)-D_i(\mus v) } {D_i(\mus v)} \right|&\le 
 \left |  \frac {g_i^2 \left(\cosh[A_i(u)]-\cosh  [A_i(v)]\right )}                          
  {D_i (\mus v)} \right|    +
\left |g_i \frac{R_i(\mus u) -R_i (\mus v)}{ D_i(\mus v) } 
  \right |\cr   &
  \le    g_i  K(J,\th)   \|u (t,.)- v (t,.) \|_1  +
\left |  R_i (\mus u) - R_i (\mus v)    
  \right |. }   \Eq (BB.11)  $$
    To estimate the  second term on the right hand side of   \eqv (BB.11), we  apply   \eqv (E0), \eqv (E2), \eqv (BB.11a) and obtain 
    $$ \eqalign { & \left| R_i (\mus u)- R_i (\mus v)   \right |  =
  \left |\frac {\left [R_i (\mus u) \right ]^2  -   \left [R_i (\mus v) \right  ]^2 }{R_i (\mus u) +  R_i (\mus v) }   \right |  \cr &
   \le 
   \left | \frac {g_i^2 \left(\cosh^2 [A_i(u)]- \cosh^2 [A_i(v)]\right)} { R_i (\mus u) +  R_i (\mus v)} \right| +
   \left|   \frac {v^2_i - u^2_i} {R_i (\mus u)  + R_i (\mus v)} \right |  
       \cr & \le
 \frac {g_i^2K(J,\th)   \| u(t,\cdot)-  v(t,\cdot) \|_1} 
{R_i (\mus u)+ R_i (\mus v)} +
\left | \frac{v^2_i - u^2_i}{R_i(\mus u)  + R_i (\mus v)}\right | \cr &\le 
  |g_i| K(J,\th)   \| u(t,\cdot)-  v(t,\cdot) \|_1  +  \frac  { p_i }  {\sqrt{p_i\d_0} }   |u_i-  v_i|. 
   } \Eq (BB.11b)  $$ 
Combining \eqv (BB.1), \eqv (BB.11b) we obtain \eqv (PP.3). 
Next we  estimate the last term of  \eqv (EJ2).  Taking into  account  \eqv (BB.11a) and  \eqv (BB.11b) we  have  
  $$ \eqalign { &   \left |   \frac { R_i(\mus v)}{\cosh [A_i(v)]} -  \frac { R_i(\mus u)} {\cosh [A_i(u)]}\right |  =  \left|  \frac { R_i(\mus v)- R_i(\mus u)} {\cosh[A_i(v)]}  +   R_i (\mus u)  \frac {\cosh [A_i(u)]- \cosh [A_i(v)]} {\cosh [A_i(v)]\cosh [A_i(u)] }   \right | \cr
   &   \le  
    \left |   R_i (\mus v)- R_i (\mus u) \right| +  R_i (\mus u) \left| \cosh [A_i(u)]- \cosh [A_i(v)]  \right| 
  \cr &  
   \le  \sqrt{\frac {p_i}{\d_0} }   |u_i-  v_i|
    +   \left(K(J,\th)  |g_i|+ K(J, \th) [|g_i| +1]\right)
    \|  v(t,\cdot)-   u(t,\cdot) \|_1.
   }  \Eq (BB.6) $$
Finally, combining \eqv(T.3), \eqv(E1), \eqv(PP.3), \eqv(BB.6) 
yields \eqv (T.1).   

We now derive \eqv (T.2a) in a similar way: 
    $$ \eqalign { &   2\left[  H_i (\mus u (t,x),   g_i(t,x)) - H_i (\mus u (t,x),   h_i(t,x)) \right] \cr &=
   [ h_i(t,x)-   g_i(t,x)] \left(A_i(u)(t,x)+\log\frac 1{ p_i- u_i (t,x)}\right) 
    +  g_i(t,x) \log  {D_i (\mus u,g_i) (t,x)} -  h_i(t,x) \log {D_i (\mus u,h_i) (t,x)}   \cr & 
 + \frac {R_i(\mus u, h_i)(t,x)- R_i (\mus u, g_i)(t,x) }{\cosh  [A_i(u)(t,x)] }.    }   
$$
We have, restricting ourselves to $g_i>0,h_i>0$,   see  \eqv (EJ1)--\eqv (DD2ZZ),  and using first
$|\log(1+a)| \le |a|$ as in \eqv(BB.1), then \eqv(BB.10),
$$  \eqalign{ \left|  g_i  \log  {
D_i ( g_i)  }  -  h_i   \log  {
D_i ( h_i) } \right| &\le  |g_i  - h_i | \left| \log  {
D_i (g_i)}\right| +  h_i \left|\log {D_i ( g_i) }- \log {D_i ( h_i) } \right| \cr &
\le |g_i  - h_i |\left|  \log  {
D_i (g_i)}\right| + 
 h_i \left|\frac {D_i ( g_i) - D_i ( h_i)  }  {D_i ( h_i)}  \right|
\cr &
\le   |g_i  - h_i |  \left(\left|\log  {
D_i (g_i)}\right|+\frac{h_i\cosh  [A_i(u) ]}{D_i ( h_i)}\right) + h_i\left|\frac{R_i( g_i)  - R_i( h_i)}{D_i ( h_i)}\right|\cr &
\le   |g_i  - h_i |  ( g_i+1) K(J,\th) + | R_i( g_i)  - R_i( h_i)|.
}$$
Then, as in \eqv(BB.11b),
$$ R_i(g_i)- R_i( h_i)  = \frac { R_i^2( g_i)- R_i^2( h_i)  } {R_i(g_i)+ R_i( h_i) }  =
\frac {  \left ( g_i^2 -  h_i^2 \right ) \left(\cosh [A_i(u) ]\right)^2  }   {R_i(g_i)+ R_i(h_i)  }.  
  $$ 
Therefore, using that $|u_i|\le p_i-\d_0$, and \eqv(E2),
$$ \eqalign { &  2\left|   H_i (\mus u,   g_i) - H_i (\mus u , h_i)\right|   \le |h_i- g_i |\left(K(J,\th) +  \log \frac 1{\d_0}+( g_i+1) K(J,\th)+2[K(J,\th)]^2\right).
   }  
   $$
\qed

   \smallskip
   \noindent {\bf  Proof of Lemma  \eqv (J40)} 
 We   exploit that  $\mus \psi \in \ac([0,T], B_{p_1,p_2})$ and $\mus\psi$ is differentiable in time in $(T,T+1]$ (see \eqv(S.2z)), 
hence $\displaystyle{\frac {\partial    \psi_i } {\partial t}  \in L^1([0,T+1] \times \L)}$. Therefore for $A >0$   and
  $$D_A = \{  x \in \L :  \sup_{t \in [0,T+1]}   \sum_{i=1}^2  |\frac {\partial    \psi_i } {\partial t}  (t,x)| >A \} . 
$$
we have for all $s\in[0,1]$,
      $$    \lim_{A \to \infty}   \sum_{i=1}^2 \int_0^T  \int_{\L}  |\frac {\partial    \psi_i } {\partial t}  (t+s,x)| \1_{D_A} (x)\,dxdt =0\, . 
$$  
    By   \eqv (T.1) of  Lemma  \eqv (3A4),  we obtain, for $2\nu_0=  \min \{\d_1, \d_2\} $,  splitting $ \L = D_A \cup D^c_A $, 
 $$  \eqalign {   |W_1| &    
   \le    \int_{\R} \Psi_{\e_0} (s)       \int_{\L}   \Phi_
   {\e_1} (y)   \int_{2\eta}^T \int_{\L}
   \sum_{i=1}^2  |\frac {\partial    \psi_i } {\partial t}  (t+s
   ,x-y)|(\1_{D_A} (x-y)+\1_{D_A^c} (x-y)) ( 1+ \frac K {\nu_0})  \cr
&\qquad\left \{ |  \psi^\e_{i}    (t,x)
 -  \psi_{i}      (t,x)| +    \|  \psi^\e_{i}  (t, \cdot )
 -  \psi_{i}      (t,\cdot )\|_1     \right \} \,dxdtdyds
   \cr &     \le
 A     ( 1+ \frac K {\nu_0})  \sup_{t \in [0,T] } \sum_{i=1}^2 \left ( \int_{\L} | \psi^\e_{i}    (t,x) 
 -  \psi_{i}      (t,x)|\,dx  +  \|  \psi^\e_{i}   (t, \cdot )
 -  \psi_{i}      (t,\cdot )\|_1   \right ) \cr &   +  4 
    \sum_{i=1}^2 \int_{2\eta}^T \int_{\L}|\frac {\partial    \psi_i } {\partial t}  (t+s,x)| \1_{D_A}(x)\,dxdt, } 
 $$ 
 where we noticed that $\left \{ |  \psi^\e_{i}    (t,x)
 -  \psi_{i}      (t,x)| +    \|  \psi^\e_{i}  (t, \cdot )
 -  \psi_{i}      (t,\cdot )\|_1     \right \}\le 4$.
 Letting first $ \e \to 0$ then $ A \to \infty $,  we get 
 $     \lim_{\e \to 0}  |W_1| =0.$ 
Next we estimate $W_2$. We apply   \eqv (T.2)  of   Lemma \eqv  (3A4).  
 More precisely
 $$\eqalign{
 &\left|
 \HH (\mus \psi (t+s,x-y) , \frac {\partial   \mus \psi } {\partial t}(t+s,x-y)) -\HH(\mus \psi     (t,x),  \frac {\partial   \mus \psi } {\partial t}  (t+s ,x-y)) \right|\cr &
  \le  \sum_{i=1}^2 \left|\frac {\partial    \psi_i } {\partial t}(t+s,x-y)\right|
\left(1 +  \frac {K} {\nu_0 }\right)    
 \left\{
  \left( 
   |{\psi_i}(t+s,x-y)- {\psi_i}(t,x) |\right)
 + |y|    \right\}.   }$$
  As  before take $A>0$ large enough,     split $ \L = D_A \cup D^c_A $, to get
  $$\eqalign{ |W_2| \le 
 &  A \left(1 +  \frac {K} {\nu_0 }\right)  \sum_{i=1}^2 \int_{\R}\Psi_{\e_0} (s)      \int_{2\eta}^T\int_{\L}   \Phi_ {\e_1} (y)\int_{\L}     
 \left\{
  \left(     |{\psi_i}(t+s,x-y)- {\psi_i}(t,x) |\right)
 + |y|    \right\}\,dxdydtds \cr & + C
 \sum_{i=1}^2  \int_{2\eta}^T \int_{\L}  \int_{\R}   \Psi_{\e_0} (s)         \1_{D_A} (x)\left|\frac {\partial     \psi_i } {\partial t}(t+s,x)\right| \,dsdxdt. 
  }$$
  Since
  $$\lim_{\e_1 \to 0 } \int_{\L}   \Phi_ {\e_1} (y) |y| dy =0\, ,\quad
  \lim_{\e_0 \to 0}  \int_{\R}   x\Psi_{\e_0} (s)    \int_{\L}  
    |{\psi_i}(t+s,x-y)- {\psi_i}(t,x) |\,dxds   =0,  $$
  letting   $\e \to 0$  and  then $A\to \infty$ we  obtain 
  $ \lim_{\e \to 0} |W_2| =0$. 
\qed 
   \bigskip
   \noindent {\bf  Proof of Lemma  \eqv (EA10)} 
We have 
 $$ \eqalign { & \int_{ \eta }^{2\eta} \int_{\L} \left[\HH (\mus \psi^\e,  \frac {\partial   \mus \psi^\e} {\partial t}) (t,x)  -    \HH (\mus \psi,  \frac {\partial   \mus \psi} {\partial t}) (t,x)\right]\,dxdt\cr &  =
 \int_{ \eta }^{2\eta}\int_{\L}\left[ \HH (\mus \psi^\e,  \frac {\partial   \mus \psi^\e} {\partial t}) (t,x)  - 
  \HH (\mus \psi,  \frac {\partial   \mus \psi^\e} {\partial t})
 (t,x)\right]\,dxdt 
\cr & +  \int_{ \eta }^{2\eta}\int_{\L}\left[ \HH (\mus \psi,  \frac {\partial   \mus \psi^\e} {\partial t}) (t,x)- 
  \HH (\mus \psi,  \frac {\partial   \mus \psi} {\partial t}) (t,x)\right]\,dxdt. }  \Eq (EA2)$$  
 The first term is estimated by applying \eqv (T.1) of Lemma \eqv(3A4). 
We have   
 $$ \eqalign { &  \left | \int_{ \eta }^{2\eta}\int_{\L} \left[\HH (\mus \psi^\e,  \frac {\partial   \mus \psi^\e} {\partial t}) (t,x)  - 
  \HH (\mus \psi,  \frac {\partial   \mus \psi^\e} {\partial t}) (t,x)\right]\,dxdt\right | \cr  & \le \sum_{i=1,2}  \int_{ \eta }^{2\eta}\int_{\L}   \left(1 +  \frac {K} {\nu_0 }\right)  | \frac {\partial   \psi^\e_i(t,x)}  {\partial t}| 
 \left(  \|  \mus  \psi^\e (t,\cdot)-    \mus  \psi (t,\cdot) \|_1         + 
|\psi_i(t,x)-  \psi^\e_i(t,x)| \right)\,dxdt,  }   \Eq  (EA1)
 $$
 where  $\nu_0=  \min \{\d_1, \d_2\} $. 
  Note that   $$\frac {\partial     \psi^{\e}_i  } {\partial t}   (t,\cdot)=
\chi_1 (t) \frac {\partial     \psi_i  }  {\partial t}  (t, \cdot) +
  \chi_2 (t)
\int_{ \R } \Psi_{\e_0}(s)   ( \Phi_{\e_1} *(\theta_s \frac {\partial     \psi_i  } {\partial t}) (t, \cdot)  +   \chi_1'(t)   [\psi _i(t,\cdot)- \psi^{\e}_i(t,\cdot) ]\,ds,     \quad t \in [0,T], 
$$ 
where  we denote $\chi_i'(t)= \frac {d} {{d}t} \chi_i(t)$, for $ i=1,2$, and we use that    $\chi_2'(t)= -   \chi_1'(t)$.  
 Since  for $t \in (0, 3 \eta)$,  $ \mus \psi (t)= \mus R(t)$  solves  \eqv (G2b)  we have
 $$  \sup_{t \in [\eta, 2 \eta]} \sup_{x \in \L} | \frac {\partial     \psi^{\e}_i  } {\partial t}   (t,\cdot)| \le  \frac 6 \eta.  \Eq (par2)$$
  Therefore 
for all  $\eta$ letting  $ \e \to 0$  the term in the right hand side of \eqv (EA1) goes to zero.   
For the second  term in the r.h.s. of \eqv (EA2), applying  \eqv (T.2a)   
in  Lemma \eqv (3A4) and taking into account  \eqv (par2), we get the result.  
   \qed

 \bigskip
\noindent {\bf Acknowledgements}. We thank Errico Presutti for helpful discussions. E. O. acknowledges the warm hospitality of  Universit\'e de Rouen and Universit\'e Paris Descartes.  Part of this work was done during the authors' stay at Institut Henri Poincar\'e, Centre
\'Emile Borel (whose hospitality is acknowledged), for the semester
``Interacting Particle Systems, Statistical Mechanics and
Probability Theory''.


   \bigskip
 \goodbreak
 \noindent 
 {\bf References}
  
  \item{[BBI]}   Bianchi, A.,  Bovier,  A., Ioffe, D.  Sharp Asymptotics for metastability in the random field Curie-Weiss model mean-field models,  {\sl Electron. J. Probab.}   {\bf 14}(2009), 1541--1603.

 \item{[Bill]} 
    Billingsley,  P.  {\sl Convergence of probability measures.} Wiley and Sons, New York, 1968. 

 \item{[Bo]}   Bovier, A.  {\sl Statistical Mechanics of Disordered Systems. A Mathematical Perspective.} Cambridge University Press,  2006.
 
 \item{[BEGK]}    Bovier,  A.,  Eckhoff, M.,   Gayrard, V.,  Klein, M. Metastability in stochastic dynamics of disordered mean-field models,  {\sl Probab. Theory Relat. Fields}   {\bf 119}(2001), 99--161.

\item{[CDS]} Collet, F., Dai Pra, P., Sartori, E.    A simple mean field model for social interactions: dynamics, fluctuations, criticality. {\sl J. Stat. Phys.} {\bf 139}(2010), 820-858.

  \item{[COP]}  Cassandro, M., Orlandi, E., Picco, P.
 Typical configurations for one dimensional 
 random field  Kac model.
{\sl Ann. Probab.}, {\bf 27}, no. 3(1999),  1414--1467.    

 \item{[COPV]}  Cassandro, M., Orlandi, E., Picco, P., Vares,  M. E. 
  Typical configurations for one dimensional 
 random field  Kac model: localization of the phases. 
 {\sl Electron. J. Probab.},  {\bf 20}(2005),  786--864. 
         
  \item{[COP4]}  Cassandro, M., Orlandi, E., Picco, P.  The optimal interface profile for a non-local model of phase separation. {\sl   Nonlinearity} {\bf 15}(2002), 1621--1651.
 
 \item{[C]} Comets, F.  Nucleation for a long range magnetic model. {\sl  Ann. Inst.
            Henri Poincar\'e, Probab. Statist.} {\bf 23}, no. 2(1987), 135--178.
 
  \item{[CE]} Comets, F., Eisele, T.  Asymptotic dynamics, non-critical and critical fluctuations for a geometric long-range interacting model.
  {\sl Commun. Math. Phys.} {\bf 118}(1988), 531--567.

 \item{[DD]} Dai Pra, P., den Hollander, F.    McKean-Vlasov limit for interacting random processes in random media.
  {\sl J. Stat. Phys.} {\bf 84}(1996), 735--771.
   
 \item{[DOPT]} De Masi, A., Orlandi, E., Presutti, E. and Triolo, L.   Glauber evolution with Kac potentials I. Mesoscopic and macroscopic limits, interface dynamics. {\sl Nonlinearity}, {\bf 7}(1996), 287--301.
 \item{   }  Glauber evolution with Kac potentials II. Fluctuation. {\sl Nonlinearity}, {\bf 9}(1996), 27--51.
 \item{   } Glauber evolution with Kac potentials. III. Spinodal decomposition. {\sl  Nonlinearity}, {\bf 9}(1996), 53--114.
 
 \item{[DV]} Donsker, M.D., Varadhan S.R.S.    Large deviations from an hydrodynamical limit. {\sl Comm. Pure Appl. Math.} {\bf 42}(1989), 243--270.

  \item{[DS]} Dunford, N., Schwartz,  J.T.  {\sl  Linear operators.} Wiley and Sons, New York, 1963. 
     
    \item {[FLM]}
 Farfan, J.S., Landim, C. and Mourragui, M. 
 Hydrostatics and dynamical large deviations of boundary driven gradient symmetric exclusion processes. 
 {\sl Stoch. Proc. and their Appl.} {\bf 121}(2011),  725--758.
 
    \item{[FMP]} Fontes, L.R., Mathieu, P., Picco, P.  On the averaged dynamics of the random field Curie-Weiss model. 
  {\sl  Ann. Appl. Probab.} {\bf 10}(2000), 1212--1245.
  
  \item{[HS]} Holley,  R.A, Stroock, D.W.  A martingale approach to infinite systems of interacting processes. {\sl Ann. Probab.}, {\bf 4},   vol. 2(1976), 196--223.

 \item{[K]} Koukkous, A. Hydrodynamic behavior of symmetric zero-range processes with  random rates. {\sl Stoch. Proc. Appl.} {\bf 84}(1999), 297--312.
 
 \item{[KL]}  Kipnis, C., Landim, C. {\sl Hydrodynamic limit of interacting particle systems}.  Springer-Verlag, 1999.
 
  \item{[KUH]} Kac, M., Uhlenbeck,  G., and Hemmer, P.C.
On the van der Waals theory of vapour-liquid equilibrium. 
I. Discussion of a
one-dimensional model. {\sl J.  Math. Phys.} {\bf 4}(1963), 216--228.

  \item{[LP]} Lebowitz J.  and Penrose O.  Rigorous treatment of
the Van der Waals Maxwell theory of the liquid-vapour
transition. {\sl J. Math. Phys.} {\bf 7}(1966), 98--113.

   \item{[MP]}   Mathieu, P., Picco, P.    Metastability and convergence to equilibrium for the random field Curie-Weiss model,   {\sl J. Stat. Phys.}   
   {\bf 91}(1998), 679--732.
   
    \item{[MO]}   Mourragui, M., Orlandi, E. 
Large deviations from a macroscopic scaling limit for particle systems with   Kac interactions  and random potential. 
{\sl  Ann. Inst.
            Henri Poincar\'e, Probab. Statist.} {\bf 43}(2007), 677--715.
 
  \item{[OP]}  Orlandi, E., Picco, P.  
 Weak large deviation principle for one dimensional 
 random field  Kac model. {\sl Electron. J. Probab.},  {\bf 14}(2009),  1372--1416. 
 
 \item{[QRV]}  Quastel,  J., Rezakhanlou, F. and Varadhan, S.R.S.  Large deviations for the symmetric simple
exclusion process in dimension $d\ge 3$. {\sl   Probab. Theory Relat. Fields},  {\bf 113}(1999),  1--84. 
 
 \item{[P]}  Presutti, E. {\sl   Scaling Limits in Statistical Mechanics and Microstructures in Continuum Mechanics.}
Springer Berlin Heidelberg, Series: Theoretical and Mathematical Physics,
2009. 

\item{[Sp]} Spohn, H. {\sl Large scale dynamics of interacting particles}.
Springer, Berlin, 1991.

  \item{[V]}  Varadhan, S.R.S. {\sl Large deviations and Applications}  CBMS-NSF regional conference Series in Applied Mathematics, {\bf 46}(1984), SIAM,  Philadelphia.   
  
\bye